\pdfoutput=1
\documentclass[11pt]{amsart}

\usepackage[margin=1in]{geometry}
\usepackage{amsmath,amssymb,amsthm}
\usepackage{mathrsfs}
\usepackage{xcolor}
\usepackage{booktabs,tabularx,longtable,array,xltabular}
\usepackage{enumitem}
\usepackage{needspace}
\usepackage{tikz}
\usetikzlibrary{arrows.meta,positioning,backgrounds}
\numberwithin{table}{section}
\numberwithin{figure}{section}
\numberwithin{equation}{section}
\usepackage[colorlinks=true,linkcolor=blue,citecolor=blue,urlcolor=blue]{hyperref}

\hypersetup{pdftitle={Centres and semisimple deformations of contraction algebras}, pdfauthor={Joseph Karmazyn and Shengyong Pan}}

\newtheorem{thm}{Theorem}[section]
\newtheorem{lem}[thm]{Lemma}
\newtheorem{prop}[thm]{Proposition}
\newtheorem{cor}[thm]{Corollary}
\theoremstyle{definition}
\newtheorem{ex}[thm]{Example}
\newtheorem{rem}[thm]{Remark}
\theoremstyle{plain}
\newtheorem*{thmA}{Theorem A}
\newtheorem*{thmB}{Theorem B}
\newtheorem*{thmC}{Theorem C}
\newtheorem*{propD}{Proposition D}
\newtheorem*{propE}{Proposition E}
\newtheorem*{conjJS}{Conjecture}

\newcommand{\Acon}{A_{\mathrm{con}}}
\newcommand{\cC}{\mathcal{C}}
\newcommand{\cA}{\mathcal{A}}
\newcommand{\cO}{\mathcal{O}}
\newcommand{\cX}{\mathcal{X}}
\newcommand{\cY}{\mathcal{Y}}
\newcommand{\h}{\mathfrak{h}}
\newcommand{\bul}{\bullet}
\newcommand{\tGamma}{\widetilde{\Gamma}}
\newcommand{\SGV}{S_{\mathrm{GV}}}
\newcommand{\HH}{\mathrm{HH}}
\newcommand{\hHH}{\widehat{\mathrm{HH}}}
\newcommand{\Z}{\mathbb{Z}}
\newcommand{\Q}{\mathbb{Q}}
\newcommand{\C}{\mathbb{C}}
\newcommand{\A}{\mathbb{A}}
\newcommand{\sAlg}{\mathfrak{A}}
\DeclareMathOperator{\Spec}{Spec}
\DeclareMathOperator{\End}{End}
\DeclareMathOperator{\uEnd}{\underline{End}}
\DeclareMathOperator{\Hom}{Hom}
\DeclareMathOperator{\Ext}{Ext}
\DeclareMathOperator{\Tor}{Tor}
\DeclareMathOperator{\Der}{Der}
\DeclareMathOperator{\Hilb}{Hilb}
\DeclareMathOperator{\Alg}{Alg}
\DeclareMathOperator{\rank}{rank}
\DeclareMathOperator{\ord}{ord}
\DeclareMathOperator{\pd}{pd}
\DeclareMathOperator{\Frac}{Frac}

\DeclareMathOperator{\Per}{Per}
\DeclareMathOperator{\Ass}{Ass}
\newcommand{\sEnd}{\mathscr{E}\!nd}
\makeatletter
\renewcommand\paragraph{\@startsection{paragraph}{4}{\z@}{\medskipamount}{-.5em}{\normalfont\bfseries}}
\makeatother
\DeclareMathOperator{\Jac}{Jac}
\newcommand{\tnu}{\tilde\nu}
\newlist{thmlist}{enumerate}{1}
\setlist[thmlist]{label=\textup{(\arabic*)}, ref=\arabic*, leftmargin=*, topsep=3pt, itemsep=2pt}

\title[Centres of contraction algebras]{Centres and semisimple deformations of contraction algebras}
\author{Joseph Karmazyn}
\email{j.h.karmazyn@gmail.com}
\author{Shengyong Pan}
\address{School of Mathematics and Statistics, Beijing Jiaotong University, Beijing 100044, China}
\email{shypan@bjtu.edu.cn}
\date{}

\begin{document}

\begin{abstract}
Let $A_{\mathrm{con}}$ be the contraction algebra of a $3$-fold flop of an irreducible curve, with genus-zero Gopakumar--Vafa invariants $n_1, \dots, n_l$. We show that $\prod_j M_j(\mathbb{C})^{n_j}$ is the only semisimple algebra to which $A_{\mathrm{con}}$ deforms. Hence $A_{\mathrm{con}}$ encodes its Gopakumar--Vafa invariants intrinsically, answering a question of Wemyss. Hua and Toda showed that this deformation exists, and Jelisiejew and Smoktunowicz proved its uniqueness in types A and D. The key step in our proof is that the centre of $A_{\mathrm{con}}$ has dimension $n_1 + \dots + n_l$. To prove this, we lift the contraction algebra of the universal flop to the Cartan subalgebra, where it becomes a corner of a deformed preprojective algebra. We show that the universal contraction algebra and its zeroth Hochschild homology are free over a polynomial ring, with presentations whose determinants are products of roots. Pulling these presentations back along the arc that classifies the flop gives the count, and recovers Toda's formula for $\dim A_{\mathrm{con}}$. We also prove Hua and Toda's flatness conjecture, and determine which Gopakumar--Vafa invariants occur in lengths $3$ to $6$.
\end{abstract}

\maketitle

\section{Introduction}\label{sec:intro}

Contraction algebras are remarkably strong invariants of flops: although they are only finite-dimensional algebras, their representation theory recovers the underlying singularity \cite{HK, JKM, KLW}. This gives the opportunity to study flops, and their geometric invariants, through finite-dimensional algebras that are accessible to computation and to the wealth of techniques of noncommutative algebra. The Gopakumar--Vafa (GV) invariants are one such geometric invariant, and it is expected that the contraction algebra encodes them in a particularly natural way: that it deforms flatly to exactly one semisimple algebra, whose numbers of matrix factors of each size are the GV invariants. That it deforms to this algebra follows from work of Hua and Toda \cite{HT}; in this paper we prove that it deforms to no other. We do so by reducing the question to intrinsic invariants of the algebra, such as its centre and zeroth Hochschild homology, which we compute through the roots of a Dynkin diagram.

\paragraph{Contraction algebras and GV invariants.} Let $f\colon Y \to X = \Spec R$ be a $3$-fold flopping contraction of an irreducible rational curve $C$, with $Y$ smooth and $R$ complete local and isolated compound Du Val (cDV). Donovan and Wemyss \cite{DW} associate to $f$ its contraction algebra $\Acon$, the algebra representing noncommutative deformations of $C$. It is a finite-dimensional symmetric algebra, it determines $R$ up to isomorphism and the length $l$ of the flop ($1 \le l \le 6$), and it is related to the genus-zero GV invariants $n_j$ of $f$ by Toda's formula \cite{T}
\[
\dim \Acon = \sum_{j=1}^{l} j^2\, n_j .
\]

\needspace{5\baselineskip}
\paragraph{Wemyss's question.} Wemyss asked (see \cite{JS, Smo2}) whether $\Acon$ always deforms flatly to the semisimple algebra
\[
\SGV = \prod_{j=1}^{l} M_j(k)^{\,n_j}
\]
and whether this is the only semisimple algebra to which it deforms.

From the geometric point of view, $\Acon$ should deform to $\SGV$, since deforming the flop to a generic one splits $C$ into $n_j$ rigid $(-1,-1)$-curves of class $j[C]$. Hua and Toda showed that for a generic one-parameter deformation of the flop the contraction algebras of the nearby fibres are exactly $\SGV$, and observed that, with Toda's formula, this family is flat \cite[Thm~4.2 and the remark after it]{HT}. Indeed, the contraction algebra is a quotient of the endomorphism algebra of a tilting bundle (\S\ref{ss:flops}), so it commutes with base change, and the special fibre $\Acon$ has the same dimension as $\SGV$.\footnote{Hua and Toda define the contraction algebra through their own tilting bundle, and their algebra is isomorphic to $\Acon$ (Remark~\ref{rem:HT}). In the background part of its introduction, \cite{K} attributes the uniqueness of the semisimple deformation to Hua and Toda \cite[Prop~1.7]{K}, but this is an inaccurate citation, as their theorem concerns only the generic deformation (Appendix~\ref{app:corrections}).}

\paragraph{Results.} In this paper we answer Wemyss's question for every flop of an irreducible curve. The main result is Theorem~C: $\Acon$ deforms to $\SGV$ and to no other semisimple algebra. The key step is Theorem~B(1), that $\dim Z(\Acon) = n_1 + \dots + n_l$, the number of curves into which a generic deformation of the flop splits $C$. We prove it by passing to the universal flop of each length, whose contraction algebra Theorem~A describes through the roots of a Dynkin diagram. We also prove Hua and Toda's conjecture that the contraction algebra, its zeroth Hochschild homology and its centre are flat in every deformation of the flop (Theorem~B(3)), and determine the GV invariants that occur in lengths $3$--$6$ (Proposition~E). Finally, we settle the conjecture of Jelisiejew and Smoktunowicz quoted below. Theorem~C also shows that $[\Acon]$ is a smooth point of the scheme $\Ass_d$ of unital associative algebra structures on $k^d$, where $d = \dim \Acon$, and that $\dim Z$ is constant near it, which proves all of the conjecture except smoothness of the centre map. That smoothness fails for the Laufer flop when the target is a moduli space of algebras (Proposition~D), and holds for every flop when the target is a Grassmannian (\S\ref{ss:centremap}).

\begin{thmC}[Algebra structures]
Let $f$ be a flop as above, $\Acon$ its contraction algebra and $d = \dim \Acon$. Then $[\Acon]$ is a smooth point of $\Ass_d$, and lies on a single irreducible component of $\Ass_d$, the closure of the orbit of $\SGV$. Hence $\SGV$ is the only semisimple algebra to which $\Acon$ deforms, and $\dim Z$ is constant on a neighbourhood of $[\Acon]$.
\end{thmC}

Theorem~C answers Wemyss's question in the form stated in \cite[Question~1.2]{JS}. We now say what was known before, and then how we prove our results.

\paragraph{Previous work.} Whilst explicit semisimple deformations have been constructed for several families, uniqueness was previously known only in special cases. Jelisiejew and Smoktunowicz answered Wemyss's question in types A and D, where such algebras deform to at most one semisimple algebra and the deformations are explicit \cite[\S3, Thm~3.7]{JS}. For the contraction algebra $L_1$ of the Laufer flop, Bowman, Pulji\'c and Smoktunowicz had constructed an explicit deformation to $k^5 \times M_2(k)$ \cite[Thm~4.2]{BPS}. Smoktunowicz showed that one formal deformation cannot be flat to two non-isomorphic semisimple algebras \cite[Thm~2]{Smo2}, which gives uniqueness within a family but not across families. In type E, Brown and Wemyss announced an infinite family of $\mathsf{E}_6$ flops with normal forms $\mathsf{E}_{6,n}$ \cite[\S1.5, Table~1]{BW}, now proved in \cite[Thms~1.2, 1.3]{BW26} (see also \cite{Wem}). Explicit deformations of their Jacobi algebras to the expected semisimple algebras are constructed in \cite{SS}, extending a method of \cite[\S11]{Smo1}, assuming the expected dimensions, and in \cite{Pan}, which proves those dimensions, realises the algebras as contraction algebras of flops and shows that the deformation is unique.

Jelisiejew and Smoktunowicz also proposed a sharper conjecture, and verified its first part by computer in small dimensions. We quote it in our notation.

\begin{conjJS}[{\cite[Conj~1.4]{JS}}]
Let $\Acon$ be a contraction algebra and $d = \dim_k \Acon$.
\begin{thmlist}
\item A point $[\Acon]$ corresponding to $\Acon$ in the moduli space $\Ass_d$ of associative unital algebras is smooth; in other words, $\Acon$ is unobstructed.
\item The dimension of the centre is constant in an open neighbourhood of $[\Acon]$, and the associated map $A \mapsto Z(A)$ is smooth.
\end{thmlist}
\end{conjJS}

\begin{figure}[t]
\centering
\begin{tikzpicture}[
  font=\small,
  box/.style={draw, rounded corners=2pt, fill=white, align=center, inner sep=3pt,
              minimum height=0.85cm, minimum width=3.0cm},
  hbox/.style={box, minimum width=2.2cm},
  lab/.style={font=\scriptsize, align=center, fill=white, inner sep=1.2pt},
  arr/.style={-{Stealth[length=5pt]}, semithick},
]
\def\xa{0} \def\xb{4.7} \def\xc{9.4}
\def\yr{0} \def\yh{-2.0} \def\ya{-4.0} \def\yg{-6.0} \def\yb{-7.5}
\def\hoff{0.4}   %
\def\lane{-1.15} %

\node[font=\bfseries] at (\xa,1.0) {Lie theory};
\node[font=\bfseries] at (\xb,1.0) {Universal flop};
\node[font=\bfseries] at (\xc,1.0) {One flop};

\node[box] (R1) at (\xa,\yr) {roots $\alpha$, \ $\alpha_\bul$};
\node[box] (R2) at (\xb,\yr) {$f_j = \prod_{\alpha_\bul = j} \alpha$, \ $D_j$};
\node[box] (R3) at (\xc,\yr) {$n_j = \ord_w f_j(\nu)$};
\node[hbox] (H1) at (\xa+\hoff,\yh) {$\HH_0(\Pi_\bul)$};
\node[hbox] (H2) at (\xb+\hoff,\yh) {$\HH_0(\cC_l)$};
\node[hbox] (H3) at (\xc+\hoff,\yh) {$Z(\Acon)$};
\node[box] (A1) at (\xa,\ya) {$\Pi_\bul = e_\bul\Pi^\lambda(\Gamma)e_\bul$};
\node[box] (A2) at (\xb,\ya) {$\cC_l$ over $H_l$};
\node[box] (A3) at (\xc,\ya) {$\Acon \ \rightsquigarrow\ \SGV$\\[-1pt]
  {\scriptsize Thm~B(2), Thm~C}};
\node[box] (G1) at (\xa,\yg) {$\cY_\Gamma \to \h$};
\node[box] (G2) at (\xb,\yg) {$\cY_l \to \cX_l$};
\node[box] (G3) at (\xc,\yg) {$Y \to X$};

\node (B1) at (\xa,\yb) {$\h$};
\node (B2) at (\xb,\yb) {$\Spec H_l$};
\node (B3) at (\xc,\yb) {$\Spec k[[w]]$};
\draw[arr] (B1) -- node[lab, above] {$\pi$} (B2);
\draw[arr] (B3) -- node[lab, above] {$\nu$} (B2);
\draw[arr, dashed] (B3.south) to[out=-165, in=-15, looseness=0.45]
  node[lab, below] {$\tnu$ \ (Lemma~\ref{lem:lift}, Cor.~\ref{prop:rootformula})} (B1.south);
\foreach \c in {1,2,3} \draw[densely dotted, gray] (G\c.south) -- (B\c.north);

\draw[arr] (R1) -- (R2);
\draw[arr] (R2) -- node[lab, above] {Prop.~\ref{prop:GV}} (R3);
\draw[arr] (A2) -- node[lab, above] {Thm~A(1)} (A1);
\draw[arr] (A2) -- node[lab, above] {Lemma~\ref{lem:basechange}} (A3);
\draw[arr] (G1) -- (G2);
\draw[arr] (G3) -- (G2);

\foreach \c/\x/\la/\lb in {%
    1/\xa/{Prop.~\ref{prop:generic}}/{Thm~\ref{thm:squeeze}},
    2/\xb/{Thm~A(2),(3)}/{Thm~A(3)},
    3/\xc/{Cor.~\ref{cor:rootsum}}/{Thm~B(1)}}{
  \draw[arr] (\x+\lane,\yr-0.425) -- node[lab, pos=0.22] {\la} (\x+\lane,\ya+0.425);
  \draw[arr] (\x+\hoff,\yr-0.425) -- node[lab, right] {\lb} (H\c);
  \draw[arr, gray] (A\c.north -| H\c) -- (H\c);
  \draw[arr, gray] (G\c) -- (A\c);}

\begin{scope}[on background layer]
  \foreach \y/\t in {\yr/roots, \yh/$\HH_0$, \ya/algebras, \yg/geometry}{
    \fill[black!4] (-2.0,\y-0.7) rectangle (\xc+1.9,\y+0.7);
    \node[font=\scriptsize\itshape, rotate=90, gray] at (-1.8,\y) {\t};}
\end{scope}
\end{tikzpicture}
\caption{The objects of the paper, defined below, and how they are related. The columns lie over $\h$, $\Spec H_l$ and $\Spec k[[w]]$. The grey arrows from the geometry take the endomorphism algebra of a tilting bundle and kill $e_0$, and those from the algebras take $\HH_0$, which base changes as the algebras do. The arrow $\rightsquigarrow$ means `deforms to'. As $\Acon$ is symmetric, $\dim Z(\Acon) = \dim \HH_0(\Acon)$. Figure~\ref{fig:laufer} follows the Laufer flop through the figure.}\label{fig:story}
\end{figure}

\paragraph{The approach.} The basic idea is to pass from the roots of a Dynkin diagram to the GV invariants of a flop via the universal flops, translating between algebra and geometry along the way, as Figure~\ref{fig:story} shows. Katz and Morrison \cite{KM} showed that every flop of length $l$ is pulled back along a classifying arc from the universal flop of length $l$, and that the length itself comes from the roots: they attach to the flop a Dynkin diagram $\Gamma$ with a flopping vertex $\bul$, one pair $(\Gamma, \bul)$ for each length, and $l$ is the coefficient of the highest root at $\bul$. We show that the GV invariants are pulled back from the roots in the same way, and that the universal contraction algebras detect them. On the geometric side, Bryan, Katz and Leung described the GV invariants through the same roots and universal flops \cite{BKL}, and in these terms $n_j$ is the order of vanishing of $f_j$, the product of the roots with coefficient $j$ at $\bul$, along the classifying arc (Proposition~\ref{prop:GV}). On the algebraic side, the universal contraction algebra $\cC_l$ and its zeroth Hochschild homology are free over a polynomial subring of the base, and are cokernels of square matrices over the base with determinants $\prod_j f_j^{\,j^2}$ and $f_1 \cdots f_l$ (Theorem~A). Pulling these presentations back along the arc gives Theorem~B, as we outline after its statement, and a tangent-space count turns Theorem~B into Theorem~C.

In \cite{K} the universal flop is translated into noncommutative algebra: it carries a tilting bundle whose endomorphism algebra, the universal flopping algebra, has an explicit presentation. There the universal flop serves as a tool for constructing examples of contraction algebras, and here we also use it to prove results about every flop of a given length at once. This leaves just one algebra to consider for each length, the universal contraction algebra, from which the contraction algebra of each flop is obtained by base change. In order to compute with it we lift to the Cartan subalgebra, where it becomes a corner of a deformed preprojective algebra of Dynkin type. That algebra is obtained from the endomorphism algebra of a tilting bundle on the simultaneous resolution of the Kleinian singularity by killing one vertex.

\paragraph{The universal setting.} Let $\cC_l = A_l/A_l e_0 A_l$ denote the universal contraction algebra of length $l$, where $A_l$ is the universal flopping algebra over the base $H_l$ of the universal flop and $e_0$ is the idempotent of the summand $\cO$ of its tilting bundle. Let $\Gamma$ be the Dynkin diagram with flopping vertex $\bul$ and vertex idempotent $e_\bul$, let $\h$ be its Cartan subalgebra, and let $W_J$ be the Weyl group generated by the simple reflections of the set $J$ of vertices other than $\bul$ (whose curves are contracted in a generic hyperplane section of the flop). Then $H_l = k[\h]^{W_J}$ by the theorem of Katz and Morrison (\S\ref{ss:universal}). Moreover $H_l = k[t] \otimes P$, where $t$ is a $W_J$-invariant linear form and $P = k[\h_J]^{W_J} \subset H_l$, the invariants of the span $\h_J$ of the coroots of $J$, is a polynomial ring (\S\ref{ss:lie}). For a positive root $\alpha$ we write $\alpha_\bul$ for its coefficient at $\bul$, and $\Phi^+_\bul$ for the set of positive roots with $\alpha_\bul \neq 0$. For $j = 1, \dots, l$ we put
\[
f_j = \prod_{\alpha > 0,\ \alpha_\bul = j} \alpha ,
\]
the product of the positive roots whose coefficient at $\bul$ is $j$. Since $W_J$ fixes this coefficient, it permutes these roots, so $f_j \in H_l$. We write $D_j = \{f_j = 0\} \subset \Spec H_l$.

\begin{thmA}[The universal contraction algebra]
Let $\cC_l$ be the universal contraction algebra of length $l$, over $H_l = P[t]$.
\begin{thmlist}
\item $\cC_l \otimes_{H_l} k[\h] \cong e_\bul \Pi^\lambda(\Gamma) e_\bul$, where $\Pi^\lambda(\Gamma)$ is the deformed preprojective algebra of Crawley-Boevey and Holland \cite{CBH}, regarded as an algebra over $k[\h]$ by taking its parameters $\lambda$ as coordinates on $\h$.
\item $\cC_l$ is supported on $D_1 \cup \dots \cup D_l$, and at a general point of $D_j$ its fibre is the matrix algebra $M_j(k)$.
\item $\cC_l$ and $\HH_0(\cC_l)$ are free over $P$, of ranks $\sum_{\alpha > 0} \alpha_\bul^2$ and $|\Phi^+_\bul|$, and the characteristic polynomials of multiplication by $t$ on them are
\[
\prod_{j=1}^{l} f_j^{\,j^2} \qquad\text{and}\qquad f_1 f_2 \cdots f_l
\]
respectively, up to nonzero scalars. Hence each is the cokernel of a square matrix over $H_l$ with this determinant.
\end{thmlist}
\end{thmA}

Part (1) is a direct consequence of \cite{K}, and (2) follows from Crawley-Boevey's description of the simple modules of $\Pi^\lambda(\Gamma)$ \cite{CB} and the Ext formula of \cite{CBK}. The statements about $\cC_l$ in (3) follow from (1) and the results of Etingof and Rains \cite{ER} (their Cohen--Macaulay theorem and Hilbert series). Both $\cC_l$ and $\HH_0(\cC_l)$ are supported on $D_1 \cup \dots \cup D_l$, so they are torsion over $H_l$. Hence the freeness in (3) is over the polynomial subring $P$, not over $H_l$. Theorem~\ref{thm:Avertex} states all three parts at every vertex $\bul$ of every Dynkin diagram, and Theorem~A is its case of the flopping vertex. In terms of the deformed preprojective algebra alone, the corner $e_\bul\Pi^\lambda(\Gamma)e_\bul$, its zeroth Hochschild homology and its centre are free over $k[\h_J]$ (Proposition~\ref{prop:ER}, Theorem~\ref{thm:squeeze} and Corollary~\ref{cor:centre}), and $\HH_0(\Pi^\lambda(\Gamma))$ is the module of vector fields on $\h$ modulo those tangent to every root hyperplane (Proposition~\ref{prop:wholeHH}).

Much of the argument brings together known material. What is new is that we prove these results once for the universal contraction algebra of each length, read the determinants of the presentations in (3) off the roots, and show that $\HH_0(\cC_l)$, which controls the centre, is free.

\paragraph{Freeness of $\HH_0$.} The main new step is the freeness of $\HH_0(\cC_l)$ in Theorem~A(3), which we now explain. After base change to $\h$, $\cC_l$ is generically the matrix algebra $M_{\alpha_\bul}(k)$ along each root hyperplane $\alpha^\perp$ with $\alpha_\bul \neq 0$, so $\HH_0(\cC_l)$ is generically one copy of the functions on each such hyperplane. Freeness is a Cohen--Macaulay property: it says how these copies are glued where the root hyperplanes meet. This is what the count of the centre needs, because the classifying arc of every flop passes through the origin, where all these hyperplanes meet.

The gluing comes from the whole deformed preprojective algebra, through a classical fact about Weyl groups: the vector fields on $\h$ tangent to every root hyperplane form a free module, with basis the gradients of the basic invariants, the homogeneous generators of $k[\h]^W$. Hence the vector fields modulo these form a Cohen--Macaulay module, with one copy of the functions on each root hyperplane, and we show that this module is the zeroth Hochschild homology of the whole algebra. Applying a vector field on $\h$ to the right-hand side $\sum_i \lambda_i e_i$ of the defining relation of $\Pi^\lambda(\Gamma)$ gives a surjection $\kappa$ from the vector fields onto $\HH_0$, an algebraic counterpart of the Kodaira--Spencer map. It kills every vector field that lifts to a derivation of the algebra. As we know the vector fields tangent to the root hyperplanes explicitly, we can lift them, first to the simultaneous resolution of the Kleinian singularity and then to derivations of the endomorphism algebra of a tilting bundle there, which descend to the algebra. Hence they lie in the kernel of $\kappa$, and a length count at the root hyperplanes shows that they are all of it. Finally $\HH_0(\cC_l)$ is a direct summand of this module, so it inherits the gluing. Restricting this description to a line recovers the computation of Etingof, Latour and Rains \cite{ELR} without case analysis (Remark~\ref{rem:ELR}).

\paragraph{From roots to GV invariants.} Pulling Theorem~A back to a single flop and its deformations gives the following.

\begin{thmB}[Families of flops]
Let $f$ be a flop of an irreducible curve as above, of length $l$, with GV invariants $n_1, \dots, n_l$ (\S\ref{ss:flops}).
\begin{thmlist}
\item $\dim Z(\Acon) = \dim \HH_0(\Acon) = n_1 + \dots + n_l$.
\item $\Acon$ deforms to $\SGV$ (\S\ref{ss:moduli}), through a family obtained by deforming a truncation of its classifying arc.
\item Let $\mathcal Y \to \mathcal X = \Spec R_O$ be a flat projective deformation of $f$ over a complete local noetherian $k$-algebra $O$ with residue field $k$, with $R_O$ complete local noetherian, and let $A_O = \End_{\mathcal Y}(\cO \oplus \mathcal M_O)$ for the relative tilting bundle of Van den Bergh. Then the relative contraction algebra $A_O/A_Oe_0A_O$, its zeroth Hochschild homology and its centre are free $O$-modules, of ranks $\sum_j j^2 n_j$, $\sum_j n_j$ and $\sum_j n_j$, with fibres $\Acon$, $\HH_0(\Acon)$ and $Z(\Acon)$ at the closed point.
\end{thmlist}
\end{thmB}

Having checked in examples that $\dim \HH_0(\Acon) = n_1 + \dots + n_l$, Hua and Toda conjectured that $\Acon$ and its zeroth Hochschild homology and cohomology (the latter is the centre) are flat in every deformation of the flop \cite[Conj~4.3]{HT}. Part (1) confirms their count, and part (3) proves the conjecture, which Remark~\ref{rem:HT} deduces in their own setting, near the singular point of the central fibre. Part (2) recovers their flat deformation without Toda's formula, and the same argument recovers Toda's formula itself (Corollary~\ref{cor:toda}).

We now outline the proof of Theorem~B(1). Let $\nu\colon \Spec k[[w]] \to \Spec H_l$ be the classifying arc of the flop. In Lemmas~\ref{lem:basechange} and~\ref{lem:arcs} we show that $\Acon \cong \cC_l \otimes_{H_l} k[[w]]$. The hypersurface $D_j$ is the image of the root hyperplanes $\alpha^\perp$ with $\alpha_\bul = j$, and over a general point of it the fibre of the universal family, a surface, contains a single smooth $(-2)$-curve, of class $j[C]$. Hence $D_j$ is where curves of class $j[C]$ survive, and in Proposition~\ref{prop:GV} we show, following \cite{BKL}, that $n_j = \ord_w f_j(\nu)$, the order of contact of the arc with $D_j$. As the length of the cokernel of a square matrix over $k[[w]]$ is the order of its determinant, and $\HH_0$ commutes with base change, pulling back the presentations of Theorem~A(3) gives
\[
\dim \Acon = \ord_w \prod_{j} f_j(\nu)^{j^2} = \sum_{j} j^2 n_j, \qquad \dim \HH_0(\Acon) = \ord_w (f_1 \cdots f_l)(\nu) = \sum_{j} n_j .
\]
The first is Toda's formula, and the second, with the symmetry of $\Acon$, is Theorem~B(1). For example, the Laufer flop of Proposition~D below has $(n_1, n_2) = (5, 1)$, so $\dim \Acon = 5 + 4 = 9$ and $\dim Z(\Acon) = 6$ (Lemmas~\ref{lem:laufer} and~\ref{lem:lauferSmooth}). Example~\ref{ex:lauferarc} realises this flop by an explicit arc, using \cite{JKM,KLW}, and follows it through the proofs.

\paragraph{From the centre count to Theorem~C\@.} We deduce Theorem~C from Theorem~B by a tangent-space count. As $\Acon$ deforms to $\SGV$, the closure of the orbit of $\SGV$ is an irreducible component of $\Ass_d$ through $[\Acon]$. Lemma~\ref{lem:tangent} shows that the tangent space at $[\Acon]$ exceeds this component in dimension by $(\dim Z(\Acon) - n_1 - \dots - n_l) + (h^2 - h^1)$, where $h^i = \dim \HH^i(\Acon, \Acon)$. Contraction algebras have $h^1 = h^2$ by Lemma~\ref{lem:pan}, so by Theorem~B(1) the point $[\Acon]$ is smooth and lies on this component only (Proposition~\ref{prop:criterion}). In particular $\Acon$ determines each $n_j$ algebraically, as the number of factors $M_j(k)$ of its only semisimple deformation. This recovers \cite[Cor~4.6]{HT}, which answers \cite[Question~1.1]{JS} through the Donaldson--Thomas invariants of the noncommutative Hilbert schemes of $\Acon$ \cite[Thm~4.4]{HT} (see \cite{Dav} for refined invariants). In types A and D, Theorem~C recovers the uniqueness results of \cite{JS} by a method that does not depend on the type, and there part (1) of the conjecture of Jelisiejew and Smoktunowicz already follows from their centre computations and Lemma~\ref{lem:pan} (Remark~\ref{rem:AD}). The remaining part of their conjecture is smoothness of the map $A \mapsto Z(A)$ near $[\Acon]$.

\begin{propD}
Let $L_1 = k\langle x,y\rangle/(xy + yx,\ x^3 + y^2)$ be the contraction algebra of the Laufer flop $\mathsf{D}_{2,\infty}$ \cite[Table~1]{BW} (\S\ref{ss:laufer}). Then $[L_1]$ is a smooth point of $\Ass_9$, but its centre
\[
Z(L_1) \cong k[z_1,z_2,z_3]/(z_1z_3, z_2z_3, z_3^2, z_2^2, z_1^3, z_1^2z_2)
\]
is, as a quotient of $k[z_1,z_2,z_3]$ by these relations, a singular point of $\Hilb^6(\A^3)$. Hence the map $A \mapsto Z(A)$ is not smooth at $[L_1]$ as a map to the stack of $6$-dimensional commutative algebras, to the scheme $\Alg_6$ of commutative algebra structures on $k^6$, or to $\Ass_6$ (\S\ref{ss:centremap}), and the second half of part (2) of the conjecture of Jelisiejew and Smoktunowicz fails for these targets.
\end{propD}

\paragraph{GV invariants in each length.} The identity $n_j = \ord_w f_j(\nu)$ also determines which GV invariants occur.

\begin{propE}
Let $f$ be a flop of length $l \ge 3$. If $l = 3$, its GV invariants are $(6, 3, 1)$ or $(7, n_2, 1)$ with $n_2 \ge 4$. If $l = 4$, $5$ or $6$, they are $(6, 5, 2, 1)$, $(8, 6, 4, 2, 1)$ or $(6, 6, 4, 3, 2, 1)$ respectively. Each of these tuples occurs.
\end{propE}

This is Proposition~\ref{prop:typeE}. We prove it in \S\ref{ss:typeE} by a direct computation of $n_j$ along the lift of the classifying arc to $\h$, and we realise the tuples in lengths $4$--$6$ by the arcs of \cite[Examples~5.4, 5.5]{K}. In lengths $3$--$6$, Collinucci, De~Marco, Sangiovanni and Valandro \cite{CDSV} computed the GV invariants of one-parameter families of flops linear in the arc parameter, and Proposition~E determines those of all flops. In lengths $4$--$6$ the invariants, and hence $\dim \Acon$ and $\dim Z(\Acon)$, depend only on the length. Brown and Wemyss also determine the GV invariants and $\dim \Acon$ of every flop, as a corollary of their classification of flops \cite[Cor~1.6]{BW26}; in lengths $3$--$6$ they agree with Proposition~E.

\paragraph{Scope.} The results concern contraction algebras of flops. Brown and Wemyss show that every finite-dimensional Jacobi algebra of a potential is the contraction algebra of a flop of a smooth $3$-fold, except those isomorphic to the Jacobi algebra of a potential in two variables of the form $x^3 + O_4$ ($x^3$ plus terms of order at least four) \cite[Cor~1.12]{BW}. Theorem~C applies to all the others. For the exceptions, Theorem~C applies once they are known to be contraction algebras, as \cite[Cor~4.12]{Pan} proves for the Jacobi algebras of $x^3 + xy^3 + p(y)$ with $4 \le \ord p < \infty$ (\S\ref{ss:typeE}). Flops of reducible curves, whose contraction algebras have several vertices, are not treated here. We expect the methods to extend to them: this would need the analogue of Sections~\ref{sec:cartan} and~\ref{sec:squeeze} at a set of vertices, with the GV invariants defined by counting curves in a generic perturbation, as in \cite[\S3]{NW}. We leave the details to later work.

\paragraph{Structure of the paper.} In Section~\ref{sec:prelim} we collect background and prove, at every vertex, the tilting result of \cite{K} that Sections~\ref{sec:cartan}--\ref{sec:flops} use (Proposition~\ref{prop:tilting}), and in Section~\ref{sec:reduction} we reduce Theorem~C to Theorem~B. In Sections~\ref{sec:cartan} and~\ref{sec:squeeze} we work over the Cartan subalgebra, at any vertex of a Dynkin diagram. In Section~\ref{sec:cartan} we show that the corner of the deformed preprojective algebra is free and compute its rank, and in Section~\ref{sec:squeeze} we show that its zeroth Hochschild homology is free and prove Theorem~A. Of Section~\ref{sec:prelim}, Section~\ref{sec:cartan} uses only \S\S\ref{ss:commalg}, \ref{ss:lie}, \ref{ss:dpa} and Lemma~\ref{lem:corner}, and Section~\ref{sec:squeeze} also uses \S\ref{ss:moduli} and the geometry of \S\S\ref{ss:kleinian}--\ref{ss:universal}. In Section~\ref{sec:flops} we pull Theorem~A back along the classifying arc of a flop and prove Theorems~B and~C, using Sections~\ref{sec:cartan} and~\ref{sec:squeeze} only through Theorem~A and Proposition~\ref{prop:support}. In Section~\ref{sec:centremap} we prove Proposition~D, in Section~\ref{sec:examples} we write each GV invariant as a sum over roots, prove Proposition~E and work out explicit flops, and in Appendix~\ref{app:corrections} we correct some slips in \cite{K}. Theorems~A, B and~C and Propositions~D and~E use no computer calculation. The notation is summarised at the end of \S\ref{ss:universal}.

\paragraph{Conventions.} We fix the following conventions.
\begin{description}[leftmargin=0pt, itemsep=3pt, font=\normalfont\itshape]
\item[Base field] We work over $k = \C$ throughout, as do the sources on flops and contraction algebras that we use. The algebras $\Pi^\lambda(\Gamma)$, their corners and $\cC_l$ are defined over $\Q$, and Propositions~\ref{prop:generic} and~\ref{prop:support} hold, with the same proofs, over any algebraically closed field of characteristic zero. We use this for the generic fibre in Lemma~\ref{lem:genericarc}, and at the generic points of the root hyperplanes in Proposition~\ref{prop:wholeHH}(2) and Corollary~\ref{cor:centre}.
\item[Local and graded settings] Flops, their deformations and their contraction algebras are complete local (\S\ref{ss:flops}). The universal families of \S\S\ref{ss:kleinian}--\ref{ss:universal}, and the algebras of Sections~\ref{sec:cartan} and~\ref{sec:squeeze}, are algebraic over the affine spaces $\h$ and $\Spec H_l$, and graded by a $k^\times$-action (\S\ref{ss:simres}); for them, a semi-universal deformation means one whose completion at the origin is semi-universal. We pass from the graded setting to a single flop by base change along an arc and completion at the singular point (Lemmas~\ref{lem:basechange} and~\ref{lem:arcs}).
\item[References] A reference in the heading of a numbered statement gives its source: a bare reference means that the statement is quoted, ``extending'' that we prove it in the generality we need, and ``recovering'' that we give a new proof of a known result. Standard facts of commutative algebra and algebraic geometry are used by name: the local criterion of flatness \cite[Thm~22.3]{Mat}, the Auslander--Buchsbaum formula, Weierstrass preparation, cohomology and base change, Grothendieck's existence theorem and Zariski's main theorem.
\item[Quivers] Paths are composed left to right.
\end{description}

\section{Preliminaries}\label{sec:prelim}

In this section we fix notation and gather the background we need. Most of it is standard. It is either cited from the literature, with the source in the heading of each cited statement, or it is a known fact proved in the form we need, such as Lemma~\ref{lem:symfamily}, stated in families, and Lemma~\ref{lem:Ysmooth}, the smoothness of the semi-universal deformation of the partial resolution. The one new result is Proposition~\ref{prop:tilting} in \S\ref{ss:partial}, which extends \cite[Prop~4.21]{K} to every vertex, and Lemma~\ref{lem:corner} records what it gives for the corner. The section ranges over algebra structures and Hochschild homology, contraction algebras and their GV invariants, Lie theory and deformed preprojective algebras, and the deformations of Kleinian singularities and their partial resolutions that give the universal flops.

\subsection{Algebra structures and Hochschild homology}\label{ss:moduli}

We first recall the scheme of algebra structures, in which Theorem~C is set, and the facts about Hochschild homology that we use to compute centres. The scheme and its deformations are treated in \cite[\S2]{JS}.

\paragraph{The scheme of algebra structures.} Fix a $d$-dimensional vector space $V$ with a distinguished vector $1$. Let $\Ass_d$ be the affine scheme of associative unital multiplications on $V$ with unit $1$, and $G \subset GL(V)$ the stabiliser of $1$, so $\dim G = d^2 - d$. For an algebra $A = (V, \mu)$, the orbit $G\cdot[A]$ has dimension $\dim G - \dim \Der(A)$, and the Zariski tangent space $T_{[A]}\Ass_d$ is the space of normalised Hochschild $2$-cocycles, so
\[
\dim T_{[A]} \Ass_d = \dim (G\cdot[A]) + h^2(A), \qquad \dim \Der(A) = (d - \dim Z(A)) + h^1(A),
\]
where $h^i(A) = \dim \HH^i(A, A)$. Here $\HH^1(A, A)$ is the space of outer derivations and $\HH^2(A, A)$ that of first-order deformations modulo trivial ones.

\paragraph{Deformations.} We say $A$ \emph{deforms to} $S$ if $[A]$ lies in the closure of $G\cdot[S]$. Equivalently \cite[Thm~2.20]{JS}, there is a flat family over $k[[s]]$ with special fibre $A$ and generic fibre $S$ over an algebraic closure of $k((s))$.\footnote{In Theorem~B(2) we use only that such a family gives $[A] \in \overline{G\cdot[S]}$, which is the implication (1)$\Rightarrow$(4) of Theorem~2.20 of \cite{JS}, shown in the proof of their Proposition~2.17.} If $S$ is semisimple then $\HH^2(S) = 0$, so $\dim T_{[S]}\Ass_d = \dim (G\cdot[S])$. Hence $\Ass_d$ is smooth at $[S]$ and $G\cdot[S]$ is open \cite[Thm~2.14]{JS}.

\paragraph{Symmetric algebras.} A finite-dimensional algebra $A$ over a field $K$ is \emph{symmetric} if it has a \emph{symmetrising form}, a linear form $\operatorname{tr}\colon A \to K$ vanishing on $[A, A]$ with $(a, b) \mapsto \operatorname{tr}(ab)$ nondegenerate. The form identifies $Z(A)$ with $\HH_0(A)^* = (A/[A, A])^*$. For an algebra $A$ over a commutative ring $O$, $\HH_0(A)$ is the cokernel of the $O$-linear map $A \otimes_O A \to A$, $a \otimes b \mapsto ab - ba$, so it commutes with base change: $\HH_0(A \otimes_O O') = \HH_0(A) \otimes_O O'$.

The identification of $Z(A)$ with $\HH_0(A)^*$ also holds in families (this is used in Corollary~\ref{cor:centre} and Theorem~B(3)):

\begin{lem}\label{lem:symfamily}
Let $O$ be a local ring with residue field $K$, and $A$ an $O$-algebra that is finite and free as an $O$-module, such that $\HH_0(A) = A/[A, A]$ is a free $O$-module and $A \otimes_O K$ is a symmetric $K$-algebra. Then there is an $O$-linear map $\operatorname{tr}_O\colon A \to O$ vanishing on $[A, A]$ for which $a \mapsto \operatorname{tr}_O(a\,\cdot)$ is an isomorphism $A \to \Hom_O(A, O)$, and
\[
Z(A) \;\cong\; \Hom_O\big(\HH_0(A), O\big),
\]
a free $O$-module of the same rank as $\HH_0(A)$. This is compatible with base change: $Z(A \otimes_O O') \cong Z(A) \otimes_O O'$ for every map of local rings $O \to O'$; in particular $Z(A) \otimes_O K \cong Z(A \otimes_O K)$.
\end{lem}

\begin{proof}
Let $\operatorname{tr}\colon A \otimes_O K \to K$ be a symmetrising form. As $\HH_0$ commutes with base change, $\operatorname{tr}$ is a linear form on $\HH_0(A \otimes_O K) = \HH_0(A) \otimes_O K$, and as $\HH_0(A)$ is free it lifts to an $O$-linear map $\HH_0(A) \to O$; composing with $A \to \HH_0(A)$ gives $\operatorname{tr}_O$, with $\operatorname{tr}_O([A, A]) = 0$. The map $A \to \Hom_O(A, O)$, $a \mapsto \operatorname{tr}_O(a\,\cdot)$, is a map of free $O$-modules of the same rank that reduces modulo the maximal ideal to $a \mapsto \operatorname{tr}(a\,\cdot)$, an isomorphism since $\operatorname{tr}$ is nondegenerate. Hence this map is surjective by Nakayama's lemma, and so an isomorphism. For $z, b, c \in A$, $\operatorname{tr}_O(z[b, c]) = \operatorname{tr}_O(zbc) - \operatorname{tr}_O(bzc) = \operatorname{tr}_O([z, b]c)$, so $z$ is central if and only if $\operatorname{tr}_O(z\,\cdot)$ vanishes on $[A, A]$. This gives $Z(A) \cong \Hom_O(\HH_0(A), O)$. After a base change $O \to O'$, $\operatorname{tr}_O \otimes O'$ has the same two properties, so $Z(A \otimes_O O') \cong \Hom_{O'}(\HH_0(A \otimes_O O'), O') \cong Z(A) \otimes_O O'$, as $\HH_0$ commutes with base change.
\end{proof}

\paragraph{Hochschild homology.} For a ring $A$ and an idempotent $e \in A$,
\begin{equation}\label{eq:idem}
\HH_0(A/AeA) = \HH_0(A)\big/\operatorname{im}\HH_0(eAe), \qquad\text{since } aeb \equiv e(ba)e \bmod [A, A].
\end{equation}
If moreover $e$ is full, say $1 = \sum_k x_key_k$, then the inclusion induces an isomorphism $\HH_0(eAe) \cong \HH_0(A)$, with inverse $[a] \mapsto \sum_k [ey_kax_ke]$ (Morita invariance).

\subsection{Two lemmas from commutative algebra}\label{ss:commalg}

In Sections~\ref{sec:cartan}\nobreakdash--\ref{sec:flops} we will show that modules are free in two ways, by counting generators and through the Cohen--Macaulay property, and we will compute lengths of modules over discrete valuation rings. We record the facts we need here.

\begin{lem}\label{lem:freeness}
Let $\Lambda$ be a polynomial ring over $k$ on generators of positive degree, and $M$ a finite graded $\Lambda$-module.
\begin{thmlist}
\item (Graded Nakayama.) $\dim_k M/\Lambda_+M \ge \rank_\Lambda M$, with equality if and only if $M$ is free.
\item If $\dim M = \dim \Lambda$, then $M$ is free if and only if it is Cohen--Macaulay. The same holds for a finite module over a regular local ring $\Lambda$. If $R$ is a commutative graded $\Lambda$-algebra with $R_0 = k$ (resp.\ a commutative local $\Lambda$-algebra with $\Lambda \to R$ local) and $M$ is an $R$-module, finite over $\Lambda$, then $M$ has the same depth and dimension over $R$ as over $\Lambda$.
\item If $M$ is the cokernel of a square matrix over $\Lambda$ with nonzero determinant (graded, or with $\Lambda$ regular local), then $M$ is zero or Cohen--Macaulay of dimension $\dim \Lambda - 1$.
\end{thmlist}
\end{lem}

\begin{proof}
(1) A homogeneous basis of $M/\Lambda_+M$ lifts to a minimal generating set, giving a surjection $\Lambda^\mu \to M$ with $\mu = \dim_k M/\Lambda_+M$, so $\rank_\Lambda M \le \mu$. If equality holds, the kernel is torsion-free of rank zero, hence zero. (2) A free module is Cohen--Macaulay of dimension $\dim \Lambda$. Conversely, if $M$ is Cohen--Macaulay then $\operatorname{depth}_\Lambda M = \dim \Lambda$, so $\pd_\Lambda M = 0$ by the Auslander--Buchsbaum formula (applied, in the graded case, at the irrelevant ideal). For the last statement, $M$ is a faithful module over $\bar R = R/\operatorname{ann}_R M$, which is finite over $\Lambda$ because it embeds in $\End_\Lambda (M)$, and $M$ has the same depth and dimension over $R$ as over $\bar R$. Let $\mathfrak m$ be $\Lambda_+$ (resp.\ the maximal ideal of $\Lambda$). As $\bar R$ is finite over $\Lambda$, its irrelevant (resp.\ maximal) ideal has a power inside $\mathfrak m\bar R$, so the depth of $M$ over $\bar R$ is its depth along $\mathfrak m\bar R$, which the Koszul complex on generators of $\mathfrak m$ computes, and so equals its depth over $\Lambda$. Moreover, $\dim_{\bar R} M = \dim \bar R = \dim \Lambda /\operatorname{ann}_\Lambda M = \dim_\Lambda M$, as $\bar R$ is finite over $\Lambda /\operatorname{ann}_\Lambda M$. (3) The matrix is injective, so $\pd_\Lambda M \le 1$ and $\operatorname{depth}_\Lambda M \ge \dim \Lambda - 1$ by Auslander--Buchsbaum; and $M$ is torsion, so $\dim M \le \dim \Lambda - 1$.
\end{proof}

\begin{lem}\label{lem:dvr}
Let $V$ be a discrete valuation ring with uniformiser $w$, and $M$ a finite torsion $V$-module.
\begin{thmlist}
\item If $M/w M$ has dimension $r$ over the residue field, then $M$ has length at least $r$, with equality if and only if $w M = 0$.
\item If $M$ is the cokernel of a square matrix $\varphi$ over $V$ with $\det\varphi \neq 0$, then $M$ has length equal to the valuation of $\det\varphi$.
\end{thmlist}
\end{lem}

\begin{proof}
(1) $M$ is a sum of $r$ nonzero cyclic modules $V/w^{a_i}$, of total length $\sum a_i \ge r$, with equality if and only if every $a_i = 1$. (2) Bring $\varphi$ to diagonal form by row and column operations.
\end{proof}

\subsection{Flops and contraction algebras}\label{ss:flops}

We now recall flops, their contraction algebras and GV invariants, the facts about perverse sheaves that we use to recognise tilting bundles, and how pushforwards behave in families. Contraction algebras come from \cite{DW}, their symmetry and periodicity from \cite{JKM}, perverse sheaves and their projective generators from \cite[\S3]{VdB}, and the GV invariants from \cite{K08}, \cite{BKL} and \cite{T}. The precise statements on contraction algebras and perverse sheaves are collected in Proposition~\ref{prop:Acon} and Theorem~\ref{thm:VdB}, with their sources.

\paragraph{Flops.} Let $f\colon Y \to X = \Spec R$ be a flopping contraction of an irreducible rational curve $C$, where $Y$ is smooth, $R$ is a complete local isolated cDV (compound Du Val) singularity and $Rf_*\cO_Y = \cO_X$. The length $l$ of $f$ is the multiplicity of $C$ in the fibre over the closed point; $1 \le l \le 6$. We call such an $f$ a flop, of length $l$.

\paragraph{Contraction algebras.} By Van den Bergh, $Y$ carries a tilting bundle $\cO_Y \oplus M$ that is a projective generator of ${}^0\Per(Y/X)$ (Theorem~\ref{thm:VdB} below). Put $A = \End_Y(\cO_Y \oplus M) \cong \End_R(R \oplus f_*M)$. Let $e_0 \in A$ be the idempotent corresponding to $R$. The contraction algebra is
\[
\Acon = A / A e_0 A \;\cong\; \uEnd_R(f_*M),
\]
the stable endomorphism algebra of $f_*M$.

\paragraph{Duals and opposites.} Every statement in this paper is invariant under passing to the opposite algebra, for the following reason. Our $M^\vee$ is an extension of a line bundle $L$ with $\deg L|_C = 1$ by $\cO_Y^{\oplus(l - 1)}$. This is the convention of both \cite{K} and \cite[\S\S2.3--2.4]{DW} (the former writes $\End_Y(\cO_Y \oplus M)^{\mathrm{op}}$, and the latter use right modules with functional composition; we compose paths left to right). The ${}^{-1}\Per$ generator $\cO_Y \oplus M^\vee$ gives the opposite algebra, equivalently the same construction for the flopped contraction $Y^+ \to X$, since ${}^{-1}\Per(Y/X) \simeq {}^{0}\Per(Y^+/X)$ \cite[Thm~4.4.2]{VdB}. On $\Ass_d$, $A \mapsto A^{\mathrm{op}}$ is a $G$-equivariant automorphism, and $\SGV^{\mathrm{op}} \cong \SGV$.

We use the following facts about $\Acon$.

\begin{prop}[{Contraction algebras; \cite[Thms~1.1, 1.2]{DW}, \cite[\S1.2, Thm~3.2.5]{JKM}}]\label{prop:Acon}
Let $f$ be a flop.
\begin{thmlist}
\item (Deformations.) $\Acon$ is finite-dimensional and represents the noncommutative deformations of~$C$.
\item (Symmetry.) $\Acon$ is a symmetric algebra.
\item (Periodicity.) $\Acon$ is $4$-periodic as a bimodule with trivial twist: $\Omega^4_{\Acon^e}(\Acon) \cong \Acon$ in the stable category of $\Acon$-bimodules.
\end{thmlist}
\end{prop}

Part (3) applies because our flops are in the setting of Jasso, Keller and Muro, an isolated cDV singularity with a crepant resolution. Their theorem, in which $\Lambda = \Acon$ \cite[Setting~3.2.1]{JKM}, states that the restricted universal Massey product, an element of $\Ext^4_{\Lambda^e}(\Lambda, \Lambda)$, ``can be represented by an exact sequence $0 \to \Lambda \to P_3 \to P_2 \to P_1 \to P_0 \to \Lambda \to 0$ with projective middle terms. In particular, $\Omega^4_{\Lambda^e}(\Lambda) \cong \Lambda$ in the stable category of $\Lambda$-bimodules.'' Both ends of the sequence are $\Lambda$ itself, so the twist is trivial, which Lemma~\ref{lem:pan} needs. Jasso, Keller and Muro deduce this from \cite[Cor~4.5.18]{JM} (Cor~4.5.17 in the version they cite).

\paragraph{Perverse sheaves.} Van den Bergh also characterises the projective objects and the projective generators of ${}^{0}\Per$ and ${}^{-1}\Per$. We use this in Propositions~\ref{prop:simres} and~\ref{prop:tilting}, Lemma~\ref{lem:basechange} and Theorem~B(3). He considers a projective morphism $f'\colon Y' \to X' = \Spec R'$ of noetherian schemes with $Rf'_*\cO_{Y'} = \cO_{X'}$ and fibres of dimension at most one, and the category $\mathcal V$ of vector bundles $\mathcal P$ on $Y'$ generated by global sections with $H^1(Y', \mathcal P^\vee) = 0$.

\begin{thm}[{Van den Bergh \cite[Props~3.2.6, 3.2.7, Cor~3.2.8, Lemma~3.4.5, Prop~3.5.4, Thm~3.5.5, (3.11)]{VdB}}]\label{thm:VdB}
In this setting, the following hold.
\begin{thmlist}
\item (Projective objects.) The projective objects of ${}^{-1}\Per(Y'/X')$ are precisely the objects of $\mathcal V$, and those of ${}^{0}\Per(Y'/X')$ are precisely their duals.
\item (Projective generators.) The projective generators of ${}^{-1}\Per(Y'/X')$ are precisely the $\mathcal P$ in $\mathcal V$ with $\det \mathcal P$ ample and $\cO_{Y'}$ a direct summand of some $\mathcal P^{\oplus a}$, and those of ${}^{0}\Per(Y'/X')$ are precisely their duals. Every projective generator is a tilting bundle.
\item (Complete local base.) Suppose that $R'$ is complete local with residue field $k$. Then a line bundle is ample if and only if it has positive degree on each curve in the fibre over the closed point. The category ${}^{0}\Per(Y'/X')$ is Krull--Schmidt, being equivalent to the finite modules over an algebra finite over $R'$, and its indecomposable projective objects are $\cO$ and, for each irreducible curve $E_i$ of the fibre over the closed point, one bundle $N_i$, the dual of Van den Bergh's universal extension. The rank of $N_i$ is the multiplicity of $E_i$ in that fibre, and $\det N_i$ has degree $-\delta_{ij}$ on $E_j$.
\end{thmlist}
\end{thm}

In particular, in (3) every projective object of ${}^{0}\Per(Y'/X')$ is $\cO^{\oplus a} \oplus \bigoplus_i N_i^{\oplus b_i}$, and when that fibre is an irreducible curve $C$, a line bundle of degree one on $C$ is ample. With a contracting $k^\times$-action, the central fibre also decides ampleness. Suppose that $k^\times$ acts on $Y'$ and $X'$, compatibly with $f'$, and that the closure of every orbit in $X'$ contains the same point $0$. Then a $k^\times$-linearised line bundle on $Y'$ that is ample on the fibre over $0$ is ample. Indeed, as $f'$ is proper, the points of $X'$ over which its restriction to the fibre is ample form an open set, over which it is relatively ample. This set is $k^\times$-stable, as the bundle is linearised, and it contains $0$. Its complement is closed and $k^\times$-stable, so if it were non-empty it would contain the closure of an orbit, and hence $0$. Thus the set is all of $X'$, and as $X'$ is affine the bundle is ample.

\paragraph{GV invariants.} The genus-zero Gopakumar--Vafa invariants $n_1, \dots, n_l$ are the invariants defined by Katz from moduli of stable one-dimensional sheaves supported on $C$ \cite[Def~2.2, \S3]{K08}. For $j \le l$ let $C_j$ be the curve of multiplicity $j$ on $C$ cut out by the $j$-th symbolic power of the ideal of $C$ in a generic hyperplane section \cite[Def~1.4, Thm~1.5]{BKL}. Katz shows that $n_j$ is the multiplicity $k_j$ of $C_j$ in its Hilbert scheme, in which $[C_j]$ is an isolated point, so $k_j$ is the length of its complete local ring. Katz works on a small analytic neighbourhood of $C$ in a Calabi--Yau $3$-fold. The complete local ring of the Hilbert scheme at $[C_j]$ pro-represents the deformations of $C_j$ as a subscheme, so it depends only on the formal neighbourhood of $C$ (proof of Lemma~\ref{lem:truncation}). For complete local $f$ we therefore take $n_j := k_j$, the length of this ring, which agrees with Katz's invariant whenever $f$ is the completion of a flop in his setting. This length is finite for every flop here (Lemma~\ref{lem:truncation} and Proposition~\ref{prop:GV}). These $n_j$ are also the genus-zero GV invariants defined through Gromov--Witten theory, by Bryan, Katz and Leung. Toda proved $\dim \Acon = \sum j^2 n_j$ for flops of smooth quasi-projective $3$-folds \cite[Thm~1.1]{T}, and Corollary~\ref{cor:toda} recovers it for the complete local flops here. Finally $n_1$ is the dimension of the abelianisation of $\Acon$. Indeed the abelianisation represents the commutative deformations of $C$ \cite[Cor~3.2, Def~3.5]{DW}, so it is the complete local ring of the Hilbert scheme at $[C]$, whose length is $k_1 = n_1$.

\paragraph{Families.} The families of \S\S\ref{ss:kleinian}--\ref{ss:universal} and the flops of Section~\ref{sec:flops} are pulled back from one another along maps of bases. We record once how pushforwards behave under such base change, and the form of Nakayama's lemma that recovers $Rg_*\cO = \cO$ from the fibres.

\begin{lem}[Pushforwards in flat families]\label{lem:flatfamilies}
Let $S$ be a noetherian scheme, $g\colon \cY \to \cX$ a proper morphism of noetherian $S$-schemes with $\cY$ and $\cX$ flat over $S$, and $\mathcal E$ a vector bundle on $\cY$.
\begin{thmlist}
\item (Base change.) For every $S$-scheme $S'$, with $g'\colon \cY' \to \cX'$ the base change of $g$ to $S'$ and $\mathcal E'$ the pullback of $\mathcal E$, $Rg'_*\mathcal E' \cong Rg_*\mathcal E \otimes^{\mathbf L}_{\cO_{\cX}} \cO_{\cX'}$. In particular, if $R^ig_*\mathcal E = 0$ for $i > 0$, then $R^ig'_*\mathcal E' = 0$ for $i > 0$ and $g'_*\mathcal E' \cong g_*\mathcal E \otimes_{\cO_S} \cO_{S'}$.
\item (Rational fibres.) Suppose that every closed point of $\cX$ lies over a closed point of $S$ (as when $S$ and $\cX$ are of finite type over $k$, or are spectra of local rings with $\cX \to S$ local), and that $R(g_s)_*\cO = \cO$ for the fibre $g_s$ of $g$ over every closed point $s$ of $S$. Then $Rg_*\cO_{\cY} = \cO_{\cX}$.
\end{thmlist}
\end{lem}

\begin{proof}
(1) As $\cX$ is flat over $S$, $\cO_{\cX'} = \cO_{\cX} \otimes^{\mathbf L}_{\cO_S} \cO_{S'}$, and as $\cY$ is flat over $S$, $\cO_{\cY} \otimes^{\mathbf L}_{\cO_{\cX}} \cO_{\cX'} = \cO_{\cY} \otimes^{\mathbf L}_{\cO_S} \cO_{S'} = \cO_{\cY'}$. Hence $\cY$ and $\cX'$ are Tor-independent over $\cX$, and the isomorphism is Tor-independent base change. If $R^ig_*\mathcal E = 0$ for $i > 0$, the right-hand side has cohomology only in degrees $\le 0$ and the left-hand side only in degrees $\ge 0$, so both are the sheaf $g_*\mathcal E \otimes_{\cO_{\cX}} \cO_{\cX'} = g_*\mathcal E \otimes_{\cO_S} \cO_{S'}$ in degree $0$.

(2) Let $\mathcal Q$ be the cone of $\cO_{\cX} \to Rg_*\cO_{\cY}$, a bounded complex with coherent cohomology. By (1) with $S' = \Spec k(s)$, its derived restriction to the fibre $\cX_s$ is the cone of $\cO_{\cX_s} \to R(g_s)_*\cO$, which vanishes. If $\mathcal Q \neq 0$, let $\mathcal F$ be its highest nonzero cohomology sheaf. The restriction of $\mathcal F$ to $\cX_s$ is the highest cohomology sheaf of that derived restriction, so it vanishes for every closed point $s$. Every closed point $x$ of $\cX$ lies on such a fibre, so $\mathcal F \otimes k(x) = 0$ and $\mathcal F_x = 0$ by Nakayama's lemma. A coherent sheaf vanishing at every closed point is zero, so $\mathcal Q = 0$.
\end{proof}

\subsection{Lie-theoretic notation}\label{ss:lie}

Next we fix the Lie-theoretic notation: the base $H_J$, the equations $f_j$ of Theorem~A and the hypersurfaces $D_j$, with the facts about $H_J$ and the root hyperplanes that Sections~\ref{sec:cartan}--\ref{sec:flops} use. The facts about Weyl groups are from \cite[Ch.~V, \S5; Ch.~VI, \S1]{Bou}. Let $\Gamma$ be a Dynkin diagram with vertex set $\Gamma_0$, $n = |\Gamma_0|$, Cartan subalgebra $\h$, root lattice $Q_\Gamma$ and Weyl group $W$, and fix a vertex $\bul$. Let $J = \Gamma_0 \setminus \{\bul\}$, let $W_J \subset W$ be the subgroup generated by the simple reflections of $J$, and let $l$ be the coefficient of $\bul$ in the highest root (for the six pairs of \S\ref{ss:universal}, the length of the flop).

\paragraph{The base $H_J$.} Let $\varpi_\bul^\vee \in \h$ be the fundamental coweight dual to the simple root at $\bul$, and $\h_J$ the span of the coroots of $J$. Then $\h = k\cdot\varpi_\bul^\vee \oplus \h_J$ is a $W_J$-stable decomposition, $W_J$ fixes $\varpi_\bul^\vee$ and acts on $\h_J$ as a reflection group, so
\[
H_J := k[\h]^{W_J} = k[t] \otimes P, \qquad P = k[\h_J]^{W_J},
\]
where $t$ is a nonzero linear form vanishing on $\h_J$, and $P$ is a polynomial ring whose generator degrees have product $|W_J|$ (Chevalley's theorem). By the same theorem, $k[\h]$ is also a free $H_J$-module of rank $|W_J|$, so base change from $H_J$ to $k[\h]$ is faithfully flat. Equivalently, as $k[\h] = k[t] \otimes k[\h_J]$, $k[\h_J]$ is free of rank $|W_J|$ over $P$. Hence a finite graded $P$-module $M$ is free if and only if $M \otimes_P k[\h_J]$ is free over $k[\h_J]$, and then of the same rank. Averaging over $W_J$, $a \mapsto |W_J|^{-1}\sum_{\sigma \in W_J}\sigma\cdot a$, is an $H_J$-linear retraction of $k[\h]$ onto $H_J$, so $H_J$ is a direct summand of $k[\h]$ as an $H_J$-module.

\paragraph{Root hyperplanes.} For a positive root $\alpha$ write $\alpha_\bul$ for its coefficient at $\bul$; then $\alpha(\varpi_\bul^\vee) = \alpha_\bul$. Write $\lambda\cdot\alpha = \alpha(\lambda)$ for $\lambda \in \h$, so that roots are linear forms on $\h$ (\emph{root forms}), and $\alpha^\perp = \{\lambda : \lambda\cdot\alpha = 0\}$ for the root hyperplane, with generic point $\eta_\alpha$. The local ring of $k[\h]$ at the prime $(\alpha)$ is a discrete valuation ring with uniformiser $\alpha$, whose residue field $k(\eta_\alpha)$ is the function field of $\alpha^\perp$. If $\alpha_\bul \neq 0$ then $\alpha(\varpi_\bul^\vee) = \alpha_\bul \neq 0$, so projection along $\varpi_\bul^\vee$ restricts to an isomorphism $\alpha^\perp \cong \h_J$, and $\alpha^\perp \to \h_J/W_J = \Spec P$ is finite of degree $|W_J|$.

\paragraph{The equations $f_j$.} Write $\Phi^+$ for the set of positive roots of $\Gamma$ and $\Phi^+_\bul$ for those with $\alpha_\bul \neq 0$, and for $j = 1, \dots, l$ put
\[
f_j = \prod_{\alpha > 0,\ \alpha_\bul = j} \alpha \;\in\; H_J, \qquad D_j = \{f_j = 0\} \subset \Spec H_J;
\]
$f_j$ is $W_J$-invariant because $W_J$ fixes the coefficient at $\bul$ and so permutes these roots, and it is a product of $\#\{\alpha \in \Phi^+ : \alpha_\bul = j\}$ linear forms (Table~\ref{tab:roots} for the six pairs). Write $\pi\colon \h \to \h/W_J = \Spec H_J$ for the quotient by $W_J$, and $\Sigma \subset \h$ for the union of the pairwise intersections of distinct root hyperplanes. Then $\Sigma$ has codimension $2$ in $\h$, and so has $\pi(\Sigma)$ in $\Spec H_J$, as $\pi$ is finite.

\paragraph{Coordinates.} We identify $\h$ with $\Hom(Q_\Gamma, k)$, with coordinates $t_i$ for $i \in \Gamma_0$, where $t_i(\lambda)$ is the value of $\lambda$ on the simple root at $i$. The parameter $\lambda_i$ of the deformed preprojective algebra $\Pi^\lambda(\Gamma)$ of \S\ref{ss:dpa} at a vertex $i \in \Gamma_0$ is $t_i(\lambda)$, and in this way we regard $\Pi^\lambda(\Gamma)$ as an algebra over $k[\h]$. In these coordinates $\varpi_\bul^\vee$ is the point with $t_\bul = 1$ and $t_i = 0$ for $i \in J$, the functions $t_i$ ($i \in J$) restrict to coordinates on $\h_J$, and $t$ spans the $W_J$-invariant part of $\h^*$. For the six pairs we normalise $t$ as minus the invariant linear form called $t$ in \cite{K}, which is not $t_\bul$.

\paragraph{A sum over roots.} Let $h$ be the Coxeter number (unrelated to the $h^i$ of \S\ref{ss:moduli}), $A_\Gamma$ the Cartan matrix, and $(\ ,\ )$ the invariant form with roots of squared length $2$. Then
\[
\sum_{\alpha>0}\alpha(v)^2 = h\,(v,v) \quad (v \in \h), \qquad\text{so}\qquad \sum_{\alpha>0}\alpha_\bul^2 = h\,(A_\Gamma^{-1})_{\bul\bul}
\]
for every vertex $\bul$. Indeed, the form $\sum_{\alpha>0}\alpha(x)\alpha(v)$ is $W$-invariant, hence a multiple $c\,(x,v)$ of the invariant form, as $\h$ is an irreducible $W$-module. Comparing traces, $c\,n = \sum_{\alpha>0}(\alpha,\alpha) = 2|\Phi^+| = n h$, so $c = h$. Finally $(\varpi_\bul^\vee, \varpi_\bul^\vee) = (A_\Gamma^{-1})_{\bul\bul}$.

\subsection{Deformed preprojective algebras}\label{ss:dpa}

Here we recall the deformed preprojective algebras, and the facts about them that Sections~\ref{sec:cartan} and~\ref{sec:squeeze} use. Let $Q$ be a quiver with vertex set $Q_0$, $\overline{Q}$ its double, and $\lambda \in k^{Q_0}$. The deformed preprojective algebra \cite{CBH} is
\[
\Pi^\lambda(Q) = k\overline{Q} \Big/ \Big( \sum_{a \in Q_1} [a, a^*] - \sum_{i \in Q_0} \lambda_i e_i \Big).
\]
It depends only on the underlying graph, and we treat $\lambda$ as coordinates, obtaining an algebra over $k[\lambda]$. We keep the name $\Pi^\lambda(Q)$ for this family, and write $\Pi_\Gamma$ for it when $Q$ has type $\Gamma$, an algebra over $k[\h]$ (\S\ref{ss:lie}). For a point $\lambda$ we write $\Pi^\lambda$ for its fibre there (as in Proposition~\ref{prop:generic}(2)). We use the following facts.

\begin{prop}[{Deformed preprojective algebras; \cite[Thms~0.1, 0.4, Lemma~1.1, Thm~3.4, Cor~3.6]{CBH}, \cite[\S1, Thm~1.2]{CB}, \cite[Prop~2.6]{CBK}, \cite[\S\S2.1, 2.3, Thm~3.1]{ER}}]\label{prop:dpa}
Let $\Gamma$ be a Dynkin graph, and $\tGamma$ its extended Dynkin graph, with minimal imaginary root $\delta$.
\begin{thmlist}
\item (Affine case.) $\Pi^\lambda(\tGamma)$ is flat over $k[\lambda]$, and $e_0\Pi^\lambda e_0$ is commutative if and only if $\lambda\cdot\delta = 0$. Flatness is a PBW theorem: $\Pi^\lambda(\tGamma)$ is a corner of the algebra $k\langle x, y\rangle \rtimes G_\Gamma$ with $[x, y] = \lambda$, for the finite group $G_\Gamma \subset SL_2(k)$ of type $\Gamma$, and this algebra has a basis free over $k[\lambda]$.
\item (Dynkin case.) $\Pi^\lambda(\Gamma)$ is finite-dimensional for every $\lambda$, but the family is not flat.
\item (Simples.) A simple $\Pi^\lambda(\Gamma)$-module of dimension vector $\alpha$ exists if and only if $\alpha$ is a positive root with $\lambda\cdot\alpha = 0$ that is not a sum of two or more positive roots $\beta$ with $\lambda\cdot\beta = 0$, and the simple module is then unique. In particular every simple has dimension vector a positive root with $\lambda\cdot\alpha = 0$.
\item (Ext formula.) For any quiver $Q$, any $\lambda$, and finite-dimensional $\Pi^\lambda(Q)$-modules $M$, $N$ of dimension vectors $\alpha$, $\beta$, $\dim \Ext^1(M, N) = \dim \Hom(M, N) + \dim \Hom(N, M) - (\alpha, \beta)$, where $(\ ,\ )$ is the symmetrised Ringel form; for $\Gamma$ Dynkin this is the symmetric form of the Cartan matrix.
\item (Etingof--Rains.) With the vertices in degree $0$, the arrows in degree $1$ and the parameters $\lambda_i$ in degree $2$, $\Pi^\lambda(\Gamma)$ is a graded Cohen--Macaulay $k[\h]$-module of dimension $n - 1$, and its matrix Hilbert series, the matrix of the Hilbert series of the $e_i\Pi^\lambda(\Gamma)e_j$, is $(1 - u^{2h})(1 - u^2)^{-n}(I - uC_\Gamma + u^2I)^{-1}$, where $C_\Gamma = 2I - A_\Gamma$ is the adjacency matrix ($k[\h]$, $h$ and $A_\Gamma$ as in \S\ref{ss:lie}).
\end{thmlist}
\end{prop}

\begin{proof}
Only (5) needs an argument beyond the citations. Etingof and Rains prove that $\Pi^\lambda(\Gamma)$ is graded Gorenstein, and in particular Cohen--Macaulay in their sense: free over some central graded polynomial subalgebra $B$ over which it is finite. This $B$ is not $k[\h]$, over which $\Pi^\lambda(\Gamma)$ is not free by (2), and the notion does not depend on the choice of $B$, as follows. The subalgebra $B$ and the image of $k[\h]$ generate a commutative graded subalgebra of the centre, over which $\Pi^\lambda(\Gamma)$ is finite, so by Lemma~\ref{lem:freeness}(2) its depth and dimension over $k[\h]$ equal those over $B$. Over $k[\h]$, $\Pi^\lambda(\Gamma)$ has dimension $n - 1$, as its support is a union of hyperplanes (Proposition~\ref{prop:generic}, whose proof uses only (2)--(4) above). Hence $B$ also has dimension $n - 1$, and $\Pi^\lambda(\Gamma)$ is Cohen--Macaulay of dimension $n - 1$ over $k[\h]$.
\end{proof}

\subsection{Kleinian singularities}\label{ss:kleinian}

In this subsection we write down the Kleinian singularities and their semi-universal deformations explicitly, following \cite[\S\S6, 8]{Slo} and \cite[\S3, Tables~2, 3]{KM}, with the weights that grade everything built from them in \S\S\ref{ss:simres}--\ref{ss:universal}, and describe the discriminant. Let $\Gamma$, $\h$, $W$ and $n$ be as in \S\ref{ss:lie}. The Kleinian singularity of type $\Gamma$ is $\{F_0 = 0\} \subset \A^3$, with $F_0$ the normal form below, and its semi-universal deformation is $\cX_\Gamma = \{F = 0\} \subset \A^3 \times \h/W$, where
\[
F = F_0(x, y, z) + \sum_{j=0}^{n-1} u_j\,\mu_j(x, y, z),
\]
with $\mu_0 = 1$, monomials $\mu_j$, and $u_j$ coordinates on $\h/W$. Let $w_x, w_y, w_z$ and $d_F$ be the weights of the coordinates and of $F_0$, normalised so that a linear coordinate on $\h$ has weight $d_F/h$, where $h$ is the Coxeter number (the constant term $u_0$ has weight $d_F$ and is the basic invariant of degree $h$):
\begin{center}
\begin{tabular}{@{}lccc@{}}
\toprule
$\Gamma$ & $F_0$ & $(w_x, w_y, w_z; d_F)$ & $h$ \\
\midrule
$\mathsf{A}_{m}$ & $x^2 + y^2 + z^{m+1}$ & $(m+1, m+1, 2; 2m+2)$ & $m+1$ \\
$\mathsf{D}_m$ & $x^2 + y^2 z + z^{m-1}$ & $(m-1, m-2, 2; 2m-2)$ & $2m-2$ \\
$\mathsf{E}_6$ & $x^2 + y^3 + z^4$ & $(6, 4, 3; 12)$ & $12$ \\
$\mathsf{E}_7$ & $x^2 + y^3 + yz^3$ & $(9, 6, 4; 18)$ & $18$ \\
$\mathsf{E}_8$ & $x^2 + y^3 + z^5$ & $(15, 10, 6; 30)$ & $30$ \\
\bottomrule
\end{tabular}
\end{center}
In every type $w_x + w_y + w_z - d_F = d_F/h$. These weights define an action of $k^\times$, with positive weights, on $\A^3 \times \h/W$ preserving $\cX_\Gamma$, and on $\h$ and $\h/W_J$ compatibly with the quotient maps. On the affine varieties $\h$, $\h/W_J$ and $\cX_\Gamma$ it contracts everything to $0$ as the parameter tends to $0$; we call it the \emph{contracting action}.

The discriminant of $\cX_\Gamma \to \h/W$, over which the fibres are singular, is $\mathcal D = \{d_W = 0\}$ with $d_W = \prod_{\alpha>0}\alpha^2 \in k[\h]^W$, which is irreducible in $k[\h]^W$, and in particular reduced. Indeed, $W$ acts transitively on the roots, so a $W$-invariant factor of $d_W$ is a power of $\prod_{\alpha>0}\alpha$, which is anti-invariant. At a general point of $\mathcal D$ the fibre has a single singular point, an $\mathsf{A}_1$ singularity, and the critical value of $F$ at it has derivative $1$ with respect to $u_0$, as $\mu_0 = 1$.

\subsection{Simultaneous resolutions}\label{ss:simres}

We recall the geometry behind Proposition~\ref{prop:tilting}, following \cite[\S2.3]{K} and \cite[4.2, 8.9]{Slo}, for $\Gamma$, $\bul$ and $J$ as in \S\ref{ss:lie}. Write $\tGamma$ for the extended Dynkin diagram, with extending vertex $0$ and minimal imaginary root $\delta$, and $\cA_\Gamma = \Pi^\lambda(\tGamma)$ with $\lambda\cdot\delta = 0$; its base is identified with $\h$ through the coordinates $\lambda_i$, $i \in \Gamma_0$, since $\lambda\cdot\delta = 0$ determines $\lambda_0$. Let $\cX_\Gamma \to \h/W$ be the semi-universal deformation of \S\ref{ss:kleinian}, and $\chi\colon \cY_\Gamma \to \h$ its simultaneous resolution: $\chi$ is smooth, and there is a projective morphism $\phi\colon \cY_\Gamma \to \cX' = \cX_\Gamma \times_{\h/W} \h$ over $\h$ that is the minimal resolution on each fibre. Let $Y_0^J$ be the partial resolution of the Kleinian singularity that contracts the curves of $J$, and $\cY_J \to \h/W_J$ its semi-universal deformation, which is projective over the affine variety $\cX_J = \cX_\Gamma \times_{\h/W} \h/W_J$.

\begin{prop}[{Simultaneous resolutions; \cite[\S\S2.3, 4.3--4.4, Thm~3.3, Prop~4.15]{K}, \cite[\S\S5, 7]{CS}, \cite[Example~8.10]{SMM}, \cite[Prop~2.2]{BKL}}]\label{prop:simres}
In this notation, the following hold.
\begin{thmlist}
\item (The tilting bundle.) $\cY_\Gamma$ carries a tilting bundle $\mathcal T$ with $\cO$ as a summand, a projective generator of ${}^0\Per(\cY_\Gamma/\cX')$, whose endomorphism algebra is $\cA_\Gamma$ over $k[\h]$, with $e_0$ the idempotent of $\cO$.
\item (The blow-down.) Blowing down the curves of $J$ gives a proper morphism $\cY_\Gamma \to \cY_J \times_{\h/W_J} \h$ over $\h$, which on each fibre is the minimal resolution of the rational double points of the target.
\item (Root hyperplanes.) Near a general point of $\alpha^\perp$, $\cX' \to \h$ is, analytically and up to a unit, $x^2 + y^2 + z^2 = a^2$ with $a$ a local equation of $\alpha^\perp$, and for a general $\lambda \in \alpha^\perp$ the fibre $\chi^{-1}(\lambda)$ contains a single smooth $(-2)$-curve, of class $\alpha$ in $H_2$ of the fibre, which specialises to $H_2(\chi^{-1}(0), \Z) = Q_\Gamma$ with the exceptional curves giving the simple roots.
\item (Rational fibres.) $R\phi_*\cO_{\cY_\Gamma} = \cO$. Likewise $R\psi_*\cO = \cO$ for every partial resolution $\psi$ of a surface with rational double points by a surface with rational double points, through which the minimal resolution factors.
\end{thmlist}
\end{prop}

\begin{proof}
Part~(1) is the deformed McKay correspondence \cite[Prop~4.15]{K}, and we outline a proof. By \cite[\S\S5, 7]{CS}, $\cY_\Gamma$ is the moduli space over $\h$ of $\theta$-stable $\cA_\Gamma$-modules of dimension vector $\delta$, for a generic stability condition $\theta$. We take $\theta$ so that the stable modules are those generated at the vertex $0$, as in \cite[\S3.2]{K}. Then the tautological bundle $\mathcal T^\vee$, with $\cO$ as the summand of the vertex $0$, is generated by global sections, and the universal module gives a $k[\h]$-algebra map $\cA_\Gamma \to \End(\mathcal T)$. This map is graded for the $k^\times$-action of ``Gradings'' below, under which $\cO(\cX')$ is non-negatively graded with degree-$0$ part $k$, so graded Nakayama applies to finite graded $\cO(\cX')$-modules such as $\Ext^i(\mathcal T, \mathcal T)$. Over $0 \in \h$ the map is the isomorphism of the McKay correspondence, and $\mathcal T$ restricts to a projective generator of ${}^0\Per$ of the minimal resolution \cite[Prop~3.7, \S4.1]{K}. Base change (Lemma~\ref{lem:flatfamilies}(1)) and graded Nakayama then show, as in \cite[Lemma~5.3]{Kar18}, that $\Ext^i(\mathcal T, \mathcal T) = 0$ for $i > 0$ and that $\End(\mathcal T)$ is flat over $k[\h]$. By graded Nakayama again, the map is an isomorphism. Finally $H^1(\mathcal T) = \Ext^1(\cO, \mathcal T)$ is a summand of $\Ext^1(\mathcal T, \mathcal T) = 0$, as $\cO$ is a summand of $\mathcal T$. The line bundle $\det \mathcal T^\vee$ is ample on $\chi^{-1}(0)$, by Theorem~\ref{thm:VdB}(2) for the minimal resolution, so it is ample on the fibre of $\phi$ over the origin, and hence ample (the remark after Theorem~\ref{thm:VdB}). Thus $\mathcal T$ is a projective generator of ${}^0\Per$, and so a tilting bundle (Theorem~\ref{thm:VdB}(2)).

Of the other parts, only (4) needs an argument beyond the citations. For $\psi$ as in (4), the minimal resolution of its target is $\psi \circ \psi'$ for a resolution $\psi'$ of its source. Both surfaces have rational double points, which are rational singularities, so $R\psi'_*\cO = \cO$ and $R(\psi \circ \psi')_*\cO = \cO$, and hence $R\psi_*\cO = R\psi_*R\psi'_*\cO = \cO$. On each closed fibre $\phi$ is the minimal resolution of a surface with rational double points, so $R(\phi_\lambda)_*\cO = \cO$ for the fibre $\phi_\lambda$ of $\phi$ over $\lambda$. As $\cY_\Gamma$ and $\cX'$ are flat over $\h$, Lemma~\ref{lem:flatfamilies}(2) gives $R\phi_*\cO_{\cY_\Gamma} = \cO$.
\end{proof}

\paragraph{Gradings.} The contracting action of \S\ref{ss:kleinian} lifts to $\cY_\Gamma$ and $\cY_J$, compatibly with the maps between them, and induces, up to rescaling, the gradings of $\cA_\Gamma$, $A_J$ (Proposition~\ref{prop:tilting}) and $k[\h]$. Indeed, up to rescaling, the weights of $x$, $y$, $z$ are their degrees as invariants of $k[u, v]$ with $u$, $v$ in degree $1$, and those of the $u_j$ then follow from the homogeneity of $F$. Hence on $\cX' = \Spec e_0\cA_\Gamma e_0$ the action scales the arrows of $\tGamma$, with $\lambda$ in degree $2$ \cite[Lemma~8.2, Cor~8.12]{CBH}. The grading of $\cA_\Gamma$ in Proposition~\ref{prop:simres}(1) is a $k^\times$-action on $\cA_\Gamma$, hence on the moduli space $\cY_\Gamma$ of the proof above, and it linearises the tautological bundle $\mathcal T^\vee$ so that $\End(\mathcal T) = \cA_\Gamma$ carries that grading. This action restricts, up to rescaling, to the one above on $\cX'$, so it is the lift of that action to $\cY_\Gamma$, the lift being determined on the dense open set where $\phi$ is an isomorphism. We use this linearisation of $\mathcal T$. The grading of $A_J$ is defined by the linearisation of Step~1 in the proof of Proposition~\ref{prop:tilting}.

\subsection{The partial resolution and the corner}\label{ss:partial}

In this subsection we show that $\cY_J$ is smooth, which Proposition~\ref{prop:transfer} needs, prove at any vertex that it carries a tilting bundle whose endomorphism algebra base changes to a corner of $\cA_\Gamma$, and define the algebras $\Pi_\bul$ and $\cC_J$ that Sections~\ref{sec:cartan} and~\ref{sec:squeeze} study.

\begin{lem}\label{lem:Ysmooth}
$\cY_J$ is smooth.
\end{lem}

\begin{proof}
As $Y_0^J$ maps to an affine surface with fibres of dimension at most one, $H^2(Y_0^J, T_{Y_0^J}) = 0$ ($T$ the tangent sheaf), and $H^1$ of the sheaf $\mathscr{E}\!xt^1(\Omega_{Y_0^J}, \cO)$ of local first-order deformations vanishes, as it is supported at the finitely many singular points. Thus the map from the deformations of $Y_0^J$ to the product of the deformations of its singular points, which are rational double points, is smooth (deformations of the singular points lift to $Y_0^J$). Indeed, by the local-to-global spectral sequence for $\Ext^\bullet(\Omega_{Y_0^J}, \cO)$ its tangent map is onto, and the obstructions to lifting lie in $H^2(Y_0^J, T_{Y_0^J}) = 0$ (compare \cite[\S2]{BW74}, where the surface is projective; the argument does not use properness). Hence at a singular point of $Y_0^J$ the completed local ring of $\cY_J$ is a power series ring over that of the total space of a semi-universal deformation of a rational double point, which is smooth (it is a hypersurface linear in the constant-term parameter). At the smooth points of $Y_0^J$ the flat morphism $\cY_J \to \h/W_J$ has smooth fibre, so $\cY_J$ is smooth there too. Hence the singular locus of $\cY_J$ is closed, disjoint from $Y_0^J$ and stable under the contracting $k^\times$-action, lifted to $\cY_J$ as in ``Gradings'' above. It is empty: this action contracts the affine variety $\cX_J$ to a point, so by the valuative criterion of properness every $k^\times$-orbit in $\cY_J$ has a limit at $0$, which lies in $Y_0^J$. (See also \cite[Example~8.10]{SMM}.)
\end{proof}

Sections~\ref{sec:cartan} and~\ref{sec:squeeze} use the following proposition at every vertex. For the six flopping pairs it is \cite[Prop~4.21]{K}. We give a proof that works at every vertex, in which Step~1 constructs the $k^\times$-linearised projective generator that \cite[Rem~3.5]{K} starts from. Theorem~A uses only the six pairs, and the general case enters through Theorem~\ref{thm:Avertex} and Corollary~\ref{cor:centre}. Later sections use only the statement.

\begin{prop}[{The tilting bundle on $\cY_J$; extending \cite[Prop~4.21]{K}}]\label{prop:tilting}
$\cY_J$ carries a $k^\times$-linearised tilting bundle $\cO \oplus \mathcal M$ that is a projective generator of ${}^0\Per(\cY_J/\cX_J)$. Its endomorphism algebra $A_J$, graded by the linearisation, satisfies
\[
A_J \otimes_{H_J} k[\h] \;\cong\; e\cA_\Gamma e, \qquad e = e_0 + e_\bul,
\]
compatibly with $e_0$ and the gradings.
\end{prop}

\begin{proof}
We argue in three steps. We build a $k^\times$-linearised tilting bundle $\cO \oplus \mathcal M$ on $\cY_J$, pull $\mathcal M$ back to $\cY_\Gamma$, and identify the pullback with $\mathcal T_\bul$, the summand of vertex $\bul$ of the bundle $\mathcal T$ of Proposition~\ref{prop:simres}(1), by comparing degrees and ranks. Put $\cY_J' = \cY_J \times_{\h/W_J} \h$. Let $g_J\colon \cY_J \to \cX_J$ be the contraction and $\vartheta\colon \cY_\Gamma \to \cY_J'$ the blow-down of Proposition~\ref{prop:simres}(2), through which $\phi$ factors. These spaces and maps form the commutative diagram
\[
\begin{tikzpicture}[baseline=(current bounding box.center), x=2.3cm, y=1.25cm, ar/.style={-{Stealth[length=4pt]}}, lab/.style={font=\scriptsize}]
\node (YG) at (0,2) {$\cY_\Gamma$};
\node (YJp) at (1,2) {$\cY_J'$};
\node (YJ) at (2,2) {$\cY_J$};
\node (Xp) at (1,1) {$\cX'$};
\node (XJ) at (2,1) {$\cX_J$};
\node (XG) at (3,1) {$\cX_\Gamma$};
\node (h) at (1,0) {$\h$};
\node (hJ) at (2,0) {$\h/W_J$};
\node (hW) at (3,0) {$\h/W$};
\draw[ar] (YG) -- node[lab, above] {$\vartheta$} (YJp);
\draw[ar] (YJp) -- (YJ);
\draw[ar] (YG) -- node[lab, above right=-2pt] {$\phi$} (Xp);
\draw[ar] (YG) to[bend right=25] node[lab, below left=-2pt] {$\chi$} (h);
\draw[ar] (YJp) -- (Xp);
\draw[ar] (YJ) -- node[lab, right] {$g_J$} (XJ);
\draw[ar] (Xp) -- (XJ);
\draw[ar] (XJ) -- (XG);
\draw[ar] (Xp) -- (h);
\draw[ar] (XJ) -- (hJ);
\draw[ar] (XG) -- (hW);
\draw[ar] (h) -- node[lab, below] {$\pi$} (hJ);
\draw[ar] (hJ) -- (hW);
\end{tikzpicture}
\]
in which every square is cartesian. The maps $g_J$ and $\phi$ are $k^\times$-equivariant, as the action lifts to $\cY_J$ and to $\cY_\Gamma$. The map $\vartheta$ is equivariant too: over the smooth points of the fibres of $\cX' \to \h$ the base change $\cY_J' \to \cX'$ of $g_J$ is an isomorphism, so there $\vartheta$ is $\phi$ followed by its inverse. As the preimage of this set is dense in $\cY_\Gamma$ and $\cY_J'$ is separated, $\vartheta$ is equivariant everywhere. By Lemma~\ref{lem:flatfamilies}(2) and Proposition~\ref{prop:simres}(4), $Rg_{J*}\cO = \cO$ and $R\vartheta_*\cO = \cO$: $\cY_J$ and $\cX_J$ are flat over $\h/W_J$, $\cY_J'$ and $\cY_\Gamma$ are flat over $\h$, and on each closed fibre these maps are partial or minimal resolutions of rational double points. Hence Theorem~\ref{thm:VdB} applies to $g_J$, and by Proposition~\ref{prop:simres}(4) also to $\phi$.

\emph{Step 1 (a graded tilting bundle on $\cY_J$).} Following Remark~3.5 of \cite{K}, we first build a $k^\times$-linearised projective generator of ${}^0\Per(\cY_J/\cX_J)$, and then cut it down to one whose completion at the origin is basic. As $g_J$ is projective and $\cX_J$ is affine, $\cY_J$ has an ample line bundle $\mathcal L$ generated by global sections. Since $\cY_J$ is smooth (Lemma~\ref{lem:Ysmooth}), hence normal, and $\operatorname{Pic}(k^\times) = 0$, the line bundle $\mathcal L$ admits a $k^\times$-linearisation \cite[Thm~2.14]{Bri}, which we fix. Then homogeneous generators of the graded $\cO(\cX_J)$-module $H^0(\mathcal L)$ give a $k^\times$-equivariant exact sequence $0 \to \mathcal K \to \cO^{\oplus a} \to \mathcal L \to 0$, with each copy of $\cO$ linearised by a character. As $H^0(\cO^{\oplus a}) \to H^0(\mathcal L)$ is onto and $H^1(\cO) = 0$, we have $H^1(\mathcal K) = 0$. The dual sequence shows that $\mathcal K^\vee$ is generated by global sections and that $\det \mathcal K^\vee \cong \mathcal L$ is ample. Hence $\cO \oplus \mathcal K$ is a $k^\times$-linearised projective generator of ${}^0\Per(\cY_J/\cX_J)$ (Theorem~\ref{thm:VdB}(2)). The action of $k^\times$ contracts $\cX_J$ to its origin, over which $\cY_J$ contains a single curve $C$, the image of the exceptional curve $E_\bul$ of $\bul$ in the minimal resolution $\chi^{-1}(0)$.

We decompose $\cO \oplus \mathcal K$ into indecomposable linearised summands, each of which completes at the origin to $\cO$ or to $N_1$ (Theorem~\ref{thm:VdB}(3)), and take for $\mathcal M$ one that completes to $N_1$. One exists, since $\det \mathcal K^\vee \cong \mathcal L$ has positive degree on $C$, while a summand that completes to $\cO$ has determinant of degree $0$ there. Two facts justify the decomposition. First, $E_{\mathcal K} = \End(\cO \oplus \mathcal K)$ is a graded module finite over $\cO(\cX_J)$, which is positively graded with degree-$0$ part $k$, so each graded piece of $E_{\mathcal K}$ is finite-dimensional, although $E_{\mathcal K}$ may have pieces of negative degree. In particular the algebra $(E_{\mathcal K})_0$ of linearised endomorphisms is finite-dimensional, so $\cO \oplus \mathcal K$ is a direct sum of linearised summands $\mathcal P$ with $\End(\mathcal P)_0$ local. Second, the completion of such a $\mathcal P$ at the origin is indecomposable, hence a copy of $\cO$ or of $N_1$ by Theorem~\ref{thm:VdB}(3). Indeed, by flat base change its endomorphism algebra is the completion of $E_{\mathcal P} = \End(\mathcal P)$ at the maximal ideal $\mathfrak m$ of the origin, whose idempotents lift from $\bar E = E_{\mathcal P}/\mathfrak m E_{\mathcal P}$, so it suffices that $\bar E$ is local. The radical of $\bar E$ is stable under $k^\times$, hence graded, so $\bar E/\operatorname{rad}\bar E$ is a graded semisimple algebra whose degree-$0$ part is local, as a quotient of $\End(\mathcal P)_0$. The connected group $k^\times$ preserves each simple factor $M_n(k)$, acting on it through $PGL_n(k)$ by the Skolem--Noether theorem; on a finite cover of $k^\times$ this lifts to a torus in $GL_n(k)$, and diagonalising the torus gives $n$ orthogonal idempotents of degree $0$ in that factor. As the degree-$0$ part is local, $\bar E/\operatorname{rad}\bar E = k$, so $\bar E$ is local.

It remains to see that $\cO \oplus \mathcal M$ is a projective generator. It is projective, as a summand of $\cO \oplus \mathcal K$. The line bundle $\det \mathcal M^\vee$ has degree $1$ on $C$ (Theorem~\ref{thm:VdB}(3)), so it is ample on the fibre over the origin, and hence ample by the remark after Theorem~\ref{thm:VdB}. Thus $\cO \oplus \mathcal M$ is a projective generator of ${}^0\Per(\cY_J/\cX_J)$, and in particular a tilting bundle (Theorem~\ref{thm:VdB}(2)). Put $A_J = \End(\cO \oplus \mathcal M)$, graded by the linearisation, which we normalise at the end of Step~3.

\emph{Step 2 (pulling back to $\cY_\Gamma$).} Let $\mathcal M'$ be the pullback of $\mathcal M$ to $\cY_J'$. As $\cO \oplus \mathcal M$ is projective in ${}^0\Per(\cY_J/\cX_J)$, the bundle $\mathcal M^\vee$ is generated by global sections and $H^1(\mathcal M) = 0$ (Theorem~\ref{thm:VdB}(1)). Since $k[\h]$ is flat over $H_J$, the same holds for $\mathcal M'$, and $\End(\cO \oplus \mathcal M') \cong A_J \otimes_{H_J} k[\h]$. By the projection formula and $R\vartheta_*\cO = \cO$, these properties pass to $\vartheta^*\mathcal M'$: $\End(\cO \oplus \vartheta^*\mathcal M') \cong \End(\cO \oplus \mathcal M')$, $H^1(\vartheta^*\mathcal M') = 0$, and $\vartheta^*\mathcal M'^\vee$ is generated by global sections. Thus $\vartheta^*\mathcal M'$ is a $k^\times$-linearised projective object of ${}^0\Per(\cY_\Gamma/\cX')$, again by Theorem~\ref{thm:VdB}(1).

\emph{Step 3 (identifying the pullback).} Write $\mathcal T = \bigoplus_{i \in \tGamma_0} \mathcal T_i$ for the bundle of Proposition~\ref{prop:simres}(1), with $\mathcal T_0 = \cO$ and $\mathcal T_i$ the summand of vertex $i$. By Proposition~\ref{prop:simres}(1), $\mathcal T$ is a projective generator of ${}^0\Per(\cY_\Gamma/\cX')$, $k^\times$-linearised as in ``Gradings'' (\S\ref{ss:simres}), and it is basic, as the degree-$0$ part of $\cA_\Gamma$ is spanned by the vertex idempotents. Hence, by Theorem~\ref{thm:VdB}(3), the summands $\mathcal T_i$ for $i \in \Gamma_0$ complete, at the origin of $\cX'$, to the bundles $N_j$ of Theorem~\ref{thm:VdB}(3) for $\chi^{-1}(0)$, one for each exceptional curve $E_j$ of $\chi^{-1}(0)$. Proposition~\ref{prop:simres}(3) decides which summand goes with which curve: $\mathcal T_i$ completes to $N_i$, where $E_i$ is the curve whose class is the simple root at $i$.

To see this, let $E$ be the curve of Proposition~\ref{prop:simres}(3) whose class is the simple root at $i$, in the fibre over a general point $\lambda$ of the hyperplane $t_i = 0$, so that each $\det \mathcal T_k$ has the same degree on $E$ as on $E_i$. As each $\mathcal T_k^\vee$ is generated by global sections, $\Hom(\mathcal T, \cO_E(-1))$ is an $\cA_\Gamma$-module, killed by the ideal of $\lambda$, of dimension vector $d = (-\deg \det \mathcal T_k|_{E_i})_k$ (each $\mathcal T_k^\vee|_E$ is a sum of line bundles $\cO(a)$ with $a \ge 0$, so $\Hom(\mathcal T_k, \cO_E(-1)) = H^0(\mathcal T_k^\vee|_E(-1))$ has dimension $\deg \mathcal T_k^\vee|_E$). Taking traces in the preprojective relation gives $\sum_k \lambda_k d_k = 0$ for all such $\lambda$, and $d_0 = 0$, so $d$ is supported at~$i$. By Theorem~\ref{thm:VdB}(3), $\det N_j$ has degree $-\delta_{ij}$ on $E_i$, so $\mathcal T_i$ is the summand that completes to $N_i$.

It remains to identify $\vartheta^*\mathcal M'$ with $\mathcal T_\bul$, first up to multiplicities and then by ranks. Since $\mathcal T$ is a projective generator, the projective object $\vartheta^*\mathcal M'$ of Step~2 corresponds to a finitely generated graded projective $\cA_\Gamma$-module. As $\cA_\Gamma$ is non-negatively graded, with degree-$0$ part spanned by the vertex idempotents, graded Nakayama gives $\vartheta^*\mathcal M' \cong \bigoplus_{i \in \tGamma_0} \mathcal T_i^{\oplus m_i}$, with each summand's linearisation twisted by a character of $k^\times$. The restriction of $\det \mathcal T_i$ to $\chi^{-1}(0)$ has degree $-\delta_{ij}$ on $E_j$ for $j \in \Gamma_0$, and likewise $\det \mathcal M$ has degree $-1$ on $C$ (Theorem~\ref{thm:VdB}(3)). Now $\vartheta$ contracts $E_j$ for $j \in J$ and maps $E_\bul$ birationally onto $C$, so $\det \vartheta^*\mathcal M'$ has degree $-\delta_{\bul j}$ on $E_j$. Since $\det \bigoplus_i \mathcal T_i^{\oplus m_i}$ has degree $-m_j$ on $E_j$, it follows that $m_j = \delta_{\bul j}$ for $j \in \Gamma_0$.

Next, $\rank \mathcal T_\bul$ and $\rank \mathcal M$ are the multiplicities of $E_\bul$ and of $C$ in the fibres over the origins of $\cX'$ and of $\cX_J$ (Theorem~\ref{thm:VdB}(3)). These agree, as the fibre of $\cY_J'$ over the origin of $\cX'$ is that of $g_J$, and $\vartheta$ is an isomorphism near the generic point of $E_\bul$. As $\rank \mathcal M = \rank \vartheta^*\mathcal M' = m_0 + \rank \mathcal T_\bul$, it follows that $m_0 = 0$. Hence, after twisting the linearisation of $\mathcal M$ by a character, $\vartheta^*\mathcal M' \cong \mathcal T_\bul$. By Step~2, $A_J \otimes_{H_J} k[\h] \cong \End(\cO \oplus \mathcal T_\bul) = e\cA_\Gamma e$, compatibly with $e_0$ and the gradings.
\end{proof}

\paragraph{The corner.} With $\Pi_\Gamma$ the family of \S\ref{ss:dpa} over $k[\h]$, put
\[
\Pi_\bul = e_\bul\Pi_\Gamma e_\bul, \qquad \cC_J = A_J / A_J e_0 A_J ,
\]
an algebra over $k[\h]$ and an algebra over $H_J$. Sections~\ref{sec:cartan} and~\ref{sec:squeeze} study them through the following identification.

\begin{lem}\label{lem:corner}
$\cA_\Gamma/\cA_\Gamma e_0\cA_\Gamma \cong \Pi_\Gamma$, and $\cC_J \otimes_{H_J} k[\h] \cong \Pi_\bul$ as graded $k[\h]$-algebras; equivalently, as $H_J = P[t]$ and $k[\h] = k[\h_J][t]$, $\Pi_\bul \cong \cC_J \otimes_P k[\h_J]$.
\end{lem}

\begin{proof}
Killing $e_0$ in $\cA_\Gamma = \Pi^\lambda(\tGamma)$ kills the arrows incident to $0$ and the relation at $0$, and each other relation loses exactly its terms through $0$, leaving the preprojective relation of $\Gamma$ at that vertex; the constraint $\lambda\cdot\delta = 0$ only determines $\lambda_0$, which disappears with vertex $0$. This gives the first isomorphism. Tensoring is right exact, so by Proposition~\ref{prop:tilting}, $\cC_J \otimes_{H_J} k[\h] \cong e\cA_\Gamma e/\langle e_0\rangle$, where $e = e_0 + e_\bul$. Since $e e_0 = e_0 = e_0 e$, the ideal of $e\cA_\Gamma e$ generated by $e_0$ is $e(\cA_\Gamma e_0\cA_\Gamma)e$, so $\cC_J \otimes_{H_J} k[\h] \cong e\Pi_\Gamma e$, and the image of $e$ in $\Pi_\Gamma$ is $e_\bul$. The last isomorphism follows, as $k[\h] = H_J \otimes_P k[\h_J]$.
\end{proof}

\subsection{Universal flops}\label{ss:universal}

\begin{figure}[ht]
\centering
\begin{tikzpicture}[scale=0.62,
  v/.style={circle, draw, fill=white, inner sep=0pt, minimum size=5pt},
  f/.style={circle, draw, fill=black, inner sep=0pt, minimum size=6pt},
  c/.style={font=\scriptsize, above=1pt},
  a/.style={font=\scriptsize}]
\begin{scope}[shift={(0,0)}]
  \node[f] (b) at (0,0) {}; \node[c] at (b) {$1$};
  \node[font=\small] at (0,-1.1) {$l = 1$, $\mathsf{A}_1$};
\end{scope}
\begin{scope}[shift={(4,0)}]
  \node[v] (x1) at (-1,0) {}; \node[f] (b) at (0,0) {}; \node[v] (x3) at (1,0) {}; \node[v] (x4) at (0,-0.9) {};
  \draw (x1)--(b)--(x3) (b)--(x4);
  \node[c] at (x1) {$1$}; \node[c] at (b) {$2$}; \node[c] at (x3) {$1$}; \node[font=\scriptsize, right=1pt] at (x4) {$1$};
  \node[a] at (-1,-0.45) {$\beta$}; \node[a] at (1,-0.45) {$\gamma$}; \node[a] at (-0.4,-0.9) {$\delta$};
  \node[font=\small] at (0,-2.1) {$l = 2$, $\mathsf{D}_4$};
\end{scope}
\begin{scope}[shift={(10,0)}]
  \foreach \i/\x/\m in {1/-2/1, 3/-1/2, 4/0/3, 5/1/2, 6/2/1} \node[v] (x\i) at (\x,0) {};
  \node[f] (x4) at (0,0) {}; \node[v] (x2) at (0,-0.9) {};
  \draw (x1)--(x3)--(x4)--(x5)--(x6) (x4)--(x2);
  \foreach \i/\m in {1/1, 3/2, 4/3, 5/2, 6/1} \node[c] at (x\i) {$\m$};
  \node[font=\scriptsize, right=1pt] at (x2) {$2$};
  \node[a] at (-1.5,-0.45) {$\beta$}; \node[a] at (1.5,-0.45) {$\gamma$}; \node[a] at (-0.4,-0.9) {$\delta$};
  \node[font=\small] at (0,-2.1) {$l = 3$, $\mathsf{E}_6$};
\end{scope}
\begin{scope}[shift={(-0.5,-4)}]
  \foreach \i/\x in {1/-2, 3/-1, 5/1, 6/2, 7/3} \node[v] (x\i) at (\x,0) {};
  \node[f] (x4) at (0,0) {}; \node[v] (x2) at (0,-0.9) {};
  \draw (x1)--(x3)--(x4)--(x5)--(x6)--(x7) (x4)--(x2);
  \foreach \i/\m in {1/2, 3/3, 4/4, 5/3, 6/2, 7/1} \node[c] at (x\i) {$\m$};
  \node[font=\scriptsize, right=1pt] at (x2) {$2$};
  \node[a] at (-1.5,-0.45) {$\delta$}; \node[a] at (2,-0.45) {$\gamma$}; \node[a] at (-0.4,-0.9) {$\beta$};
  \node[font=\small] at (0.5,-2.1) {$l = 4$, $\mathsf{E}_7$};
\end{scope}
\begin{scope}[shift={(7,-4)}]
  \foreach \i/\x in {1/-2, 3/-1, 4/0, 6/2, 7/3, 8/4} \node[v] (x\i) at (\x,0) {};
  \node[f] (x5) at (1,0) {}; \node[v] (x2) at (0,-0.9) {};
  \draw (x1)--(x3)--(x4)--(x5)--(x6)--(x7)--(x8) (x4)--(x2);
  \foreach \i/\m in {1/2, 3/4, 4/6, 5/5, 6/4, 7/3, 8/2} \node[c] at (x\i) {$\m$};
  \node[font=\scriptsize, right=1pt] at (x2) {$3$};
  \node[font=\small] at (1,-2.1) {$l = 5$, $\mathsf{E}_8$};
\end{scope}
\begin{scope}[shift={(15,-4)}]
  \foreach \i/\x in {1/-2, 3/-1, 5/1, 6/2, 7/3, 8/4} \node[v] (x\i) at (\x,0) {};
  \node[f] (x4) at (0,0) {}; \node[v] (x2) at (0,-0.9) {};
  \draw (x1)--(x3)--(x4)--(x5)--(x6)--(x7)--(x8) (x4)--(x2);
  \foreach \i/\m in {1/2, 3/4, 4/6, 5/5, 6/4, 7/3, 8/2} \node[c] at (x\i) {$\m$};
  \node[font=\scriptsize, right=1pt] at (x2) {$3$};
  \node[a] at (-1.5,-0.45) {$\gamma$}; \node[a] at (2.5,-0.45) {$\delta$}; \node[a] at (-0.4,-0.9) {$\beta$};
  \node[font=\small] at (1,-2.1) {$l = 6$, $\mathsf{E}_8$};
\end{scope}
\end{tikzpicture}
\caption{The six pairs $(\Gamma, \bul)$, labelled by the coefficients of the highest root; $\bul$ (filled) is the vertex of coefficient $l$, and the other vertices form $J$. When $\bul$ is trivalent, the arms of $J$ are labelled by the loops $\beta$, $\gamma$, $\delta$ at $\bul$ in the presentation of \cite{K} (as arm labels these letters name loops, not roots); in length $5$, $\bul$ is the vertex labelled $5$, not the trivalent one, and $J$ has type $\mathsf{A}_4 \times \mathsf{A}_3$.}\label{fig:pairs}
\end{figure}
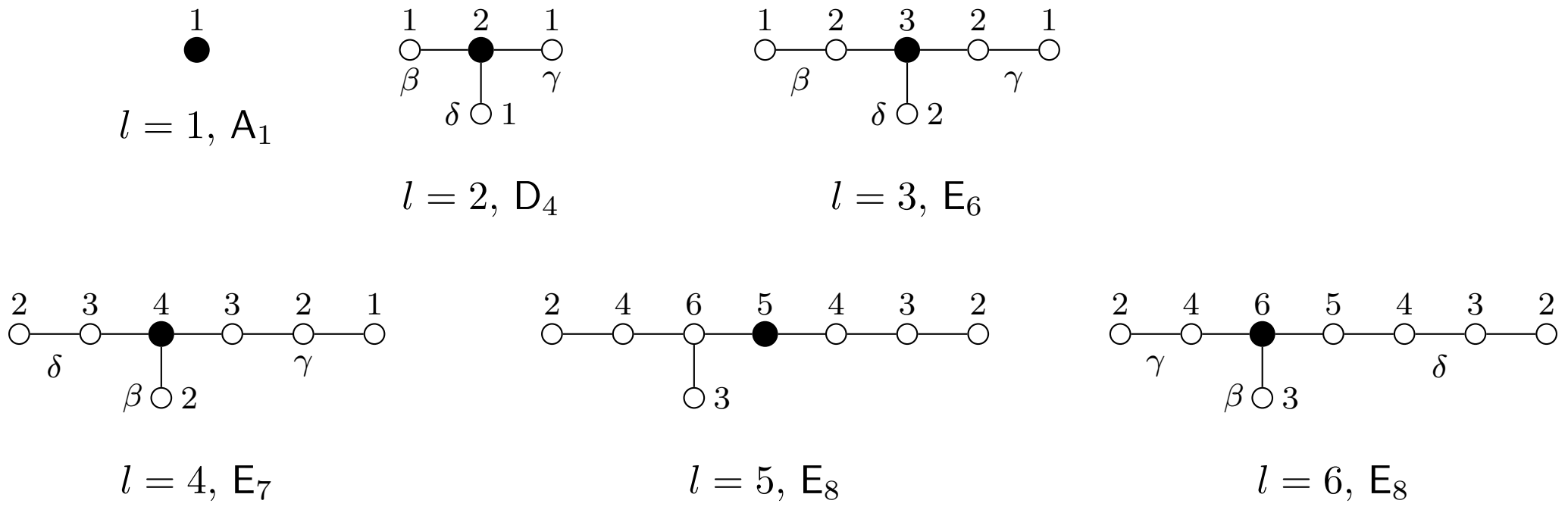

Finally we recall the universal flops, from which every flop is pulled back, and the curves in their fibres. By a theorem of Katz and Morrison \cite[\S1, Main Theorem and Corollary]{KM}, which Kawamata reproved without computer calculations \cite[Theorem]{Kaw94}, the generic hyperplane section of a flop of length $l$ determines a Dynkin diagram $\Gamma$ with a vertex $\bul$ of coefficient $l$ in the highest root, one pair $(\Gamma, \bul)$ for each length (Table~\ref{tab:roots} and Figure~\ref{fig:pairs}). With $J$ and $W_J$ as in \S\ref{ss:lie}, there is a universal flop of length $l$ over the base $\Spec H_l$, with $H_l = H_J = k[\h]^{W_J}$.

\begin{table}[ht]
\caption{Root counts for the six pairs $(\Gamma, \bul)$.}\label{tab:roots}
\begin{tabular}{@{}cclcc@{}}
\toprule
Length & $\Gamma$ & $\#\{\alpha \in \Phi^+ : \alpha_\bul = j\}$, $j = 1, \dots, l$ & $\sum_{\Phi^+} \alpha_\bul^2$ & $|\Phi^+_\bul|$ \\
\midrule
1 & $\mathsf{A}_1$ & $1$ & $1$ & $1$ \\
2 & $\mathsf{D}_4$ & $8, 1$ & $12$ & $9$ \\
3 & $\mathsf{E}_6$ & $18, 9, 2$ & $72$ & $29$ \\
4 & $\mathsf{E}_7$ & $24, 18, 8, 3$ & $216$ & $53$ \\
5 & $\mathsf{E}_8$ & $40, 30, 20, 10, 4$ & $600$ & $104$ \\
6 & $\mathsf{E}_8$ & $30, 30, 20, 15, 6, 5$ & $900$ & $106$ \\
\bottomrule
\end{tabular}
\end{table}

\paragraph{Versality.} In more detail, let $Y_0$ be the proper transform of a generic hyperplane section $X_0$ of a flop through the singular point, a partial resolution of the Kleinian singularity of type $\Gamma$ contracting the curves of $J$. By Katz and Morrison, the semi-universal deformations of $Y_0$ and $X_0$ have bases $\Spec H_l = \h/W_J$ and $\h/W$, and the universal flop is the semi-universal deformation of $Y_0$ with its contraction. These statements are made for analytic germs in \cite{KM}, and for the partial resolution $Y_0^J$ of \S\ref{ss:simres} in \cite{K} (see also \cite{BKL}). They pass to $Y_0$, the pullback of $Y_0^J$ to the completion $X_0$ of the Kleinian singularity. Restriction from $Y_0^J$ to $Y_0$ is an isomorphism on first-order deformations, which are finite-dimensional and supported on the exceptional curves. Both deformation functors are unobstructed, since the singularities are hypersurfaces and $H^2(T) = 0$. A map of unobstructed deformation functors that is onto on first-order deformations is smooth, as a lift chosen in the source can be corrected by a first-order deformation to match a given lift in the target, so the semi-universal deformation of $Y_0^J$ restricts to one of $Y_0$.

We write $g\colon \cY_l \to \cX_l = \Spec R_l$ for this universal flopping contraction; thus $\cY_l = \cY_J$ (\S\ref{ss:simres}), which is smooth (Lemma~\ref{lem:Ysmooth}). Both $\cY_l$ and $\cX_l$ are flat over $H_l$, as deformations of $Y_0$ and $X_0$ over $\Spec H_l$, and $Rg_*\cO_{\cY_l} = \cO_{\cX_l}$ by Lemma~\ref{lem:flatfamilies}(2), as the fibres of $g$ over closed points are partial resolutions of rational double points (Proposition~\ref{prop:simres}(4)).

We use versality in the following form, in Lemma~\ref{lem:arcs} and Theorem~B(3).

\begin{lem}[Classifying maps]\label{lem:classifying}
Let $B$ be a complete local noetherian $k$-algebra with residue field $k$, and $f_B\colon \mathcal Y \to \mathcal X = \Spec R$ a proper morphism over $B$ with $Rf_{B*}\cO_{\mathcal Y} = \cO_{\mathcal X}$, where $R$ is a complete local noetherian $B$-algebra and $\mathcal Y$ is flat over $B$ with closed fibre $Y_0$. Then there is a map $\Spec B \to \Spec H_l$, sending the closed point to the origin, along which $\mathcal Y \to \mathcal X$ is the pullback of the universal flop, completed at the singular point of its closed fibre.
\end{lem}

\begin{proof}
As the semi-universal deformation is versal, one can choose maps $\Spec B/\mathfrak m_B^n \to \Spec H_l$ and isomorphisms of the deformations over them inductively, each extending the one before. Together they give an isomorphism of the formal completions along $Y_0$, the $\mathfrak m_B$-adic completions. This isomorphism lies over $\mathcal X$, because $\mathcal X$ is the affinisation $\Spec \Gamma(\mathcal Y, \cO_{\mathcal Y})$ of $\mathcal Y$ (as $Rf_{B*}\cO_{\mathcal Y} = \cO_{\mathcal X}$) and $\cX_l$ is that of $\cY_l$. As $\mathcal Y$ and the completed pullback of $\cY_l$ are proper over $\mathcal X$, and $R$ is $\mathfrak m_B R$-adically complete since $\mathfrak m_B R \subseteq \mathfrak m_R$, Grothendieck's existence theorem makes it an isomorphism of schemes over $\mathcal X$.
\end{proof}

\paragraph{Generators of $H_l$.} The group $W_J$ acts on $\mathrm{span}\{t_i : i \in J\}$ as the reflection representation of a product of symmetric groups, one for each arm of $J$, that is, each connected component of $J$. An arm of type $\mathsf{A}_{m-1}$ is said to have \emph{size} $m$, and its coordinates are written $\tau_1, \dots, \tau_m$ with $\sum \tau_k = 0$, permuted by the symmetric group. Hence $P$ is generated by the elementary symmetric polynomials $\sigma_2, \dots, \sigma_m$ of each arm, of degrees $2, \dots, m$ in units of the degree of $t_i$, with $\prod m! = |W_J|$, and $H_l = P[t]$ is generated in degrees less than the Coxeter number $h$. Up to sign, the $\sigma_i$ are the coefficients of the \emph{arm polynomial} $\prod_k (x + \tau_k)$ and the parameters $T$ of the presentation of \cite{K}, which Remark~\ref{rem:Kpresentation} makes explicit.

\paragraph{Curves in the fibres.} Bryan, Katz and Leung describe the curves in the fibres of $\cY_l$, building on Katz and Morrison. Let $D^{\mathrm{curv}} \subset \Spec H_l$ be the set of points whose fibre in $\cY_l$ contains a complete curve.

\begin{prop}[{Curves in the fibres; \cite[Prop~2.3, Lemmas~2.4, 2.8]{BKL}}]\label{prop:BKL}
Let $D^{\mathrm{curv}}$ be as above.
\begin{thmlist}
\item (The discriminant.) $D^{\mathrm{curv}} = D_1 \cup \dots \cup D_l$, and $D_j$ is the image of the hyperplanes $\alpha^\perp$ with $\alpha_\bul = j$, which form a single $W_J$-orbit. Over a general point of $D_j$ the fibre contains a single complete curve, a smooth $(-2)$-curve of class $j[C]$.
\item (Generic deformations.) Let $\nu$ be an arc through the origin whose pullback is a flop and is the completion of an algebraic one-parameter deformation of $Y_0^J$ over a disc, and let $\nu'$ be a small deformation of $\nu$ transverse to $D^{\mathrm{curv}}$. Then every complete curve in the pullback of $\cY_l$ along $\nu'$ is a smooth $(-1,-1)$-curve; over each point near the origin where $\nu'$ meets $D_j$, and no other $D_i$, there is exactly one, of class $j[C]$; and the number of those of class $j[C]$ is the multiplicity $k_j$ of $C_j$ in its Hilbert scheme (\S\ref{ss:flops}).
\end{thmlist}
\end{prop}

We write $D = D_1 \cup \dots \cup D_l$, so that $D^{\mathrm{curv}} = D$ by (1). The description of the curves in (2) is the paragraph after Lemma~2.7 of \cite{BKL}, and the count is their Lemma~2.8, stated for a small deformation $\nu'$ in an analytic family of arcs deforming $\nu$, such as $\nu + s\xi$ in Proposition~\ref{prop:GV}. Both use the local model of Proposition~\ref{prop:simres}(3) near a general point of each $\alpha^\perp$. These statements are made in \cite{BKL} for analytic one-parameter deformations, with $C_j$ defined through a generic hyperplane section; the proof of Lemma~\ref{lem:truncation} shows that the section $w = 0$ of a classifying arc may be used instead.

The universal flop carries the tilting bundle $\cO \oplus \mathcal M$ of Proposition~\ref{prop:tilting}, and its endomorphism algebra $A_l = A_J$ is the universal flopping algebra of \cite[\S\S3.3, 4.4]{K}. Indeed, \cite{K} defines that algebra through a $k^\times$-linearised projective generator of ${}^0\Per(\cY_l/\cX_l)$ that completes at the origin to $\cO \oplus N_1$ (its Remarks~3.5 and~3.11), and $\cO \oplus \mathcal M$ is one (Step~1 of the proof of Proposition~\ref{prop:tilting}). Such a generator is unique up to twisting the linearisations of its summands by characters. This follows by graded Nakayama, as in Step~3 of that proof, since $A_J$ is a graded subalgebra of $A_J \otimes_{H_J} k[\h] \cong e\cA_\Gamma e$, with $e = e_0 + e_\bul$, and so is non-negatively graded with degree-$0$ part spanned by $e_0$ and $e_\bul$. By Proposition~\ref{prop:tilting}, $A_l \otimes_{H_l} k[\h] \cong e\cA_\Gamma e$, and \cite[Thm~1.3]{K} presents $A_l$ by a two-vertex quiver with relations over $H_l$.

The universal contraction algebra is
\[
\cC_l = A_l / A_l e_0 A_l = \cC_J ,
\]
so $\cC_l \otimes_{H_l} k[\h] \cong \Pi_\bul$ by Lemma~\ref{lem:corner}. The universal algebras $A_l$ and $\cC_l$ are graded, with $H_l$ positively graded. Kawamata used the universal flops of \cite{CM} and \cite{K} to describe the noncommutative deformations of the reduced fibre, which give the completion of $\cC_l$ at the origin \cite[Thm~6.2]{Kaw}.

\paragraph{Notation.} Figure~\ref{fig:story} shows how the main objects are related. A subscript $J$ marks an object defined at any vertex $\bul$ (Sections~\ref{sec:cartan} and~\ref{sec:squeeze}), and a subscript $l$ an object of the universal flop of length $l$; for the six pairs the two agree.
{
\small
\renewcommand{\arraystretch}{1.1}
\begin{xltabular}{\textwidth}{@{}>{\raggedright}p{2.6cm} X l@{}}
\multicolumn{3}{@{}l}{\emph{Algebras}} \\*
$\cA_\Gamma$ & $\Pi^\lambda(\tGamma)$ with $\lambda\cdot\delta = 0$, over $k[\h]$ & \S\ref{ss:simres} \\
$\Pi_\Gamma$, $\Pi_\bul$ & $\Pi^\lambda(\Gamma)$ over $k[\h]$, and its corner $e_\bul\Pi_\Gamma e_\bul$ & \S\S\ref{ss:dpa}, \ref{ss:partial} \\
$A_J$ & endomorphism algebra of the tilting bundle on $\cY_J$; $A_J \otimes_{H_J} k[\h] \cong e\cA_\Gamma e$ & Prop.~\ref{prop:tilting} \\
$\cC_J$ & $A_J/A_Je_0A_J$; $\cC_J \otimes_{H_J} k[\h] \cong \Pi_\bul$ & Lemma~\ref{lem:corner} \\
$A_l$, $\cC_l$ & the universal flopping and contraction algebras, $A_J$ and $\cC_J$ for the six pairs & \S\ref{ss:universal} \\
$\cC_\nu$ & $\cC_l \otimes_{H_l} B$, the pullback along $\nu\colon \Spec B \to \Spec H_l$, such as a classifying arc & Lemma~\ref{lem:basechange} \\[3pt]
\multicolumn{3}{@{}l}{\emph{Bases and roots}} \\*
$H_J$, $P$ & $H_J = k[\h]^{W_J} = k[t] \otimes P$, with $P = k[\h_J]^{W_J}$ & \S\ref{ss:lie} \\
$\alpha_\bul$, $\Phi^+_\bul$ & the coefficient of $\alpha$ at $\bul$; the positive roots with $\alpha_\bul \neq 0$ & \S\ref{ss:lie} \\
$f_j$, $D_j$ & $\prod_{\alpha_\bul = j}\alpha$ and $\{f_j = 0\} \subset \Spec H_J$ & \S\ref{ss:lie} \\
$\pi$, $\Sigma$ & the quotient $\h \to \Spec H_J$; the union of the pairwise intersections of root hyperplanes & \S\ref{ss:lie} \\
$\eta_\alpha$, $\ell_\alpha$, $q_\alpha$ & the generic point of $\alpha^\perp$; the lengths of $\Pi_\bul$ and $\HH_0(\Pi_\bul)$ there & \S\ref{ss:lie}, Cor.~\ref{cor:lower} \\
$\tau_k$, $T$ & coordinates on an arm of $J$; the parameters of \cite{K} & Rem.~\ref{rem:Kpresentation} \\[3pt]
\multicolumn{3}{@{}l}{\emph{Geometry}} \\*
$C$, $J$ & the flopping curve; the set of contracted vertices & \S\S\ref{ss:flops}, \ref{ss:lie} \\
$Y_0^J$, $Y_0$, $\chi^{-1}(0)$ & the partial resolution contracting the curves of $J$, its completion, and the minimal resolution & \S\S\ref{ss:simres}, \ref{ss:universal} \\
$g\colon \cY_l \to \cX_l$ & the universal flopping contraction & \S\ref{ss:universal} \\[3pt]
\multicolumn{3}{@{}l}{\emph{Flops}} \\*
$n_j$, $\SGV$ & the GV invariants, and $\SGV = \prod_j M_j(k)^{n_j}$ & \S\ref{ss:flops}, Section~\ref{sec:intro} \\
$\nu$, $\tnu$ & the classifying arc and its lift to $\h$ & Lemmas~\ref{lem:arcs}, \ref{lem:lift} \\
$\Ass_d$ & the scheme of $d$-dimensional algebra structures & \S\ref{ss:moduli} \\
\end{xltabular}
}

\section{Reduction to a centre count}\label{sec:reduction}

In this section we show that smoothness of a contraction algebra $A$, and uniqueness of its semisimple deformation, follow from its deformation to $\SGV$ and one numerical statement: $\dim Z(A) = n_1 + \dots + n_l$. Jelisiejew and Smoktunowicz observe that, once $A$ deforms to a semisimple algebra $S$ with $\dim Z(A) = \dim Z(S)$, as in types A and D, smoothness of $[A]$ is equivalent to $h^1(A) = h^2(A)$ \cite[\S1, Thm~3.7]{JS}. The arguments of this section come from \cite[Lemmas~2.8, 2.9, Prop~2.11]{Pan}, where they are stated for symmetric algebras with $\Omega^4(A) \cong A$ and for algebras that deform to a semisimple algebra, and applied to the Jacobi algebras of the $\mathsf{E}_6$ family. Lemma~\ref{lem:tangent} shows that the tangent space at $[A]$ exceeds the orbit $G\cdot[S]$ in dimension by $(\dim Z(A) - \dim Z(S)) + (h^2(A) - h^1(A))$, and Lemma~\ref{lem:pan} proves $h^1 = h^2$ for every contraction algebra, through \cite{JKM}, so that the centres alone decide smoothness. Proposition~\ref{prop:criterion} combines them, and adds the inequality $\dim Z(A) \ge \dim Z(S)$, the converse (smoothness forces equality), the unique component and the constancy of $\dim Z$. Hence Theorem~C reduces to Theorem~B(1) and~(2), which we prove in Section~\ref{sec:flops}.

\begin{lem}\label{lem:tangent}
Let $A$ and $S$ be $d$-dimensional algebras with $S$ semisimple. Then
\[
\dim T_{[A]}\Ass_d - \dim (G\cdot[S]) = \big(\dim Z(A) - \dim Z(S)\big) + \big(h^2(A) - h^1(A)\big).
\]
\end{lem}

\begin{proof}
By \S\ref{ss:moduli}, $\dim T_{[A]}\Ass_d = \dim (G\cdot[A]) + h^2(A)$ and $\dim (G\cdot[A]) = \dim G - \dim \Der(A) = \dim G - d + \dim Z(A) - h^1(A)$. Since $h^1(S) = 0$, $\dim (G\cdot[S]) = \dim G - d + \dim Z(S)$. Subtract.
\end{proof}

\begin{lem}[{$h^1 = h^2$; \cite[Lemma~2.8]{Pan}, \cite[Thm~3.2.5]{JKM}}]\label{lem:pan}
Let $A$ be a finite-dimensional symmetric algebra with $\Omega^4_{A^e}(A) \cong A$ in the stable category of bimodules. Then $h^1(A) = h^2(A)$. In particular this holds for every contraction algebra (Proposition~\ref{prop:Acon}(2), (3)).
\end{lem}

\begin{proof}
Write $\hHH^n(A) = \underline{\Hom}_{A^e}(\Omega^n_{A^e}A, A)$, $n \in \Z$, for Tate--Hochschild cohomology, computed in the stable category of $A$-bimodules; it agrees with $\HH^n(A)$ for $n \ge 1$. Two inputs are used. First, for a finite-dimensional symmetric algebra $B$ and finite-dimensional $B$-modules $X$, $Y$, Auslander--Reiten duality $\underline{\Hom}_B(X, Y)^* \cong \Ext^1_B(Y, \tau X)$ takes the form $\underline{\Hom}_B(X, Y)^* \cong \underline{\Hom}_B(Y, \Omega X)$, as $\tau = \Omega^2$ for symmetric $B$ and $\Ext^1_B(Y, \Omega^2 X) \cong \underline{\Hom}_B(Y, \Omega X)$ \cite[\S\S1--2]{Lin} (compare the stable Hochschild cohomology of periodic Calabi--Yau Frobenius algebras in \cite{ES}). With $B = A^e$, which is symmetric as $A$ is, $X = \Omega^n A$ and $Y = A$, this gives $\hHH^n(A)^* \cong \underline{\Hom}_{A^e}(A, \Omega^{n+1}A) \cong \hHH^{-n-1}(A)$. Second, the periodicity $\Omega^4_{A^e}(A) \cong A$ gives $\hHH^{m}(A) \cong \hHH^{m+4}(A)$ for every $m$; here the trivial twist is needed, as $\Omega^4_{A^e}(A) \cong A_\sigma$ for an automorphism $\sigma$ would only give $\hHH^{m+4}(A, A_\sigma)$. Hence $h^1(A) = \dim \hHH^1(A) = \dim \hHH^{-2}(A) = \dim \hHH^{2}(A) = h^2(A)$.
\end{proof}

\begin{prop}[{The smoothness criterion; extending \cite[Lemma~2.9, Prop~2.11]{Pan}}]\label{prop:criterion}
Let $A$ be the contraction algebra, of dimension $d$, of a flop with GV invariants $n_j$, and suppose $A$ deforms to $\SGV = \prod M_j(k)^{n_j}$. Then $\dim Z(A) \ge \sum n_j$, and the following are equivalent:
\begin{thmlist}
\item $[A]$ is a smooth point of $\Ass_d$;
\item $\dim Z(A) = n_1 + \dots + n_l$.
\end{thmlist}
If they hold, then the closure of $G\cdot[\SGV]$ is the unique irreducible component of $\Ass_d$ through $[A]$; every semisimple algebra to which $A$ deforms is isomorphic to $\SGV$; and $\dim Z$ is constant on a neighbourhood of $[A]$.
\end{prop}

\begin{proof}
For $\mu \in \Ass_d$ write $Z(V, \mu)$ for the centre of the algebra $(V, \mu)$. It is the kernel of the linear map $V \to \Hom_k(V, V)$, $a \mapsto \mu(a, -) - \mu(-, a)$, whose matrix entries are regular functions on $\Ass_d$, so the function $\Ass_d \to \Z_{\ge 0}$, $\mu \mapsto \dim Z(V, \mu)$, is upper semicontinuous. As $G\cdot[\SGV]$ is dense in its closure, which contains $[A]$, $\dim Z(A) \ge \dim Z(\SGV) = \sum n_j$. By Lemmas~\ref{lem:tangent} and~\ref{lem:pan},
\[
\dim T_{[A]}\Ass_d - \dim (G\cdot[\SGV]) = \dim Z(A) - \textstyle\sum_j n_j \;\ge\; 0.
\]
Since $G\cdot[\SGV]$ is open in $\Ass_d$ ($\HH^2(\SGV) = 0$), its closure is an irreducible component of dimension $\dim (G\cdot[\SGV])$ through $[A]$, so the local dimension at $[A]$ is at least $\dim (G\cdot[\SGV])$. If (2) holds, tangent and local dimension agree and $[A]$ is smooth. Conversely, if $[A]$ is smooth it lies on a unique component, which must be the closure of $G\cdot[\SGV]$, so the local dimension equals $\dim (G\cdot[\SGV]) = \dim T_{[A]}\Ass_d$, giving (2). If $S$ is semisimple and $A$ deforms to $S$, the closure of $G\cdot[S]$ is also a component through $[A]$, hence equal to that of $G\cdot[\SGV]$; two open orbits of one irreducible variety meet, so $S \cong \SGV$. Finally, near $[A]$ upper semicontinuity gives $\dim Z \le \dim Z(A) = \dim Z(\SGV)$. Conversely, every point near $[A]$ lies on the unique component, in the closure of $G\cdot[\SGV]$, so it has $\dim Z \ge \dim Z(\SGV)$, again by upper semicontinuity.
\end{proof}

\begin{rem}[Types A and D]\label{rem:AD}
In types A and D, Proposition~\ref{prop:criterion} gives smoothness from results of Jelisiejew and Smoktunowicz. Type A contraction algebras are commutative, $k[y]/(y^{n_1})$, and deform to $k^{n_1}$. In type D they deform to $\SGV$ and have $\dim Z(A) = \dim Z(\SGV)$ \cite[Rem~3.4, Thms~3.7, 3.13]{JS}.
\end{rem}
\section{The corner of a deformed preprojective algebra}\label{sec:cartan}

In this section and the next we work over the Cartan subalgebra $\h$, the left-hand column of Figure~\ref{fig:story}, where the structure of the corner is controlled by the roots. Throughout Sections~\ref{sec:cartan} and~\ref{sec:squeeze}, $\Gamma$ is a Dynkin diagram and $\bul$ any vertex, with $J$, $H_J = k[t] \otimes P$, $f_j$, $D_j$ and $l$ as in \S\ref{ss:lie}, $\Pi_\Gamma$ as in \S\ref{ss:dpa}, and $\Pi_\bul$ and $\cC_J$ as in \S\ref{ss:partial}, so that $\Pi_\bul \cong \cC_J \otimes_{H_J} k[\h] \cong \cC_J \otimes_P k[\h_J]$ (Lemma~\ref{lem:corner}). As $k[\h]$ is free of rank $|W_J|$ over $H_J$, and $k[\h_J]$ over $P$, supports and finiteness of $\cC_J$ can be checked on $\Pi_\bul$, and $\cC_J$ is free over $P$ if and only if $\Pi_\bul$ is free over $k[\h_J]$, of the same rank (\S\ref{ss:lie}). Proposition~\ref{prop:generic} restates the results of Crawley-Boevey and of Crawley-Boevey and Kimura quoted in Proposition~\ref{prop:dpa}(3), (4) at a general point of a root hyperplane. In \S\ref{ss:generic} we describe $\Pi_\bul$ along the root hyperplanes, which gives Theorem~A(2) at every vertex (Proposition~\ref{prop:support}), and in \S\ref{ss:free} we prove that $\cC_J$ is free over $P$, with a presentation over $H_J$ whose determinant is a product of the $f_j$. In Section~\ref{sec:squeeze} we prove the same for $\HH_0(\cC_J)$. The arguments do not use that $\bul$ is a flopping vertex, and we specialise to the six pairs, and to the universal contraction algebra, only in \S\ref{ss:proofA}. For $\Gamma = \mathsf{D}_4$ with $\bul$ the centre, Example~\ref{ex:running} lists the nine root hyperplanes involved and the ranks they predict.

\subsection{Generic structure and lower bounds}\label{ss:generic}

We first describe $\Pi_\bul$ near a general point of each root hyperplane. Here the roots become visible as the sizes of matrix algebras, and this gives the support of $\cC_J$ and lower bounds for the ranks of $\cC_J$ and $\HH_0(\cC_J)$.

\begin{prop}\label{prop:generic}
The algebra $\Pi_\Gamma$ and its corner $\Pi_\bul$ are described as follows.
\begin{thmlist}
\item $\Pi_\Gamma$ is a finitely generated graded $k[\h]$-module, supported on $\bigcup_{\alpha>0} \alpha^\perp$.
\item If $\lambda \in \alpha^\perp$ lies on no other root hyperplane, there is a unique simple $\Pi^\lambda$-module $S_\alpha$, of dimension vector $\alpha$; $\Ext^1(S_\alpha, S_\alpha) = 0$; and $\Pi^\lambda \cong M_{|\alpha|}(k)$ with $e_\bul\Pi^\lambda e_\bul \cong M_{\alpha_\bul}(k)$, where $\alpha_i$ is the coefficient of $\alpha$ at $i$ and $|\alpha| = \sum_i \alpha_i = \dim S_\alpha$.
\item Consequently $\Pi_\bul$ is supported on $\bigcup_{\alpha_\bul \neq 0} \alpha^\perp$, and at the generic point $\eta_\alpha$ of each such hyperplane the fibres of $\Pi_\bul$ and $\HH_0(\Pi_\bul)$ have dimensions $\alpha_\bul^2$ and $1$.
\end{thmlist}
\end{prop}

\begin{proof}
(1) The fibre at $\lambda = 0$ is the finite-dimensional preprojective algebra $\Pi(\Gamma) = \Pi^0(\Gamma)$; graded Nakayama gives finite generation. A nonzero fibre $\Pi^\lambda$ has a simple module, whose dimension vector is a positive root orthogonal to $\lambda$ (Proposition~\ref{prop:dpa}(3)), so $\lambda$ lies on some $\alpha^\perp$. (2) Only $\pm\alpha$ are orthogonal to $\lambda$, so by Proposition~\ref{prop:dpa}(3) every simple $\Pi^\lambda$-module has dimension vector $\alpha$, and there is exactly one, $S_\alpha$. The Ext formula (Proposition~\ref{prop:dpa}(4)) gives $\dim \Ext^1(S_\alpha, S_\alpha) = 2\dim\End(S_\alpha) - (\alpha, \alpha) = 2 - 2 = 0$. Hence every finite-dimensional module has all composition factors isomorphic to $S_\alpha$ and no self-extensions, so it is semisimple. Thus $\Pi^\lambda$ is semisimple with unique simple $S_\alpha$, so that $\Pi^\lambda \cong \End_k(S_\alpha) = M_{|\alpha|}(k)$, and $e_\bul$ acts on $S_\alpha$ with rank $\alpha_\bul$. (3) If $\lambda$ lies on no $\alpha^\perp$ with $\alpha_\bul \neq 0$, every simple $\Pi^\lambda$-module has dimension vector with $\alpha_\bul = 0$, so $e_\bul$ kills every simple, lies in the radical of the finite-dimensional algebra $\Pi^\lambda$, and is therefore zero; so the corner vanishes there. By (2), at a point of $\alpha^\perp$ on no other root hyperplane the fibre of $\Pi_\bul$ is $M_{\alpha_\bul}(k)$, and $\HH_0(M_{\alpha_\bul}(k))$ is one-dimensional. As $\alpha^\perp$ is integral, the restrictions of $\Pi_\bul$ and $\HH_0(\Pi_\bul)$ to $\alpha^\perp$ are free on a dense open subset (generic freeness), so the same values hold at $\eta_\alpha$.
\end{proof}

\begin{cor}\label{cor:lower}
$\cC_J$ and $\HH_0(\cC_J)$ are finite $P$-modules. Writing $\ell_\alpha$ and $q_\alpha$ for the lengths of $\Pi_\bul$ and $\HH_0(\Pi_\bul)$ at $\eta_\alpha$ ($\alpha > 0$, $\alpha_\bul \neq 0$),
\[
\rank_P \cC_J = \sum_{\alpha} \ell_\alpha \;\ge\; \sum_{\alpha} \alpha_\bul^2, \qquad \rank_P \HH_0(\cC_J) = \sum_{\alpha} q_\alpha \;\ge\; |\Phi^+_\bul|,
\]
with $\ell_\alpha = \alpha_\bul^2$ if and only if $\alpha$ annihilates $\Pi_\bul$ at $\eta_\alpha$.
\end{cor}

\begin{proof}
By Proposition~\ref{prop:generic}(3) the support of $\Pi_\bul$ is a union of hyperplanes $\alpha^\perp$ with $\alpha_\bul \neq 0$, each finite of degree $|W_J|$ over $\Spec P$ (\S\ref{ss:lie}), so $\Pi_\bul$ is finite over $P$, and hence so is $\cC_J$, a $P$-direct summand of $\Pi_\bul = \cC_J \otimes_{H_J} k[\h]$, as $H_J$ is a direct summand of $k[\h]$. At $\eta_\alpha$ the local ring is a discrete valuation ring with uniformiser $\alpha$, and by Proposition~\ref{prop:generic}(3) $\Pi_\bul$ localised there needs $\alpha_\bul^2$ generators. Hence, by Lemma~\ref{lem:dvr}(1), $\ell_\alpha \ge \alpha_\bul^2$, with equality if and only if $\alpha$ kills $\Pi_\bul$ at $\eta_\alpha$; similarly $q_\alpha \ge 1$. Over $K = \Frac P$, the fibre of each $\alpha^\perp$ is the single point $\eta_\alpha$ with $[k(\eta_\alpha) : K] = |W_J|$, and components of the support of smaller dimension do not dominate $\Spec P$. Hence $\rank_P \Pi_\bul = |W_J| \sum \ell_\alpha$, while $\rank_P \Pi_\bul = |W_J| \rank_P \cC_J$ as $k[\h]$ is free of rank $|W_J|$ over $H_J$. The same argument applies to $\HH_0$.
\end{proof}

\begin{prop}\label{prop:support}
$\cC_J$ is supported on $D = D_1 \cup \dots \cup D_l$. Let $p$ be a closed point of $D_j \setminus \pi(\Sigma)$, a dense open subset of $D_j$. Then $\pi$ is \'etale over $p$, $D$ is smooth at $p$ with local equation $f_j$, and $\cC_J \otimes_{H_J} k(p) \cong M_j(k)$. In particular $f_j$ is a reduced equation of $D_j$.
\end{prop}

\begin{proof}
By construction $\pi^{-1}(D_j) = \bigcup_{\alpha_\bul = j} \alpha^\perp$. Since $k[\h]$ is faithfully flat over $H_J$ (\S\ref{ss:lie}), the support of $\cC_J$ is the image of the support of $\Pi_\bul$, which is $\pi\big(\bigcup_{\alpha_\bul \neq 0} \alpha^\perp\big) = D$ by Proposition~\ref{prop:generic}(3). The image $\pi(\Sigma)$ of the locus $\Sigma \subset \h$ has codimension $2$ in $\Spec H_J$, so it contains no component of $D_j$. For $p$ as in the statement, every point $x \in \pi^{-1}(p)$ lies on exactly one root hyperplane $\alpha^\perp$, with $\alpha_\bul = j$, and in particular on no reflection hyperplane of $W_J$ (these are the $\beta^\perp$ with $\beta_\bul = 0$). By Steinberg's theorem, the stabiliser in $W_J$ of a point of $\h$ is generated by the reflections fixing it \cite[Thm~1.5]{Ste}. Hence $W_J$ acts freely on $\pi^{-1}(p)$ and $k[\h] \otimes_{H_J} k(p) \cong k^{|W_J|}$, so $\pi$ is \'etale over $p$. Since, near each point $x$ of $\pi^{-1}(p)$, $\pi^*f_j = \prod_{\alpha_\bul = j} \alpha$ is a unit times the linear form of the unique root hyperplane through $x$, $D$ is smooth at $p$ with local equation $f_j$. Moreover, for any $x \in \pi^{-1}(p)$ we have $k(x) = k(p) = k$, so $\cC_J \otimes_{H_J} k(p) \cong \Pi_\bul \otimes_{k[\h]} k(x) \cong M_j(k)$ by Proposition~\ref{prop:generic}(2). As $D_j \setminus \pi(\Sigma)$ is dense in $D_j$, $f_j$ is reduced.
\end{proof}

\subsection{\texorpdfstring{Freeness of the corner and generic lengths}{Freeness of the corner and generic lengths}}\label{ss:free}

We now use the results of Etingof and Rains to show that $\cC_J$ is free over $P$ and to compute its rank. Comparing this rank with the lower bounds of \S\ref{ss:generic} then gives the lengths along each root hyperplane $\alpha^\perp$ with $\alpha_\bul \neq 0$, and with them the determinants of the presentations in Lemma~\ref{lem:presentation}.

\begin{prop}[{Freeness of the corner; extending \cite[Thm~3.3(i)]{ER}}]\label{prop:ER}
$\cC_J$ is a graded free $P$-module of rank $\sum_{\alpha>0} \alpha_\bul^2$, and $\Pi_\bul$ is a graded free $k[\h_J]$-module of the same rank.
\end{prop}

\begin{proof}
For the nodal vertex (the branch vertex in types D and E, the middle vertex in type $\mathsf{A}_{2m-1}$), \cite[Thm~3.3(i)]{ER} proves that the corner with a central variable adjoined is free over a polynomial ring in $|\Gamma_0|$ variables, which amounts to the Cohen--Macaulay statement for that corner; the following argument works at every vertex. By Proposition~\ref{prop:dpa}(5), $\Pi_\Gamma$ is a Cohen--Macaulay $k[\h]$-module; its corner $\Pi_\bul$ is a direct summand, so also Cohen--Macaulay. Recall that $k[\h_J] = k[t_i : i \in J]$. The corner is finite over $P \subset k[\h_J]$ (proof of Corollary~\ref{cor:lower}), hence over $k[\h_J]$, and its support is a union of root hyperplanes (Proposition~\ref{prop:generic}), so its dimension is $\dim k[\h_J]$. Hence $\Pi_\bul$ is free over $k[\h_J]$ by Lemma~\ref{lem:freeness}(2). As $\Pi_\bul = \cC_J \otimes_P k[\h_J]$ (Lemma~\ref{lem:corner}), $\cC_J$ is graded free over $P$ (\S\ref{ss:lie}).

It remains to compute the rank, which the fibres in Proposition~\ref{prop:generic} bound only from below (Corollary~\ref{cor:lower}). By Proposition~\ref{prop:dpa}(5) \cite[\S2.1, Thm~3.1(ii)]{ER}, the $(\bul, \bul)$ entry of the matrix Hilbert series of $\Pi^\lambda(\Gamma)$ is
\[
\operatorname{HS}_{\Pi_\bul}(u)
=\frac{1-u^{2h}}{(1-u^2)^{n}}
\bigl((I-uC_\Gamma+u^2I)^{-1}\bigr)_{\bul\bul}.
\]
As $\Pi_\bul$ is graded free over $k[\h_J]$, a polynomial ring in $n - 1$ variables of degree $2$, its rank is the value at $u = 1$ of $(1-u^2)^{n-1}\operatorname{HS}_{\Pi_\bul}(u)$:
\[
\rank_{k[\h_J]}\Pi_\bul
=\lim_{u\to1}\frac{1-u^{2h}}{1-u^2}
\bigl((I-uC_\Gamma+u^2I)^{-1}\bigr)_{\bul\bul}
=h(A_\Gamma^{-1})_{\bul\bul}.
\]
By \S\ref{ss:lie}, $\rank_P\cC_J=\rank_{k[\h_J]}\Pi_\bul$ and $h(A_\Gamma^{-1})_{\bul\bul} = \sum_{\alpha>0}\alpha_\bul^2$.
\end{proof}

For the six pairs of \S\ref{ss:universal}, Remark~\ref{rem:noER} explains how Proposition~\ref{prop:ER} can be verified without \cite{ER}.

\begin{cor}[Generic lengths]\label{cor:trace-generic}
One has
\[
\ell_\alpha=\alpha_\bul^2,\qquad q_\alpha=1
\quad(\alpha_\bul\ne0),\qquad \rank_P \HH_0(\cC_J)=|\Phi^+_\bul|.
\]
\end{cor}

\begin{proof}
Corollary~\ref{cor:lower} gives $\ell_\alpha\ge\alpha_\bul^2$, and the sum of these inequalities is an equality by Proposition~\ref{prop:ER}. Hence $\alpha$ annihilates $(\Pi_\bul)_{\eta_\alpha}$ by the last assertion of Corollary~\ref{cor:lower}. It also annihilates its quotient $\HH_0(\Pi_\bul)_{\eta_\alpha}$, whose fibre is one-dimensional by Proposition~\ref{prop:generic}(3); so $q_\alpha=1$ (Lemma~\ref{lem:dvr}(1)). The rank formula follows from Corollary~\ref{cor:lower}.
\end{proof}

Free $P$-modules with known generic lengths, such as $\cC_J$ and (by Section~\ref{sec:squeeze}) $\HH_0(\cC_J)$, have presentations $(\ast)$ over $H_J$ whose determinants are computed by roots:

\begin{lem}[Presentations]\label{lem:presentation}
Let $M$ be a finite graded $H_J$-module that is free over $P$, such that $M \otimes_{H_J} k[\h]$ is supported on the union of the $\alpha^\perp$ with $\alpha_\bul \neq 0$ and has length $m_\alpha$ at $\eta_\alpha$, with $m_\alpha = m_j$ depending only on $j = \alpha_\bul$. Then there is an exact sequence of graded $H_J$-modules
\[
0 \to M \otimes_P H_J \xrightarrow{\;1 \otimes t - t \otimes 1\;} M \otimes_P H_J \to M \to 0, \tag{$\ast$}
\]
whose first map, in a $P$-basis of $M$, is $t\cdot I - [t]$, with $[t]$ the matrix of multiplication by $t$. In particular, $M$ has projective dimension at most $1$ over $H_J$ and is zero or Cohen--Macaulay of dimension $\dim H_J - 1$. The determinant of $(\ast)$, the characteristic polynomial of $t$ on $M$, is $c \prod_{j} f_j^{\,m_j}$ with $c \in k^\times$, and after base change to $k[\h]$ it is $c \prod_{\alpha \in \Phi^+_\bul} \alpha^{m_\alpha}$.
\end{lem}

\begin{proof}
The first map of $(\ast)$ is injective, as its determinant is monic in $t$ and hence nonzero. As $M \otimes_P H_J = M[t]$, its cokernel is $M[t]/(t - [t])M[t]$, which the action map $M \otimes_P H_J \to M$ identifies with $M$; so $(\ast)$ is exact. Hence $\pd M \le 1$, and $M$ is Cohen--Macaulay of dimension $\dim H_J - 1$ (Lemma~\ref{lem:freeness}(3)).

Let $d \in H_J$ be the determinant of $(\ast)$, a nonzero homogeneous element. The base change of $(\ast)$ to $k[\h]$ is exact, as $k[\h]$ is flat over $H_J$ (\S\ref{ss:lie}), with determinant $\pi^*d$. Localised at $(\alpha)$, $k[\h]$ is a discrete valuation ring with uniformiser $\alpha$, so by Lemma~\ref{lem:dvr}(2) $\pi^*d$ vanishes along $\alpha^\perp$ to order $m_\alpha$. Any other prime factor of $\pi^*d$ would lie in the support of the cokernel, which is the union of the $\alpha^\perp$ with $\alpha_\bul \neq 0$. As $k[\h]$ is a polynomial ring, $\pi^*d = c\prod_{\alpha}\alpha^{m_\alpha}$. Grouping the roots by $\alpha_\bul$ gives $\pi^*d = c\,\pi^*\big(\prod_j f_j^{m_j}\big)$, and $\pi^*$ is injective.
\end{proof}

\begin{ex}[Length 2]\label{ex:running}
Let $\Gamma = \mathsf{D}_4$ with $\bul$ the centre, so that $\cC_J$ is the universal contraction algebra $\cC_2$ of length $2$ (\S\ref{ss:universal}). The positive roots with $\alpha_\bul \neq 0$ are the highest root $\theta$, with $\alpha_\bul = 2$, and eight roots with $\alpha_\bul = 1$: for each $U \subseteq J$, the root $\alpha_U$ with coefficient $1$ at $\bul$ and at the vertices of $U$ and $0$ elsewhere. Thus $\Pi_\bul$ is supported on nine hyperplanes: by Proposition~\ref{prop:generic}, it is $M_1(k) = k$ along the eight $\alpha_U^\perp$ and $M_2(k)$ along $\theta^\perp = \{t = 0\}$. This predicts rank $8\cdot 1 + 1\cdot 4 = 12$ for $\cC_J$ and $8 + 1 = 9$ for $\HH_0(\cC_J)$, as Proposition~\ref{prop:ER} and Corollary~\ref{cor:trace-generic} confirm. Details are in \S\ref{ss:rootformula2}, and Example~\ref{ex:lauferarc} writes $\cC_2$ down by three relations over $H_2$.
\end{ex}

\section{\texorpdfstring{Hochschild homology of the corner}{Hochschild homology of the corner}}\label{sec:squeeze}

In this section we prove that $\HH_0(\cC_J)$ is free over $P$, which is the main new step in the proof of Theorem~A. Corollary~\ref{cor:trace-generic} computes the generic rank of $\HH_0(\cC_J)$ by roots, and freeness is a Cohen--Macaulay property. Here $\HH_0(\Pi_\bul) = \HH_0(\cC_J) \otimes_{H_J} k[\h]$, as $\HH_0$ commutes with base change (\S\ref{ss:moduli}). The support of $\HH_0(\Pi_\bul)$ is the union of the root hyperplanes $\alpha^\perp$ with $\alpha_\bul \neq 0$ (Proposition~\ref{prop:generic}). This union is finite over $k[\h_J]$, so $\HH_0(\Pi_\bul)$ is a finite $k[\h_J]$-module of dimension $n - 1$, where $n = \dim \h$. Hence $\HH_0(\Pi_\bul)$ is a Cohen--Macaulay $k[\h]$-module of dimension $n - 1$ if and only if it is free over $k[\h_J]$ (Lemma~\ref{lem:freeness}(2)), that is (\S\ref{ss:lie}), if and only if $\HH_0(\cC_J)$ is free over $P$. The Cohen--Macaulay theorem of Etingof and Rains applies to $\Pi_\bul$, but Cohen--Macaulayness need not pass to the quotient $\Pi_\bul/[\Pi_\bul, \Pi_\bul]$.

We therefore need a different argument, outlined in the introduction (``Freeness of $\HH_0$''): a classical fact about the Weyl group says that the vector fields on $\h$ tangent to every root hyperplane form a free module (Lemma~\ref{lem:logbasis}(1)), and we show that $\HH_0$ of the whole algebra $\Pi_\Gamma$ is the quotient of all vector fields by these (Proposition~\ref{prop:wholeHH}).

In order to prove freeness we compare $\HH_0(\Pi_\Gamma)$ with the vector fields on $\h$, using the simultaneous resolution $\cY_\Gamma \to \h$ of \S\ref{ss:simres} to lift them. In \S\ref{ss:wholeHH} we show that $\HH_0(\Pi_\Gamma)$ is Cohen--Macaulay, in \S\ref{ss:transfer} that $\HH_0(\cC_J)$ is a direct summand of it, and in \S\ref{ss:proofA3} we deduce that $\HH_0(\cC_J)$ is free. The only result of \S\ref{ss:wholeHH} needed for Theorem~\ref{thm:squeeze} is Proposition~\ref{prop:wholeHH}(1), that $\HH_0(\Pi_\Gamma)$ is a finite graded Cohen--Macaulay $k[\h]$-module of dimension $n - 1$. Part~(2) describes $\HH_0(\Pi_\Gamma)$ explicitly. In \S\ref{ss:proofA} we collect the results at every vertex (Theorem~\ref{thm:Avertex}) and deduce Theorem~A, and in \S\ref{ss:centre} we show that the centre of the corner is free as well.

\subsection{\texorpdfstring{$\HH_0$ of the deformed preprojective algebra}{HH0 of the deformed preprojective algebra}}\label{ss:wholeHH}

In this subsection we find the kernel of the $k[\h]$-linear map
\[
\kappa\colon \Theta_\h \to \HH_0(\Pi_\Gamma), \qquad \xi \mapsto \textstyle\sum_{i \in \Gamma_0} \xi(t_i)\,[e_i],
\]
from the module $\Theta_\h$ of vector fields on $\h$. The map $\kappa$ applies $\xi$ to the right-hand side $\sum_i \lambda_i e_i$ of the defining relation, as $\lambda_i = t_i$ (\S\ref{ss:lie}). Our aim is to show that $\ker\kappa$ is the module $\Der(-\log\mathcal H)$ of vector fields tangent to every root hyperplane (Lemma~\ref{lem:KS}), so that $\HH_0(\Pi_\Gamma)$ is the Cohen--Macaulay module $\Theta_\h/\Der(-\log\mathcal H)$ (Proposition~\ref{prop:wholeHH}).

Two properties of $\kappa$ are elementary: $\kappa$ is onto, as $\HH_0(\Pi_\Gamma)$ is generated by the $[e_i]$ (Lemma~\ref{lem:HHgen}), and $\kappa$ kills every vector field that lifts to a derivation of $\cA_\Gamma$, as one sees by applying the derivation to the defining relation (Lemma~\ref{lem:KS}). What remains is to lift the vector fields tangent to every root hyperplane (Lemma~\ref{lem:derlift}), and we do this in two steps. First we lift them from $\h$ to vector fields on the simultaneous resolution $\chi\colon \cY_\Gamma \to \h$. Then we lift these to derivations of $\cA_\Gamma = \End(\mathcal T)$, using the defining property $\Ext^1(\mathcal T, \mathcal T) = 0$ of the tilting bundle $\mathcal T$ of Proposition~\ref{prop:simres}(1). Locally $\sEnd(\mathcal T)$ is a matrix algebra and the lift is immediate, and this vanishing lets the local lifts glue. A length count at the root hyperplanes then gives equality (Lemma~\ref{lem:KS}). In geometric terms $\kappa$ plays the role of the Kodaira--Spencer map of $\chi$ (see the paragraph after Proposition~\ref{prop:wholeHH}). The family $\cY_J$ and its algebra $A_J$ enter only in \S\ref{ss:transfer}. We first collect what we need about these vector fields in Lemma~\ref{lem:logbasis}, using the equations of \S\ref{ss:kleinian}: part~(1) gives the Cohen--Macaulay property, and part~(2) is the key to the first lifting step.

\paragraph{Logarithmic vector fields.} Write $\mathcal H = \bigcup_{\alpha > 0} \alpha^\perp \subset \h$ for the union of the root hyperplanes, $\Theta_\h = \Der(k[\h]) = k[\h] \otimes \h$, and
\[
\Der(-\log \mathcal H) = \{\xi \in \Theta_\h : \xi(\alpha) \in \alpha\, k[\h] \text{ for all } \alpha > 0\}
\]
for the derivations of $k[\h]$ that preserve each ideal $(\alpha)$, that is, the vector fields tangent to every root hyperplane. As the $\alpha > 0$ are pairwise coprime, this is equivalent to preserving the ideal $(\prod_{\alpha > 0}\alpha)$ of $\mathcal H$, the condition in Saito's criterion \cite[(1.8)(ii)]{Sai}: $n$ derivations preserving the ideal of a reduced equation, whose determinant is that equation times a nonzero scalar (a unit, in Saito's analytic setting), form a basis of the module of all derivations preserving it. Saito states it for germs of analytic functions. We use it for weighted homogeneous polynomials, where it follows from the analytic statement at the origin, because the analytic local ring at $0$ is faithfully flat over the algebraic local ring at $0$, and a graded module whose localisation at the origin vanishes is zero. For hyperplane arrangements see also \cite[Thm~4.19]{OT}.

\begin{lem}\label{lem:logbasis}
Let $p_1, \dots, p_n$ be basic invariants of $W$, that is, homogeneous generators of the polynomial ring $k[\h]^W$, let $(\ ,\ )$ be the invariant form of \S\ref{ss:lie} on $\h^*$, extended $k[\h]$-bilinearly to $1$-forms, and let $\operatorname{grad} p_r$ be the derivation $f \mapsto (dp_r, df)$ of $k[\h]$.
\begin{thmlist}
\item The vector fields $\operatorname{grad} p_1, \dots, \operatorname{grad} p_n$ form a basis of $\Der(-\log\mathcal H)$, and $\Theta_\h/\Der(-\log\mathcal H)$ is a Cohen--Macaulay $k[\h]$-module of dimension $n - 1$, whose associated primes are the $(\alpha)$, $\alpha > 0$, at each of which it has length one.
\item Every vector field $\xi$ on $\h/W$ tangent to $\mathcal D$, that is, with $\xi(d_W) \in (d_W)$, lifts to a vector field $\Xi = \sum_j a_j\partial_{u_j} - \sum_i v_i\partial_{x_i}$ on $\A^3 \times \h/W$ with polynomial coefficients, where $\xi = \sum_j a_j\partial_{u_j}$ and $(x_1, x_2, x_3) = (x, y, z)$, such that $\Xi(F) = \psi F$ for a polynomial $\psi$; in particular $\Xi$ is tangent to $\cX_\Gamma$.
\end{thmlist}
\end{lem}

\begin{proof}
(1) Let $\alpha > 0$. The reflection $s_\alpha$ preserves the form and $p_r$, and $s_\alpha(\alpha) = -\alpha$, so $\operatorname{grad} p_r(\alpha) = (dp_r, \alpha)$ is anti-invariant under $s_\alpha$. Hence it vanishes on the fixed hyperplane $\alpha^\perp$ of $s_\alpha$, and lies in $(\alpha)$: so $\operatorname{grad} p_r \in \Der(-\log\mathcal H)$. Their determinant $\det(\operatorname{grad} p_r(t_j))$ is the Jacobian $\det(\partial p_r/\partial t_i)$ times the determinant of the form, so a nonzero scalar multiple of $\prod_{\alpha>0}\alpha$ \cite[Ch.~V, \S5, no.~4, Prop.~5(ii)]{Bou}, and they form a basis by Saito's criterion. Hence $\Theta_\h/\Der(-\log\mathcal H)$ is the cokernel of an $n \times n$ matrix with determinant $c\prod_{\alpha>0}\alpha$, $c \in k^\times$, so it is Cohen--Macaulay of dimension $n - 1$ (Lemma~\ref{lem:freeness}(3)). Its associated primes are therefore minimal, so they are the $(\alpha)$, and its length at $(\alpha)$ is the order of the determinant there, which is one (Lemma~\ref{lem:dvr}(2)).

(2) We lift a basis of the vector fields on $\h/W$ tangent to $\mathcal D$.

\emph{The Jacobian algebra.} Let $\Jac(F) = \cO_{\A^3 \times \h/W}/(\partial_x F, \partial_y F, \partial_z F)$. It is a complete intersection, finite over $\h/W$ by the grading, hence free over it (Lemma~\ref{lem:freeness}(2)). Its fibre over the origin is the Milnor algebra $k[x, y, z]/(\partial_x F_0, \partial_y F_0, \partial_z F_0)$, in which $F_0$ vanishes by the Euler relation, so that it is the algebra $k[x, y, z]/(F_0, \partial_x F_0, \partial_y F_0, \partial_z F_0)$ of Slodowy. The $\mu_j$ restrict to a basis of it, as the semi-universal deformation is built from representatives of such a basis \cite[2.4, 8.6]{Slo}. Hence, by graded Nakayama, $\mu_0, \dots, \mu_{n-1}$ generate $\Jac(F)$, and so form a basis of this free module of rank $n$.

\emph{Lifts.} For each index $r$, expanding $\mu_r F$ in this basis gives $a_{jr} \in k[\h]^W$ and $v_{ir} \in \cO(\A^3 \times \h/W)$ with $\mu_r F = \sum_j a_{jr}\mu_j - \sum_i v_{ir}\,\partial_{x_i}F$; thus $(a_{jr})$ is the matrix of multiplication by $F$ on $\Jac(F)$. Put $\Xi_r = \sum_j a_{jr}\partial_{u_j} - \sum_i v_{ir}\partial_{x_i}$ and $\xi_r = \sum_j a_{jr}\partial_{u_j}$. Since $\partial F/\partial u_j = \mu_j$, we get $\Xi_r(F) = \mu_r F$, so $\Xi_r$ is tangent to $\cX_\Gamma$, and it lifts $\xi_r$.

\emph{The $\xi_r$ are tangent to $\mathcal D$.} Let $I = (F, \partial_xF, \partial_yF, \partial_zF)$, an ideal defining the singular points of the fibres of $\cX_\Gamma$. Then $\Xi_r(I) \subseteq I$, as $\Xi_r(F) = \mu_rF$ and $\Xi_r(\partial_{x_i}F) = \partial_{x_i}(\mu_rF) + \sum_{i'} (\partial_{x_i}v_{i'r})\,\partial_{x_{i'}}F$, the $a_{jr}$ being constant in $x, y, z$. Hence $\xi_r$ preserves $I_0 = I \cap k[\h]^W$. The zero set of $I_0$ is $\mathcal D$, the image of the zero set of $I$, which is finite over $\h/W$ as it lies in $\Spec\Jac(F)$. Hence $I_0$ has radical $(d_W)$, and in the discrete valuation ring $R = k[\h]^W_{(d_W)}$ we have $I_0R = (d_W^m)$ for some $m \ge 1$. Then $\xi_r$ extends to $R$ by the quotient rule and preserves $I_0R$. Differentiating, $m\,d_W^{m-1}\xi_r(d_W) \in (d_W^m)$ there, so, as $m$ is invertible in $k$ and $(d_W)$ is prime, $\xi_r(d_W) \in (d_W)$.

\emph{The $\xi_r$ form a basis.} The determinant $\det(a_{jr})$ of multiplication by $F$ on $\Jac(F)$ vanishes exactly on $\mathcal D$, over which a fibre of $\cX_\Gamma$ is singular. Near a general point of $\mathcal D$ the fibre has a single singular point, an $\mathsf{A}_1$ point (\S\ref{ss:kleinian}), and $\Jac(F)$ is there the product of a factor isomorphic to the base, on which $F$ acts by the critical value, and a factor on which $F$ is a unit. As the critical value has derivative $1$ with respect to $u_0$, $\det(a_{jr})$ vanishes to order one there. As $d_W$ is prime in $k[\h]^W$ and the units of $k[\h]^W$ are the scalars, $\det(a_{jr})$ is $d_W$ up to a nonzero scalar, so the $\xi_r$ form a basis of the vector fields tangent to $\mathcal D$ by Saito's criterion.

Hence $\xi = \sum_r b_r\xi_r$ with $b_r$ polynomial, and $\Xi = \sum_r b_r\Xi_r$ lifts $\xi$, with $\Xi(F) = \psi F$ for $\psi = \sum_r b_r\mu_r$.
\end{proof}

\begin{lem}[{Generators of $\HH_0$; extending \cite[Thm~1.4(a)]{MOV}}]\label{lem:HHgen}
$\HH_0(\Pi_\Gamma)$ is a finite graded $k[\h]$-module, generated by the classes $[e_i]$, $i \in \Gamma_0$. In particular $\kappa$ is onto.
\end{lem}

\begin{proof}
$\Pi_\Gamma$ is a finite graded $k[\h]$-module, with $k[\h]$ in positive degrees (Proposition~\ref{prop:generic}(1)), and hence so is its quotient $\HH_0(\Pi_\Gamma)$. As $\HH_0$ commutes with base change, by graded Nakayama it suffices to show that the fibre $\HH_0(\Pi^0(\Gamma))$ at $0$ is spanned by the $[e_i]$. This is known \cite[Thm~1.4(a)]{MOV} (see also \cite[Thm~3.0.3]{EE}), with a case-by-case proof; we give a uniform one.

Let $G_\Gamma \subset SL_2(k)$ be the finite group of the Kleinian singularity and $\mathcal S = k[x, y] \rtimes G_\Gamma$. By \cite[Thm~3.4]{CBH}, $\Pi^0(\tGamma) \cong e'\mathcal Se'$ with $e_i \mapsto e'_i$, where $e' = \sum_{i \in \tGamma_0} e'_i$ and $e'_i \in kG_\Gamma$ is a primitive idempotent for the irreducible representation $\rho_i$ of $G_\Gamma$ at the vertex $i$ of the McKay graph. As $\rho_0$ is trivial, $e'_0 = \frac{1}{|G_\Gamma|}\sum_{\gamma \in G_\Gamma}\gamma$. As $kG_\Gamma e'kG_\Gamma = kG_\Gamma$, $e'$ is full, so the inclusion induces $\HH_0(e'\mathcal Se') \cong \HH_0(\mathcal S)$ (\S\ref{ss:moduli}). As $\Pi^0(\Gamma) = \Pi^0(\tGamma)/(e_0)$ (Lemma~\ref{lem:corner}), by \eqref{eq:idem} $\HH_0(\Pi^0(\Gamma))$ is the quotient of $\HH_0(\mathcal S)$ by the image of $\HH_0(e'_0\mathcal Se'_0)$. We show that
\[
\HH_0(\mathcal S) = k[x, y]^{G_\Gamma} \oplus \textstyle\bigoplus_{C \neq \{1\}} k\,[\gamma_C],
\]
the sum over the nontrivial conjugacy classes $C$ of $G_\Gamma$, with $\gamma_C \in C$, where $a \in k[x, y]^{G_\Gamma}$ stands for its class $[a]$.

Let $a \in k[x, y]$ and $\gamma \in G_\Gamma$. If $\gamma \neq 1$, then $\gamma - 1$ is invertible on $kx + ky$, as $\gamma \in SL_2(k)$ has finite order, so $x$ and $y$ lie in the ideal generated by the $b - \gamma(b)$, $b \in k[x, y]$; as $[b, c\gamma] = (b - \gamma(b))c\gamma$ for $b, c \in k[x, y]$, this gives $[a\gamma] = a(0)[\gamma]$, and $[\gamma] = [\gamma_C]$ for $\gamma \in C$ by conjugation. If $\gamma = 1$, then $[\gamma', a\gamma'^{-1}] = \gamma'(a) - a$ for $\gamma' \in G_\Gamma$ gives $[a] = [\bar a]$, with $\bar a \in k[x, y]^{G_\Gamma}$ the average of the $\gamma'(a)$. Hence the right-hand side spans. The sum is direct, as the linear maps sending $\sum_\gamma a_\gamma\gamma$ to $\bar a_1$ and to $\sum_{\gamma \in C} a_\gamma(0)$ vanish on commutators $[b\delta, c\varepsilon] = b\delta(c)\,\delta\varepsilon - c\varepsilon(b)\,\varepsilon\delta$.

Now $e'_0\mathcal Se'_0 = k[x, y]^{G_\Gamma}e'_0$, and $[\bar ae'_0] = |G_\Gamma|^{-1}\big([\bar a] + \bar a(0)\sum_{\gamma \neq 1}[\gamma]\big)$. Thus the quotient $\HH_0(\Pi^0(\Gamma))$ is spanned by the classes of $kG_\Gamma$, and these are spanned by the $[e'_i]$, as $kG_\Gamma = \bigoplus_i M_{\dim\rho_i}(k)$.
\end{proof}

\begin{lem}\label{lem:derlift}
For every $\xi \in \Der(-\log\mathcal H)$ there is a $k$-linear derivation $\partial$ of $\cA_\Gamma$ that restricts to $\xi$ on $k[\h]$.
\end{lem}

\begin{proof}
As $k[\h]$ is central in $\cA_\Gamma$, the vector fields on $\h$ that are restrictions of derivations of $\cA_\Gamma$ form a $k[\h]$-submodule. Hence, by Lemma~\ref{lem:logbasis}(1), it suffices to treat $g = \operatorname{grad} p_r$.

\emph{From $\h$ to the simultaneous resolution.} As $p_r$ and the form are $W$-invariant, $g$ is $W$-equivariant, so it preserves $k[\h]^W$ and induces a vector field $\bar g$ on $\h/W$. It is tangent to $\mathcal D$: $g$ preserves each ideal $(\alpha)$, hence also the ideal $(d_W) = (\prod\alpha^2)$ of $k[\h]$, and $g(d_W)/d_W$ is $W$-invariant, so $g(d_W) \in d_W\,k[\h]^W$. Let $\Xi$ be the lift of $\bar g$ given by Lemma~\ref{lem:logbasis}(2); as $\Xi(F) = \psi F$, it is a derivation of $\cO(\cX_\Gamma)$ restricting to $\bar g$ on $k[\h]^W$. Write $\cX' = \cX_\Gamma \times_{\h/W} \h$. Since $\Xi$ and $g$ agree on $k[\h]^W$, $\Xi'(a \otimes b) = \Xi(a) \otimes b + a \otimes g(b)$ is a well-defined derivation of $\cO(\cX') = \cO(\cX_\Gamma) \otimes_{k[\h]^W} k[\h]$, restricting to $g$ on $k[\h]$.

To carry $\Xi'$ to $\cY_\Gamma$, let $\operatorname{Exc} \subset \cY_\Gamma$ be the union of the positive-dimensional fibres of $\phi\colon \cY_\Gamma \to \cX'$. It lies over $\mathcal H$ and meets each fibre of $\chi$ in exceptional curves, so $\dim \operatorname{Exc} \le n$, and $\operatorname{Exc}$ has codimension at least $2$ in $\cY_\Gamma$, which has dimension $n + 2$. Off $\operatorname{Exc}$, $\phi$ is an isomorphism onto $\cX' \setminus \phi(\operatorname{Exc})$ by Zariski's main theorem. Indeed, as $\operatorname{Exc}$ is a union of fibres, $\phi$ is proper and quasi-finite over $\cX' \setminus \phi(\operatorname{Exc})$, hence finite. It is birational, as $\cY_\Gamma$ is smooth and connected (the $k^\times$-action contracts it onto $\chi^{-1}(0)$) and $\phi$ is an isomorphism over $\h \setminus \mathcal H$, where it is a fibrewise isomorphism of flat families. Finally, $\cX'$ is normal: it is a hypersurface whose singular locus, finitely many points in each fibre over $\mathcal H$, has dimension at most $n - 1$. Through $\cY_\Gamma \setminus \operatorname{Exc} \cong \cX' \setminus \phi(\operatorname{Exc})$, $\Xi'$ gives a vector field on $\cY_\Gamma \setminus \operatorname{Exc}$, which extends to a vector field $V$ on $\cY_\Gamma$ as $T_{\cY_\Gamma}$ is locally free and $\operatorname{Exc}$ has codimension at least $2$. Then $V(\chi^*b) = \chi^*g(b)$ for $b \in k[\h]$, as this holds off $\operatorname{Exc}$.

\emph{From $\cY_\Gamma$ to $\cA_\Gamma = \End(\mathcal T)$.} Recall that $\cA_\Gamma = \End(\mathcal T)$, with $k[\h]$ acting through $\chi^*$ (Proposition~\ref{prop:simres}(1)), and let $\mathcal E = \sEnd(\mathcal T)$, where $\mathcal T$ has rank $N$. Cover $\cY_\Gamma$ by finitely many affine open subsets $U_j$ on which $\mathcal T$ is trivial, and let $V_j$ be the derivation of $\mathcal E|_{U_j} \cong M_N(\cO_{U_j})$ that applies $V$ to each matrix entry, so that $V_j(bs) = V(b)s + bV_j(s)$ for a function $b$. On $U_j \cap U_{j'}$ the difference $V_j - V_{j'}$ is an $\cO$-linear derivation of a matrix algebra over a commutative ring, hence inner: $V_j - V_{j'} = [c_{jj'}, -]$, with $c_{jj'}$ a section of $\mathcal E/\cO$. These form a \v{C}ech cocycle, and $\mathcal E/\cO$ is a direct summand of $\mathcal E$, split by $\frac1N\operatorname{tr}$, so its class lies in $H^1(\cY_\Gamma, \mathcal E/\cO) \subseteq H^1(\cY_\Gamma, \mathcal E) = \Ext^1(\mathcal T, \mathcal T) = 0$. As $\cY_\Gamma$ is separated, \v{C}ech cohomology of the affine cover computes this group, so $c_{jj'} = c_{j'} - c_j$ modulo $\cO$ for sections $c_j$ of $\mathcal E$ over $U_j$, and the derivations $V_j + [c_j, -]$ glue to a derivation of $\mathcal E$. On global sections it gives a derivation $\partial$ of $\cA_\Gamma$, with $\partial(b) = V(\chi^*b) = g(b)$ for $b \in k[\h]$.
\end{proof}

\begin{lem}\label{lem:KS}
$\ker\kappa = \Der(-\log\mathcal H)$, so $\kappa$ induces an isomorphism $\Theta_\h/\Der(-\log\mathcal H) \cong \HH_0(\Pi_\Gamma)$ with $\partial/\partial t_i \mapsto [e_i]$.
\end{lem}

\begin{proof}
\emph{$\Der(-\log\mathcal H) \subseteq \ker\kappa$.} Let $\xi \in \Der(-\log\mathcal H)$, and let $\partial$ be as in Lemma~\ref{lem:derlift}. With $c = \sum_{i \in \tGamma_0} \partial(e_i)e_i$, differentiating $e_je_i = \delta_{ij}e_j$ gives $[c, e_j] = \partial(e_j)$, and $[c, -]$ vanishes on the central $k[\h]$. Hence, replacing $\partial$ by $\partial - [c, -]$, we may assume that $\partial(e_i) = 0$ for all $i \in \tGamma_0$. Applying $\partial$ to the defining relation $\sum_a [a, a^*] = \sum_{i \in \tGamma_0} \lambda_i e_i$ of $\cA_\Gamma$ gives
\[
\textstyle\sum_{i \in \tGamma_0} \xi(\lambda_i)\, e_i = \sum_a \big([\partial a, a^*] + [a, \partial a^*]\big) \in [\cA_\Gamma, \cA_\Gamma].
\]
In $\Pi_\Gamma = \cA_\Gamma/\cA_\Gamma e_0\cA_\Gamma$ (Lemma~\ref{lem:corner}) this becomes $\sum_{i \in \Gamma_0} \xi(t_i)e_i \in [\Pi_\Gamma, \Pi_\Gamma]$, that is, $\kappa(\xi) = 0$.

\emph{Equality.} Hence $\kappa$ induces a surjection $\bar\kappa\colon \bar\Theta = \Theta_\h/\Der(-\log\mathcal H) \to \HH_0(\Pi_\Gamma)$ (Lemma~\ref{lem:HHgen}). By Lemma~\ref{lem:logbasis}(1), $\bar\Theta$ is Cohen--Macaulay, its associated primes are the $(\alpha)$, $\alpha > 0$, and $\bar\Theta$ has length $1$ at each of them. On the other hand $\HH_0(\Pi_\Gamma)$ is nonzero at the generic point of each $\alpha^\perp$, as $\HH_0$ commutes with base change and $\Pi^\lambda \cong M_{|\alpha|}(k)$ there (Proposition~\ref{prop:generic}(2)). Hence $\bar\kappa$ is an isomorphism at every associated prime of $\bar\Theta$, and its kernel, a submodule of $\bar\Theta$, has no associated primes, so it is zero.
\end{proof}

\begin{prop}\label{prop:wholeHH}
The zeroth Hochschild homology of $\Pi_\Gamma$ is described as follows.
\begin{thmlist}
\item $\HH_0(\Pi_\Gamma) \cong \Theta_\h/\Der(-\log\mathcal H)$ as $k[\h]$-modules (up to a shift of grading), and this is a finite graded Cohen--Macaulay $k[\h]$-module of dimension $n - 1$.
\item Let $\alpha > 0$, write $k[\alpha^\perp] = k[\h]/(\alpha)$, and for $\lambda \in \alpha^\perp \setminus \Sigma$ let $S_\alpha(\lambda)$ be the simple $\Pi^\lambda$-module of Proposition~\ref{prop:generic}(2). For $a \in \Pi_\Gamma$ the function $\lambda \mapsto \operatorname{tr}(a \mid S_\alpha(\lambda))$ on $\alpha^\perp \setminus \Sigma$ extends to an element $\operatorname{tr}_\alpha(a)$ of $k[\alpha^\perp]$, and the traces define an injective $k[\h]$-linear map
\[
\operatorname{tr} = (\operatorname{tr}_\alpha)_{\alpha > 0}\colon \HH_0(\Pi_\Gamma) \hookrightarrow \bigoplus_{\alpha > 0} k[\alpha^\perp], \qquad \operatorname{tr}[e_i] = (\alpha_i)_{\alpha > 0},
\]
which is an isomorphism over $\h \setminus \Sigma$. The isomorphism of (1) sends $[e_i]$ to the constant vector field $\partial/\partial t_i$, and $\operatorname{tr}$ is evaluation on the roots, $\xi \mapsto (\xi(\alpha)|_{\alpha^\perp})_{\alpha > 0}$.
\end{thmlist}
\end{prop}

\begin{proof}
(1) Lemmas~\ref{lem:KS} and~\ref{lem:logbasis}(1). (2) Let $R_\alpha = k[\alpha^\perp]$, a polynomial ring with fraction field $k(\eta_\alpha)$ (\S\ref{ss:lie}), and $B = \Pi_\Gamma \otimes_{k[\h]} R_\alpha$, a finite $R_\alpha$-module. Over $U = \alpha^\perp \setminus \Sigma$ the fibres of $B$ are $M_{|\alpha|}(k)$ (Proposition~\ref{prop:generic}(2)), of constant dimension, so $B|_U$ is locally free and $\lambda \mapsto \operatorname{tr}(a \mid S_\alpha(\lambda)) = |\alpha|^{-1}\operatorname{tr}(a\cdot \mid \Pi^\lambda)$ is regular on $U$. It extends to $\alpha^\perp$: $a$ is integral over $R_\alpha$, so the eigenvalues of $a$ on the simple module of $\Pi_\Gamma \otimes_{k[\h]} \overline{k(\eta_\alpha)} \cong M_{|\alpha|}(\overline{k(\eta_\alpha)})$ (Proposition~\ref{prop:generic}(2) over that field) are integral over $R_\alpha$, and so is their sum, which lies in $k(\eta_\alpha)$. As $R_\alpha$ is integrally closed, $\operatorname{tr}_\alpha(a) \in R_\alpha$. The map $\operatorname{tr}_\alpha$ is $k[\h]$-linear and kills commutators, and $\operatorname{tr}_\alpha(e_i) = \alpha_i$ as $S_\alpha$ has dimension vector $\alpha$.

\emph{$\operatorname{tr}$ is evaluation.} Since $\alpha = \sum_j \alpha_j t_j$, $\operatorname{tr}_\alpha(e_i) = \alpha_i = \partial\alpha/\partial t_i$. Hence, as the isomorphism of (1) sends $[e_i]$ to $\partial/\partial t_i$ (Lemma~\ref{lem:KS}) and both maps are $k[\h]$-linear, $\operatorname{tr}$ is $\xi \mapsto (\xi(\alpha)|_{\alpha^\perp})_\alpha$, whose kernel is $\Der(-\log\mathcal H)$; so $\operatorname{tr}$ is injective. Finally, near a point of $\h \setminus \Sigma$ on $\alpha^\perp$, every other root is a unit and the constant vector field $\xi_0$ with $\xi_0(\alpha) = 1$ maps to $1$, so the map is onto there.
\end{proof}

\paragraph{Geometric motivation (a sketch).} The following heuristic lies behind Lemma~\ref{lem:KS}. It is not used in the proofs, and we do not make it precise. Through the derived equivalence of Proposition~\ref{prop:simres}(1) and the Hochschild--Kostant--Rosenberg theorem, $\HH_0(\Pi_\Gamma)$ corresponds to $H^1(\cY_\Gamma, \Omega^1_{\cY_\Gamma/\h})$, which the relative symplectic form identifies with $H^1(\cY_\Gamma, T_{\cY_\Gamma/\h})$. Then $\kappa$ becomes the Kodaira--Spencer map of $\chi\colon \cY_\Gamma \to \h$. It is onto because $\chi$ is a versal deformation of the minimal resolution, and its kernel is the vector fields on $\h$ that lift to $\cY_\Gamma$, which are those tangent to the root hyperplanes. The trace on $S_\alpha$ in Proposition~\ref{prop:wholeHH}(2) is the normal component of a vector field along $\alpha^\perp$, the hyperplane along which a single $(-2)$-curve survives.

\subsection{From the whole algebra to the corner}\label{ss:transfer}

In this subsection we show that $\HH_0(\cC_J)$ is a direct summand of $\HH_0(\Pi_\Gamma)$, an $H_J$-module through $H_J \subset k[\h]$, so that the Cohen--Macaulay property of \S\ref{ss:wholeHH} passes to the corner. We use a transfer: the inclusion of $A_J$ in $\cA_\Gamma$ induces a map $\operatorname{inc}_*$ on $\HH_0$, and the trace of left multiplication on $\cA_\Gamma e$, a perfect right $A_J$-module because $\cY_J$ is smooth (Lemma~\ref{lem:Ysmooth}), gives a map $\operatorname{Tr}$ the other way with $\operatorname{Tr} \circ \operatorname{inc}_* = |W_J|$. Both are compatible with killing $e_0$.

\begin{prop}\label{prop:transfer}
$\HH_0(\cC_J)$ is a direct summand of $\HH_0(\Pi_\Gamma)$ as an $H_J$-module.
\end{prop}

\begin{proof}
\emph{The two maps.} Let $r = |W_J|$ and $e = e_0 + e_\bul$. By Proposition~\ref{prop:tilting}, $A_J \subset A_J \otimes_{H_J} k[\h] \cong e\cA_\Gamma e \subset \cA_\Gamma$, compatibly with $e_0$, where the first map is injective as $H_J$ is a direct summand of $k[\h]$ (\S\ref{ss:lie}). This inclusion induces an $H_J$-linear map $\operatorname{inc}_*\colon \HH_0(A_J) \to \HH_0(\cA_\Gamma)$. Conversely, $\cA_\Gamma e$ is a finitely generated right $A_J$-module: $\cA_\Gamma$ is finite over $\cO(\cX_\Gamma \times_{\h/W} \h)$, as $\phi$ is projective (\S\ref{ss:simres}), and this ring acts on $\cA_\Gamma e$ through its image in $e\cA_\Gamma e \cong A_J \otimes_{H_J} k[\h]$, which is finite over $A_J$. As $\cY_J$ is smooth (Lemma~\ref{lem:Ysmooth}), $A_J$ has finite global dimension \cite[Cor~3.2.8 and the proof of Cor~3.2.11]{VdB}, so $\cA_\Gamma e$ is a perfect right $A_J$-module. Sending $a \in \cA_\Gamma$ to the trace of left multiplication by $a$ on $\cA_\Gamma e$ defines a map $\operatorname{Tr}\colon \HH_0(\cA_\Gamma) \to \HH_0(A_J)$. This trace is the alternating sum of the Hattori--Stallings traces of a lift of left multiplication by $a$ to a finite resolution by finitely generated projective $A_J$-modules. It is independent of the choices and vanishes on commutators \cite[\S1, Thm~2.2]{Sta}, \cite[Thm~1]{Han}. The Hattori--Stallings trace is additive on direct sums and linear over central elements, and $H_J$ acts on $\cA_\Gamma e$ through central elements on either side, so this map is $H_J$-linear.

\emph{$\operatorname{Tr} \circ \operatorname{inc}_* = r$.} For $x \in A_J$, $x = exe$ kills $(1 - e)\cA_\Gamma e$ and acts on $e\cA_\Gamma e \cong A_J \otimes_{H_J} k[\h] \cong A_J^{\oplus r}$ (Chevalley, \S\ref{ss:lie}) by left multiplication on each summand; so $\operatorname{Tr} \circ \operatorname{inc}_* = r$.

\emph{Killing $e_0$.} Both maps preserve the images of $\HH_0(e_0(-)e_0)$: for $\operatorname{inc}_*$ because $e_0A_Je_0 \subset e_0\cA_\Gamma e_0$, and for $\operatorname{Tr}$ because $a \in e_0\cA_\Gamma e_0$ kills $(1 - e_0)\cA_\Gamma e$ and acts on $e_0\cA_\Gamma e \cong (e_0A_J)^{\oplus r}$ by an $r \times r$ matrix over $e_0A_Je_0$, whose trace is the sum of its diagonal entries. By \eqref{eq:idem} and Lemma~\ref{lem:corner}, $\operatorname{inc}_*$ and $\operatorname{Tr}$ descend to $H_J$-linear maps $\HH_0(\cC_J) \to \HH_0(\Pi_\Gamma) \to \HH_0(\cC_J)$ whose composite is $r$, so $\HH_0(\cC_J)$ is a direct summand of $\HH_0(\Pi_\Gamma)$. We do not claim that $\operatorname{Tr}$ is homogeneous; Theorem~\ref{thm:squeeze} uses only the ungraded statement, as depth is computed by $\Ext$.
\end{proof}

\subsection{\texorpdfstring{Freeness of $\HH_0$}{Freeness of HH0}}\label{ss:proofA3}

We now put \S\S\ref{ss:wholeHH} and~\ref{ss:transfer} together to deduce that $\HH_0(\cC_J)$ is free over $P$.

\begin{thm}\label{thm:squeeze}
Let $\Gamma$ be a Dynkin diagram and $\bul$ any vertex. Then $\HH_0(\cC_J)$ is a free $P$-module of rank $|\Phi^+_\bul|$, and $\HH_0(\Pi_\bul)$ is a graded Cohen--Macaulay $k[\h]$-module of dimension $n - 1$, free over $k[\h_J]$ of the same rank.
\end{thm}

\begin{proof}
By Lemma~\ref{lem:freeness}(2), depth is unchanged by passing to a subring over which the module is finite. Hence $\HH_0(\Pi_\Gamma)$, which is finite over $H_J \subset k[\h]$, has depth $n - 1$ over $H_J$ (Proposition~\ref{prop:wholeHH}(1)). Its direct summand $\HH_0(\cC_J)$ (Proposition~\ref{prop:transfer}) has depth at least $n - 1$, as depth is computed by $\Ext$ and $\Ext$ is additive. As $\HH_0(\cC_J)$ is finite over $P$ (Corollary~\ref{cor:lower}), its depth over $P$ is at least $n - 1 = \dim P$, so it is free over $P$ by Lemma~\ref{lem:freeness}(2). Its rank is $|\Phi^+_\bul|$ by Corollary~\ref{cor:trace-generic}. Then $\HH_0(\Pi_\bul) = \HH_0(\cC_J) \otimes_P k[\h_J]$ (Lemma~\ref{lem:corner}, as $\HH_0$ commutes with base change) is free over $k[\h_J]$ (\S\ref{ss:lie}), hence Cohen--Macaulay of dimension $n - 1$ (Lemma~\ref{lem:freeness}(2)).
\end{proof}

\subsection{The universal contraction algebra}\label{ss:proofA}

We now collect the results of Sections~\ref{sec:cartan} and~\ref{sec:squeeze} into one statement at an arbitrary vertex, and specialise it to the universal contraction algebra.

\begin{thm}[Theorem~A at every vertex]\label{thm:Avertex}
Let $\Gamma$ be a Dynkin diagram and $\bul$ any vertex, with $J$, $H_J = P[t]$, $l$, $f_j$ and $D_j$ as in \S\ref{ss:lie}, and $\cC_J$ as in \S\ref{ss:partial}.
\begin{thmlist}
\item $\cC_J \otimes_{H_J} k[\h] \cong \Pi_\bul = e_\bul\Pi^\lambda(\Gamma)e_\bul$.
\item $\cC_J$ is supported on $D_1 \cup \dots \cup D_l$, and at a general point of $D_j$ its fibre is $M_j(k)$.
\item $\cC_J$ and $\HH_0(\cC_J)$ are free over $P$, of ranks $\sum_{\alpha > 0} \alpha_\bul^2$ and $|\Phi^+_\bul|$, and the characteristic polynomials of multiplication by $t$ on them are $\prod_{j=1}^{l} f_j^{\,j^2}$ and $f_1 f_2 \cdots f_l$ respectively, up to nonzero scalars. Hence each is the cokernel of a square matrix over $H_J$ with this determinant.
\end{thmlist}
\end{thm}

\begin{proof}
Part (1) is Lemma~\ref{lem:corner}, and part (2) is Proposition~\ref{prop:support}. In part (3), the freeness and ranks are Proposition~\ref{prop:ER} and Theorem~\ref{thm:squeeze}. The characteristic polynomials and presentations follow from Lemma~\ref{lem:presentation}: for $\cC_J$ with $m_j = j^2$, by Proposition~\ref{prop:ER} and Corollary~\ref{cor:trace-generic}, and for $\HH_0(\cC_J)$ with $m_j = 1$, as $\HH_0(\Pi_\bul)$ has length $q_\alpha = 1$ at $\eta_\alpha$ (Corollary~\ref{cor:trace-generic}) and is supported on the union of the $\alpha^\perp$ with $\alpha_\bul \neq 0$ (Proposition~\ref{prop:generic}).
\end{proof}

\begin{proof}[Proof of Theorem~A]
For the six pairs of Table~\ref{tab:roots}, $l$ is the length, $H_J = H_l$ and $\cC_J = \cC_l$ is the universal contraction algebra (\S\ref{ss:universal}). Hence Theorem~A is Theorem~\ref{thm:Avertex} at the flopping vertex.
\end{proof}

\begin{rem}[Central fibres and computer checks]\label{rem:noER}
For the six pairs, a computer calculation gives a second proof of Theorem~\ref{thm:squeeze} and Proposition~\ref{prop:ER}, independent of \S\S\ref{ss:wholeHH}--\ref{ss:transfer} and of \cite{ER}. Let $\overline{\cC}_l = \cC_l/P_+\cC_l$ be the central fibre of $\cC_l$ over the origin of $\Spec P$, where $P_+$ is the positive-degree part of $P$. It is a finite-dimensional graded algebra, the algebra of \cite[Thm~1.3, Examples~4.23--4.27]{K} with $e_0$ killed and every parameter $T$ set to zero (Remark~\ref{rem:Kpresentation}). For trivalent $\bul$, $\overline{\cC}_l \cong k[y]\langle b, c\rangle/(b^p, c^q, (y - b - c)^r)$, with $y$ central and $p, q, r$ the sizes of the arms \cite[Cor~3.4]{ER}, a presentation that the computation below does not use. As $\HH_0$ commutes with base change, Corollary~\ref{cor:lower} and Lemma~\ref{lem:freeness}(1) show that $\HH_0(\cC_l)$ is free if (and, by Corollary~\ref{cor:trace-generic}, only if) $\dim \HH_0(\overline{\cC}_l) = |\Phi^+_\bul|$, that is, $1, 9, 29, 53, 104, 106$ for $l = 1, \dots, 6$, and $\cC_l$ is free of rank $\sum_{\alpha>0}\alpha_\bul^2$ if $\dim \overline{\cC}_l = 1, 12, 72, 216, 600, 900$. We have carried out this computation for the six pairs. On the line $k\varpi_\bul^\vee \subset \h$, on which $P_+$ vanishes and $t$ is a coordinate, $\Pi_\bul$ restricts to $\overline{\cC}_l$ (Lemma~\ref{lem:corner}). We computed this algebra and its $\HH_0$ degree by degree from the relations of $\Pi^\lambda(\Gamma)$ on this line, over $\mathbb F_{1000003}$. We stopped once $e_\bul$ times the restriction of $\Pi^\lambda(\Gamma)$ to this line vanished in two consecutive degrees. This suffices, as that restriction is generated over its degree-$0$ part in degrees $1$ and~$2$. Each graded piece of the corresponding algebra over $\Z$ is a finitely generated abelian group, so the dimensions over $\mathbb F_{1000003}$ bound those over $\Q$, and hence over $k$, from above. Corollary~\ref{cor:lower} and Lemma~\ref{lem:freeness}(1) bound them from below, and both bounds equal the values above. This remark is not used elsewhere in the paper.
\end{rem}

\subsection{The centre of the corner}\label{ss:centre}

We return to an arbitrary vertex $\bul$ of a Dynkin diagram, and show that the centre of the corner is free as well, of the same rank as its zeroth Hochschild homology. This is the universal analogue of the equality $\dim Z(\Acon) = \dim \HH_0(\Acon)$ in Theorem~B(1), which for a single flop follows from the symmetry of $\Acon$. Corollary~\ref{cor:centre} is not used elsewhere in the paper.

\begin{cor}[The centre of the corner]\label{cor:centre}
Let $\Gamma$ be a Dynkin diagram and $\bul$ any vertex, as in Sections~\ref{sec:cartan} and~\ref{sec:squeeze}. Then $Z(\cC_J)$ is a graded free $P$-module of rank $|\Phi^+_\bul|$, and $Z(\Pi_\bul) = Z(\cC_J) \otimes_P k[\h_J]$ is free over $k[\h_J]$ of the same rank. In particular the centre of the universal contraction algebra $\cC_l$ is free over $P$, of the same rank $1, 9, 29, 53, 104, 106$ as $\HH_0(\cC_l)$, for $l = 1, \dots, 6$.
\end{cor}

\begin{proof}
It suffices to show that $Z(\Pi_\bul)$ is free over $k[\h_J]$ of rank $|\Phi^+_\bul|$: as $\Pi_\bul = \cC_J \otimes_P k[\h_J]$ (Lemma~\ref{lem:corner}) and the centre commutes with the flat base change $P \to k[\h_J]$, $Z(\Pi_\bul) = Z(\cC_J) \otimes_P k[\h_J]$, so then $Z(\cC_J)$ is free over $P$ of the same rank (\S\ref{ss:lie}).

The idea is to restrict to the line $k\rho^\vee$ defined next, which meets the root hyperplanes only at $0$. There $\Pi_\Gamma$ becomes the central extension of \cite{ELR}, a symmetric algebra, so Lemma~\ref{lem:symfamily} gives the freeness of the centre near the line, and graded Nakayama lifts it to $k[\h_J]$. Let $\rho^\vee \in \h$ be the point with $t_i(\rho^\vee) = 1$ for all $i \in \Gamma_0$ (half the sum of the positive coroots). A positive root takes at $\rho^\vee$ the sum of its coefficients, which is positive, so the line $k\rho^\vee$ meets each root hyperplane only at $0$. Let $\Lambda \subset k[\h]$ be the polynomial ring on the linear forms vanishing on $k\rho^\vee$, so that $k[\h]/\Lambda_+k[\h]$ is the coordinate ring of this line. The supports of $\Pi_\bul$ and $\HH_0(\Pi_\bul)$ lie in the union of the root hyperplanes (Proposition~\ref{prop:generic}), so modulo $\Lambda_+$ both are finite-dimensional, and by graded Nakayama both are finite over $\Lambda$. As they are Cohen--Macaulay of dimension $n - 1$ (Proposition~\ref{prop:ER}, Theorem~\ref{thm:squeeze}), they are free over $\Lambda$ (Lemma~\ref{lem:freeness}(2)).

Modulo $\Lambda_+$, $\Pi_\bul$ becomes the corner $e_\bul\Pi_\rho e_\bul$ of the restriction $\Pi_\rho = \Pi_\Gamma \otimes_{k[\h]} k[\h]/\Lambda_+k[\h]$ of $\Pi_\Gamma$ to the line, which is finite-dimensional for the same reason. As $t_i(\rho^\vee) = 1$ for every $i \in \Gamma_0$, at the points $\lambda = s\rho^\vee$ of the line every parameter $\lambda_i$ equals $s$, so the defining relation of $\Pi^\lambda(\Gamma)$ becomes $\sum_a [a, a^*] = s$. Thus $\Pi_\rho$ is the central extension of the preprojective algebra of $\Gamma$ for the weight $\rho$, with central variable $s$, and it is a symmetric algebra \cite[\S2.1, Prop~2.1]{ELR}. Hence so is $e_\bul\Pi_\rho e_\bul$, as a symmetrising form restricts to one on a corner.

Let $O$ be the local ring of $\Lambda$ at $\Lambda_+$. Then $\Pi_\bul \otimes_\Lambda O$ and $\HH_0(\Pi_\bul \otimes_\Lambda O) = \HH_0(\Pi_\bul) \otimes_\Lambda O$ are free over $O$, and $\Pi_\bul \otimes_\Lambda k = e_\bul\Pi_\rho e_\bul$ is symmetric, so $Z(\Pi_\bul \otimes_\Lambda O)$ is free over $O$ (Lemma~\ref{lem:symfamily}). The centre commutes with flat base change, so $Z(\Pi_\bul) \otimes_\Lambda O$ is free over $O$. As $Z(\Pi_\bul)$ is a finite graded $\Lambda$-module, like $\Pi_\bul$, it is free over $\Lambda$ by Lemma~\ref{lem:freeness}(1), since localising at $\Lambda_+$ changes neither its rank nor its quotient by $\Lambda_+$. Hence $Z(\Pi_\bul)$ is a Cohen--Macaulay $k[\h]$-module of dimension $n - 1$. It is finite over $k[\h_J]$, as $\Pi_\bul$ is (Proposition~\ref{prop:ER}), and hence free over $k[\h_J]$ (Lemma~\ref{lem:freeness}(2)).

Finally, $Z(\Pi_\bul) \subset \Pi_\bul$ is supported on the hyperplanes $\alpha^\perp$ with $\alpha_\bul \neq 0$, each of which maps isomorphically onto $\h_J$ (\S\ref{ss:lie}), so its rank over $k[\h_J]$ is the sum of its lengths at the generic points $\eta_\alpha$. At $\eta_\alpha$, $\Pi_\bul$ is an algebra over the residue field $k(\eta_\alpha)$ (proof of Corollary~\ref{cor:trace-generic}) that becomes $M_{\alpha_\bul}(\overline{k(\eta_\alpha)})$ over an algebraic closure (Proposition~\ref{prop:generic}(2)), so its centre is $k(\eta_\alpha)$, of length $1$, and the rank is $|\Phi^+_\bul|$.
\end{proof}

\begin{rem}[Central extensions]\label{rem:ELR}
Proposition~\ref{prop:wholeHH}(1) also computes $\HH_0$ of the central extension $\Pi_\rho$ in the proof above. With $p_1, \dots, p_n$ as in Lemma~\ref{lem:logbasis}, let $m_i = \deg p_i - 1$, the exponents of $W$, and let $M$ be the matrix whose $i$th column is $\operatorname{grad} p_i$ in a basis of the constant vector fields, so that $\Theta_\h/\Der(-\log\mathcal H)$ is the cokernel of $M$. As $\HH_0$ commutes with base change, $\HH_0(\Pi_\rho)$ is the cokernel of $M(s\rho^\vee)$ over $k[s]$. The entries of the $i$th column of $M$ are homogeneous of degree $m_i$, so $M(s\rho^\vee) = M(\rho^\vee)\operatorname{diag}(s^{m_1}, \dots, s^{m_n})$, and $\det M(\rho^\vee) = c\prod_{\alpha>0}\alpha(\rho^\vee) \neq 0$ (proof of Lemma~\ref{lem:logbasis}). Hence $\HH_0(\Pi_\rho) \cong \bigoplus_{i=1}^n k[s]/(s^{m_i})$: $s$ acts on it with Jordan blocks of sizes $m_1, \dots, m_n$, and its dimension is $\sum_i m_i = |\Phi^+|$, the degree of $\det M$. As $\Pi_\rho$ is symmetric, $Z(\Pi_\rho)$ is dual to $\HH_0(\Pi_\rho)$, with $s$ self-adjoint \cite[Cor~2.2, \S2.3]{ELR}, so the same holds for the centre. The computation of $\HH_0$ applies to any point of $\h$ off the root hyperplanes in place of $\rho^\vee$. This recovers \cite[Thm~3.1]{ELR}, for $\rho$ and for generic weights, which is proved there case by case, with a computer calculation in types $\mathsf{E}_6$--$\mathsf{E}_8$.
\end{rem}

\section{\texorpdfstring{Flops (Theorems~B and~C)}{Flops (Theorems B and C)}}\label{sec:flops}

In this section we pull Theorem~A back from the universal flop to each individual flop, along its classifying arc, and prove Theorems~B and~C. In Figure~\ref{fig:story} this is the passage from the middle column to the right-hand one. The geometric input is the versality of the universal flop (Lemma~\ref{lem:classifying}) and the curves in its fibres (Proposition~\ref{prop:BKL}). Corollary~\ref{cor:toda} and Theorem~B(2) recover results of Toda and of Hua and Toda. In \S\ref{ss:basechange} we first show that the contraction algebra of a flop is the universal one pulled back along the arc. Then in \S\ref{ss:dims} the presentations of Theorem~A(3) compute $\dim \Acon$ and $\dim \HH_0(\Acon)$ as orders of vanishing of the $f_j$ along the arc. In \S\ref{ss:GV} we show that these orders are the GV invariants. In \S\ref{ss:proofB} we combine these to prove Theorem~B(1), (2), in \S\ref{ss:defflop} we prove Theorem~B(3), on deformations of the flop, and in \S\ref{ss:proofC} we deduce Theorem~C. In Section~\ref{sec:examples} we express the same numbers as sums over roots, and follow explicit flops through these steps, including the Laufer flop in Example~\ref{ex:lauferarc}.

Fix a flop $f\colon Y \to X = \Spec R$ of length $l$, with $R$ complete local isolated cDV.

\subsection{Base change and classifying arcs}\label{ss:basechange}

In order to identify the contraction algebra of a flop with the pullback of the universal one, we show that the tilting bundle of \cite{K} on the universal flop pulls back to Van den Bergh's tilting bundle on the flop. We state the base change over a more general base, as Lemma~\ref{lem:truncation} and the proof of Theorem~B(3) use it for other maps to $\Spec H_l$. We then show that every flop arises in this way (Lemma~\ref{lem:arcs}), and say which arcs classify flops (Lemma~\ref{lem:arcflop}).

\begin{lem}[Base change]\label{lem:basechange}
Let $B$ be $k[[w]]$ or $k[[w]] \mathbin{\hat\otimes} O$, with $O$ a complete local noetherian $k$-algebra with residue field $k$, and let $\nu\colon \Spec B \to \Spec H_l$ send the closed point to the origin. Let $X = \Spec R$ be the completion of $\cX_l \times_{H_l} \Spec B$ at the singular point of its closed fibre, and $f\colon Y \to X$ the pullback of $\cY_l$. Then the pullback $\mathcal T_X$ of the universal tilting bundle $\cO \oplus \mathcal M$ of Proposition~\ref{prop:tilting} is a tilting bundle on $Y$, isomorphic to $\cO_Y \oplus N_1$, where $\cO_Y$ and $N_1$ are the indecomposable projective objects of ${}^0\Per(Y/X)$, and
\[
\End_Y(\mathcal T_X)/(e_0) \;\cong\; \cC_\nu := \cC_l \otimes_{H_l} B ,
\]
with $B$ an $H_l$-algebra through $\nu$ and $e_0$ the idempotent of $\cO_Y$. In particular, if $B = k[[w]]$ and $f$ is a flop, then $\mathcal T_X$ is the basic projective generator of \cite[\S2.4, Def~2.12]{DW} and $\Acon \cong \cC_\nu$.
\end{lem}

\begin{proof}
Let $g\colon \cY_l \to \cX_l = \Spec R_l$ be the universal flopping contraction (\S\ref{ss:universal}) and $\mathcal T = \cO \oplus \mathcal M$ the universal tilting bundle, so that $\End(\mathcal T) = A_l$ and $\Ext^i(\mathcal T, \mathcal T) = 0$ for $i > 0$. We argue in five steps. In Step~1 we place $f$ in the setting of Theorem~\ref{thm:VdB}, and in Step~2 we compute the endomorphism algebra of $\mathcal T_X$ and show that its higher self-extensions vanish. In Step~3 we show that $\mathcal T_X$ generates, so is a tilting bundle, and that it is a projective generator of ${}^0\Per(Y/X)$, isomorphic to $\cO_Y \oplus N_1$. Then in Step~4 we show that $\cC_\nu$ is complete, and in Step~5 we identify the quotient. Completion at the singular point is flat, so the base-change statements of Steps~1--3 hold after it, and we omit it from the notation until Step~4.

\emph{Step 1 (Van den Bergh's setting).} As $\cY_l$ and $\cX_l$ are flat over $H_l$ (\S\ref{ss:universal}), Lemma~\ref{lem:flatfamilies}(1) applies to $g$ and the base change $\nu\colon \Spec B \to \Spec H_l$. With $Rg_*\cO_{\cY_l} = \cO_{\cX_l}$ it gives $Rf_*\cO_Y = \cO_X$, as $\cO_X = R_l \otimes_{H_l} B$ before completion; and $f$ is projective with fibres of dimension at most one, being a pullback of $g$. Both persist after completion. Hence Theorem~\ref{thm:VdB} applies to $f$.

\emph{Step 2 (endomorphisms and self-extensions).} Let $\mathcal T_X = \cO_Y \oplus \mathcal M_X$ be the pullback of $\mathcal T$ to $Y$. As $\cX_l$ is affine and $\Ext^i(\mathcal T, \mathcal T) = 0$ for $i > 0$, the vector bundle $\sEnd(\mathcal T)$ has $R^ig_*\sEnd(\mathcal T) = 0$ for $i > 0$ and $g_*\sEnd(\mathcal T) = A_l$; so by Lemma~\ref{lem:flatfamilies}(1), $R^if_*\sEnd(\mathcal T_X) = 0$ for $i > 0$ and $f_*\sEnd(\mathcal T_X) \cong A_l \otimes_{H_l} B$. As $X$ is affine, $\Ext^i_Y(\mathcal T_X, \mathcal T_X) = 0$ for $i > 0$ and $\End_Y(\mathcal T_X) \cong A_l \otimes_{H_l} B$.

\emph{Step 3 (generation, ${}^0\Per$ and basicness).} The projection $p\colon Y \to \cY_l$ is affine. If $\mathcal F \in D_{\mathrm{qc}}(Y)$ satisfies $R\Hom(p^*\mathcal T, \mathcal F) = 0$, then by adjunction $R\Hom(\mathcal T, Rp_*\mathcal F) = 0$, so $Rp_*\mathcal F = 0$ because $\mathcal T$ generates $D_{\mathrm{qc}}(\cY_l)$, and hence $\mathcal F = 0$ because affine pushforward is conservative. Hence $\mathcal T_X = p^*\mathcal T$ is a tilting bundle.

To identify the contraction algebra one must also check that $\mathcal T_X$ is a projective generator of ${}^0\Per(Y/X)$. Put $\mathcal N = \mathcal M^\vee$, and $\mathcal N_Y = \mathcal M_X^\vee$ for its pullback to $Y$. By Proposition~\ref{prop:tilting}, $\cO \oplus \mathcal M$ is a projective generator of ${}^{0}\Per(\cY_l/\cX_l)$. By Van den Bergh's characterisation of projective generators (Theorem~\ref{thm:VdB}(2)), its dual $\cO \oplus \mathcal N$ is one of ${}^{-1}\Per(\cY_l/\cX_l)$, so $\mathcal N$ is globally generated, $H^1(\cY_l, \mathcal N^\vee) = 0$, and $\det \mathcal N$ is relatively ample. On the fibre over the origin of $\Spec H_l$ it has degree one on the curve. Global generation pulls back, as the evaluation surjection does. Next, $H^1(Y, \mathcal N_Y^\vee) = 0$: as $R^ig_*\mathcal N^\vee = 0$ for $i > 0$ (the fibres of $g$ have dimension at most one and $H^1(\mathcal N^\vee) = 0$), Lemma~\ref{lem:flatfamilies}(1) gives $R^1f_*\mathcal N_Y^\vee = 0$, and $X$ is affine. Finally $\det \mathcal N_Y$ has degree one on the curve $C$ over the closed point, so it is relatively ample after completion (Theorem~\ref{thm:VdB}(3)). Hence $\cO \oplus \mathcal N_Y$ is a projective generator of ${}^{-1}\Per(Y/X)$, and its dual $\mathcal T_X$ is one of ${}^{0}\Per(Y/X)$.

Finally, $\mathcal T_X \cong \cO_Y \oplus N_1$, with $N_1$ the dual of Van den Bergh's extension (the bundle $M$ of \S\ref{ss:flops} when $f$ is a flop). Indeed, as $X$ is complete local, every projective object of ${}^0\Per(Y/X)$ is $\cO_Y^{\oplus a} \oplus N_1^{\oplus b}$ (Theorem~\ref{thm:VdB}(3)). For $\mathcal T_X$ we have $a \ge 1$. Also $b \ge 1$, since otherwise $\mathcal T_X$ would be a sum of copies of $\cO_Y$, which does not generate, as $R\Hom(\cO_Y, \cO_C(-1)) = R\Gamma(C, \cO_C(-1)) = 0$. Comparing ranks, $a + lb = l + 1$: $N_1$ has rank the multiplicity of $C$ in the fibre of $f$ over the closed point (Theorem~\ref{thm:VdB}(3)). As $\nu$ sends the closed point to the origin, this fibre is that of $Y_0 \to X_0$ over its singular point, in which $C$ has multiplicity $l$ (\S\ref{ss:universal}). Hence $a = b = 1$.

\emph{Step 4 (completeness).} $\cC_l$ is finite over $H_l$ (Theorem~A(3)), so $\cC_\nu$ is a finite $B$-module, and hence complete for the adic topology of the maximal ideal $\mathfrak m_B$ of $B$. Its fibre $\cC_\nu \otimes_B k$, the fibre of $\cC_l$ over the origin of $\Spec H_l$, is finite-dimensional and supported at the singular point of the fibre of $\cX_l$ there. Hence a power of the ideal of the singular point maps $\cC_\nu$ into $\mathfrak m_B\cC_\nu$. As $\mathfrak m_B$ lies in that ideal, the topology of the singular point on $\cC_\nu$ is the $\mathfrak m_B$-adic one, and $\cC_\nu$ is complete for it.

\emph{Step 5 (the quotient).} By Step~2, $\End_Y(\mathcal T_X)$ is the completion of $A_l \otimes_{H_l} B$ at the singular point, compatibly with $e_0$. Since killing $e_0$ commutes with $\otimes$ and with completion, $\End_Y(\mathcal T_X)/(e_0)$ is the completion of $\cC_\nu$, which is $\cC_\nu$ itself by Step~4. If $B = k[[w]]$ and $f$ is a flop, $\End_Y(\mathcal T_X) \cong \End_R(R \oplus f_*\mathcal M_X)$ with $\mathcal T_X \cong \cO_Y \oplus N_1$ is the algebra from which Donovan and Wemyss define $\Acon$, so $\Acon \cong \cC_\nu$.
\end{proof}

\begin{lem}[Classifying arcs]\label{lem:arcs}
For every flop $f\colon Y \to X$ of length $l$ there is a map $\nu\colon \Spec k[[w]] \to \Spec H_l$, sending the closed point to the origin, from which $f$ arises as in Lemma~\ref{lem:basechange}, so that $\Acon \cong \cC_\nu$. Moreover $\nu$ meets $D = D_1 \cup \dots \cup D_l$ only at the closed point, so each $f_j(\nu)$ is nonzero.
\end{lem}

\begin{proof}
Let $X_0$ be a generic hyperplane section through the singular point and $Y_0$ its proper transform, a partial resolution of the Kleinian singularity of type $\Gamma$ contracting the curves of $J$ (the completion of $Y_0^J$). With $w$ the equation of $X_0$, $Y$ is flat over $k[[w]]$, as $w$ is not a zero divisor on $\cO_Y$, and $Rf_*\cO_Y = \cO_X$. Hence, by Lemma~\ref{lem:classifying} with $B = k[[w]]$, $Y \to X$ is the pullback of the universal flop along a classifying arc $\nu\colon \Spec k[[w]] \to \Spec H_l$ sending the closed point to the origin, completed at the singular point as in Lemma~\ref{lem:basechange}. By Lemma~\ref{lem:basechange}, $\Acon \cong \cC_\nu$. As $\Acon$ is finite-dimensional (Proposition~\ref{prop:Acon}(1)), $\cC_\nu$ has finite length, and as the support of $\cC_l$ is $D$ (Theorem~A(2)), $\nu$ meets $D$ only at the closed point.
\end{proof}

\paragraph{Arcs that classify flops.} Conversely, we now say which arcs in the base of the universal flop classify flops. Let $\nu\colon \Spec k[[w]] \to \Spec H_l$ send the closed point to the origin, and let $f\colon Y \to X$ be the pullback of $g\colon \cY_l \to \cX_l$ along $\nu$. The partial resolution $Y_0 \to X_0$ contracts a single curve $C$, the curve of $\bul$, and the singular points of $Y_0$ lie on $C$, one $\mathsf{A}_{m-1}$ point for each arm of $J$ of type $\mathsf{A}_{m-1}$. At such a point, the completed local ring of $\cY_l$ is pulled back from the semi-universal deformation $x^2 + y^2 + z^m + u_{m-2}z^{m-2} + \dots + u_0 = 0$ of the $\mathsf{A}_{m-1}$ singularity (the family $F$ of \S\ref{ss:kleinian} for type $\mathsf{A}_{m-1}$) along a formally smooth map of bases (proof of Lemma~\ref{lem:Ysmooth}), through which the $u_i$ become elements of the completion of $H_l$ at the origin, vanishing there.

\begin{lem}[Arcs that classify flops]\label{lem:arcflop}
In this notation:
\begin{thmlist}
\item $Y$ is smooth along $C$ if and only if $\ord_w u_0(\nu) = 1$ at each singular point of $Y_0$. In particular this depends only on $\nu$ modulo $w^2$.
\item Suppose that $Y$ is smooth along $C$ and that $f_j(\nu) \neq 0$ for every $j$, that is, $\nu$ meets $D = D_1 \cup \dots \cup D_l$ only at the closed point. Then $\nu$ classifies a flop: completing $Y \to X$ at the singular point of its closed fibre gives a flop in the sense of \S\ref{ss:flops}.
\end{thmlist}
\end{lem}

\begin{proof}
(1) As $Y_0$ is a Cartier divisor in $Y$, $Y$ is smooth along $C$ except possibly at the singular points of $Y_0$. At each of them, the completed local ring of $Y$ is pulled back along $\nu$ from the local model above, so by the Jacobian criterion it is smooth there if and only if $\ord_w u_0(\nu) = 1$.

(2) After completion every closed point of $Y$ lies on $C$, so $Y$ is smooth. The singularity $X$ is a hypersurface, as $\cX_l$ is, hence Gorenstein, and it is cDV in Reid's sense, as its section $w = 0$ is the rational double point $X_0$. It is singular: the constant term of the equation $F$ of $\cX_\Gamma$ (\S\ref{ss:kleinian}) has degree $h \ge 2$ in $k[\h]^W \subseteq H_l$, while $H_l$ is generated in degrees less than $h$ (by $t$ and the elementary symmetric polynomials of the arms; see ``Generators of $H_l$'' in \S\ref{ss:universal}). Hence this constant term lies in the square of the ideal of the origin, so its order along $\nu$ is at least $2$, and every partial derivative of the equation of $X$ vanishes at the closed point. Next, $Rf_*\cO_Y = \cO_X$ by base change along $\nu$ (Lemma~\ref{lem:flatfamilies}(1)). The exceptional locus of $f$ is $C$. Indeed, $f$ contracts only $C$ on $Y_0$, and it is an isomorphism over the generic point of $\Spec k[[w]]$, where the arc misses $D$ and so the fibres of $g$ contain no complete curve. The points of $\Spec H_l$ over which the fibres of $g$ do contain a complete curve form a constructible set, the image of the exceptional locus of $g$. Its closed points lie in $D$ (Proposition~\ref{prop:BKL}(1)), so it lies in $D$, as $\Spec H_l$ is a Jacobson scheme. Hence $f$ is an isomorphism over the complement of the closed point of $X$: there it is proper and quasi-finite, hence finite, and $f_*\cO_Y = \cO_X$. As $Y$ is smooth, $X$ is therefore an isolated singularity. Finally $f$ is crepant, being a small morphism to the Gorenstein $X$.
\end{proof}

\subsection{\texorpdfstring{Dimensions from the presentations}{Dimensions from the presentations}}\label{ss:dims}

Pulling back the presentations of Theorem~A(3) along the arc, we now compute the dimensions of $\Acon$ and $\HH_0(\Acon)$ as orders of vanishing.

\begin{cor}\label{cor:rootsum}
Let $\nu$ be the classifying arc of Lemma~\ref{lem:arcs}. Then
\[
\dim \HH_0(\Acon) = \ord_w (f_1 \cdots f_l)(\nu), \qquad \dim \Acon = \sum_{j} j^2 \ord_w f_j(\nu).
\]
\end{cor}

\begin{proof}
$\HH_0$ commutes with base change, so by Lemma~\ref{lem:arcs}, $\HH_0(\Acon) = \HH_0(\cC_l) \otimes_{H_l} k[[w]]$. By Theorem~A(3), $\HH_0(\cC_l)$ is the cokernel of a square matrix over $H_l$ with determinant $c'f_1 \cdots f_l$, $c' \in k^\times$; pulling it back along $\nu$ gives a square matrix over $k[[w]]$ with cokernel $\HH_0(\Acon)$ and determinant $c'\,(f_1 \cdots f_l)(\nu)$, which is nonzero because $\nu$ meets the support of $\cC_l$ only at the closed point (Lemma~\ref{lem:arcs}). By Lemma~\ref{lem:dvr}(2), $\HH_0(\Acon)$ has length $\ord_w (f_1 \cdots f_l)(\nu)$. The same argument applies to $\cC_l$, with determinant $c \prod_j f_j^{j^2}$, $c \in k^\times$.
\end{proof}

\subsection{Identifying the GV invariants}\label{ss:GV}

We now show that these orders of vanishing are the GV invariants (Proposition~\ref{prop:GV}). This is the step where the geometry of the curves in the universal family enters. The count of \cite{BKL} (Proposition~\ref{prop:BKL}(2)) needs a polynomial arc, so we first compare an arc with its truncations.

\begin{lem}[Truncation]\label{lem:truncation}
Let $\nu$ be the classifying arc of Lemma~\ref{lem:arcs}, and for $r \ge 1$ let $\nu_r$ be its truncation modulo $w^r$, a polynomial map. For $r$ large, each $f_j$ has the same order of vanishing along $\nu_r$ as along $\nu$, and $\nu_r$ classifies a flop with the same contraction algebra as $f$. Let $n_1', \dots, n_l'$ be the GV invariants of this truncated flop. If moreover $r \ge j(n_j' + 1)$ for each $j$, then $f$ has the same GV invariants.
\end{lem}

\begin{proof}
\emph{Orders.} Each $\ord_w f_j(\nu)$ is finite (Lemma~\ref{lem:arcs}), so it does not change once $r$ exceeds it.

\emph{A flop.} The pullback of $\cY_l$ along $\nu_r$ is still smooth along $C$, as this depends only on the arc modulo $w^2$ (Lemma~\ref{lem:arcflop}(1)). By the orders above, $\nu_r$ meets $D$ only at the closed point. Hence $\nu_r$ classifies a flop (Lemma~\ref{lem:arcflop}(2)).

\emph{The curves $C_j$.} In $f$, the curve $C_j$ is cut out in the partial resolution $Y_0 = \{w = 0\}$ of type $(\Gamma, \bul)$, the proper transform of a generic hyperplane section (Lemma~\ref{lem:arcs}), which is also the section $w = 0$ of the truncated flop. The surface determines $C_j$, as its unique one-dimensional subscheme supported on $C$ without embedded points and of multiplicity $j$ at the generic point \cite[Lemma~3.2(i)]{K08}.

\emph{Katz's argument for the section $w = 0$.} We check that the GV invariants $n_j'$ of the truncated flop are the lengths of its Hilbert scheme at these $[C_j]$, although $w = 0$ need not be a generic hyperplane section of its base. The proof of \cite[Prop~3.3]{K08} uses the generic section only as a Cartier divisor through $C$ whose proper transform is a partial resolution with rational double points. A stable sheaf of class $j[C]$ and Euler characteristic $1$ is then a module over that surface, so the kernel of a section of it cuts out a one-dimensional subscheme of the surface supported on $C$, and by the same lemma and stability the sheaf is $\cO_{C_j}$. This also identifies $C_j$ with the curve that \S\ref{ss:flops} defines through a generic hyperplane section of the truncated flop, which is a flop. The structure sheaf of that curve is stable of class $j[C]$ with Euler characteristic $1$, by the same argument for that section, hence isomorphic to $\cO_{C_j}$, and isomorphic cyclic $\cO_Y$-modules have the same annihilator, so the two subschemes coincide. (Katz also intersects over all such sections, to place the support in the fibre of $f$ over the closed point of $X$; this argument does not need that step.) As $Y_0$ is the fibre of $w$ in the truncated flop, that proof applies to it with the same curves $C_j$. The truncated flop is the completion of an algebraic flop, so this gives Katz's results for it: $[C_j]$ is an isolated point of its Hilbert scheme, and the complete local ring $O_r$ there has length $n_j'$. Nothing here uses that the arc is a truncation, so the same holds for the flop classified by any polynomial arc.

\emph{GV invariants: dependence on the truncation.} A deformation of $C_j$ over an Artinian local base of length $s$ lies in $Y \times_X \Spec R/\mathfrak m_R^{\,js}$, as its structure sheaf has a filtration of length $s$ whose quotients are $\cO_{C_j}$-modules, and $\mathfrak m_R^{\,j}\cO_Y \subseteq I_C^{\,j} \subseteq I_{C_j}$ for the ideal sheaves $I_C$ and $I_{C_j}$. As $w \in \mathfrak m_R$, this subscheme of $Y$ depends only on $\nu$ modulo $w^{js}$, so for $r \ge js$ it is also a closed subscheme of the truncated flop.

\emph{The Hilbert scheme of $Y$.} Now suppose that $r \ge j(n_j' + 1)$, and let $Y' = Y \times_X \Spec R/\mathfrak m_R^{\,j(n_j' + 1)}$, a closed subscheme of $Y$ and of the truncated flop, so that its Hilbert scheme is closed in theirs. The complete local ring $O'$ of its Hilbert scheme at $[C_j]$ is $O_r$, as the universal deformation over $O_r$, which has length $n_j'$, lies in $Y'$. It is also the complete local ring $O$ at $[C_j]$ of the Hilbert scheme of $Y$. Indeed $O$ surjects onto $O'$, and if the surjection were not an isomorphism, $O$ would have a quotient of length $n_j' + 1$ mapping onto $O'$. The corresponding deformation of $C_j$ would lie in $Y'$, so that quotient would be a quotient of $O'$, which has length $n_j'$. Hence $[C_j]$ is an isolated point of the Hilbert scheme of $Y$, and $n_j = n_j'$.

\emph{The contraction algebra.} Write $\cC_{\nu_r} = \cC_l \otimes_{H_l} k[[w]]$ via $\nu_r$. As in the proof of Corollary~\ref{cor:rootsum}, $\cC_\nu$ and $\cC_{\nu_r}$ are cokernels of square matrices over $k[[w]]$ with determinants $c\prod_j f_j^{j^2}$ evaluated along $\nu$ and $\nu_r$, so by the orders above both have length $L = \sum_j j^2 \ord_w f_j(\nu)$ (Lemma~\ref{lem:dvr}(2)). Hence for $r \ge L$ both are killed by $w^r$, and as $\cC_l \otimes_{H_l} k[w]/(w^r)$ depends only on the arc modulo $w^r$, $\cC_{\nu_r} \cong \cC_\nu \cong \Acon$. By Lemma~\ref{lem:basechange}, $\cC_{\nu_r}$ is the contraction algebra of the truncated flop.
\end{proof}

Next we move an arc into general position. Proposition~\ref{prop:GV} uses part~(2) of the lemma, and Theorem~B(2) uses part~(1).

\begin{lem}\label{lem:genericarc}
Let $\nu$ be the germ at $0$ of a curve $\A^1 \to \Spec H_l$ with $\nu(0) = 0$ and $m_j := \ord_w f_j(\nu) < \infty$ for each $j$. For $\xi \in \Spec H_l$ put $\hat\nu(w, s) = \nu(w) + s\xi$, using the addition of the affine space $\Spec H_l$ in a fixed choice of homogeneous generators of $H_l$ (the lemma holds for every such choice), and let $\bar K$ be an algebraic closure of $k((s))$. Then for general $\xi$ the following hold.
\begin{thmlist}
\item Over $\bar K$ the arc $\hat\nu$ meets each $D_j$ near the origin in exactly $m_j$ points, all outside $\pi(\Sigma)$, namely at the closed points of $\Spec\big(k[[w,s]]/(f_j(\hat\nu)) \otimes_{k[[s]]} \bar K\big)$, and $f_j$ vanishes to order one along the arc at each of them. At each such point $p$, $f_j$ is the only $f_i$ vanishing, and $\cC_l \otimes_{H_l} k(p) \cong M_j(\bar K)$.
\item For all but finitely many $s \in k$ the curve $\nu + s\xi$ misses $\pi(\Sigma)$ and meets $D_1 \cup \dots \cup D_l$ transversally (the pullback of each $D_j$ along it is reduced), and for $0 < |s|$ small it meets each $D_j$ near the origin in exactly $m_j$ points.
\end{thmlist}
\end{lem}

\begin{proof}
We first move the arc into general position, then count its points on each $D_j$, with $D = D_1 \cup \dots \cup D_l$ as in \S\ref{ss:universal}. By Proposition~\ref{prop:support}, at every point $p \in D_j \setminus \pi(\Sigma)$, $D$ is smooth with local equation $f_j$, and $\cC_l \otimes_{H_l} k(p) \cong M_j(k)$; by the Conventions this also holds over $\bar K$. We apply Kleiman's transversality theorem \cite[Thm~2]{Kl} to the morphism $\nu\colon \A^1 \to \Spec H_l$ and the transitive action of $\Spec H_l$ on itself by translation. Part~(i), for $\pi(\Sigma)$, which has codimension $2$ in $\Spec H_l$ (\S\ref{ss:lie}), and part~(ii), for the regular subvarieties $D_j \setminus \pi(\Sigma)$, show that the set of $\xi$ for which the curve $\nu + \xi$ meets $\pi(\Sigma)$ or fails to meet $D$ transversally is constructible and not dense. Let $\mathcal B$ be its closure, a proper closed subset. Both conditions are geometric (they are preserved by extending the field). For general $\xi$ the line $k\xi$ is not contained in $\mathcal B$, so its generic point, the point $s\xi$ with $s$ transcendental over $k$, lies outside $\mathcal B$. As $\mathcal B$, $\pi(\Sigma)$ and $D$ are defined over $k$, the curve $\nu + s\xi$ over $\bar K$ also misses $\pi(\Sigma)$ and meets $D$ transversally.

It remains to count the points. By Weierstrass preparation, $k[[w, s]]/(f_j(\hat\nu))$ is finite and free over $k[[s]]$, with special fibre of length $m_j$. Its points over $\bar K$ are the roots of positive valuation of $f_j(\nu(w) + s\xi)$, as a power series in $w$ over $\bar K$, with the same multiplicities, since the unit factor of the Weierstrass decomposition does not vanish at them. By transversality these roots are simple (the pullback of $D_j$ along the moved arc is reduced). Hence over $\bar K$ the arc meets $D_j$ in exactly $m_j$ points near the origin, at each of which $f_j$ vanishes to order one along the arc. For (2), $\mathcal B$ meets the line $k\xi$ in finitely many points, which gives its first part. Weierstrass preparation of $f_j(\nu(w) + s\xi)$ in $\C\{w, s\}$ shows that for $0 < |s|$ small it has exactly $m_j$ zeros near $w = 0$, counted with multiplicity, which are simple by transversality.
\end{proof}

\begin{prop}\label{prop:GV}
Let $\nu$ be a classifying arc of $f$, as in Lemma~\ref{lem:arcs}. Then $n_j = \ord_w f_j(\nu)$.
\end{prop}

\begin{proof}
The geometric input is Proposition~\ref{prop:BKL}; compare \cite[\S3.1, (3.A)]{NW}, where GV invariants of flops of several curves are defined by this count. Steps~1 and~2 treat polynomial arcs, and Step~3 reduces any arc to them by truncation. Step~1 is the count of \cite{BKL}, and Step~2, where we identify the equation of $D_j$, is new.

\emph{Step 1 (GV invariants as intersection numbers).} In Steps~1 and~2 we suppose that $\nu$ is a polynomial map, so that they apply to every flop classified by a polynomial arc, such as the truncated flop of Lemma~\ref{lem:truncation}. Then $f$ is the completion of an algebraic, hence analytic, one-parameter deformation of $Y_0^J$ over a disc. Let $\nu' = \nu + s\xi$ for general $\xi$ and small $s \neq 0$, which is transverse to $D^{\mathrm{curv}}$ near the origin by Lemma~\ref{lem:genericarc}(2). By Proposition~\ref{prop:BKL}(2), over each point near the origin where $\nu'$ meets $D_j$ there is exactly one complete curve, of class $j[C]$ (matching Theorem~A(2), by which the fibre of $\cC_l$ at a general point of $D_j$ is $M_j(k)$), and there are $k_j$ such points; here $k_j = n_j$ (\S\ref{ss:flops}, and the proof of Lemma~\ref{lem:truncation} for the section $w = 0$). By conservation of number for the intersection of the moving arc with $D_j$ near the origin, where it is proper (the Weierstrass preparation argument in the proof of Lemma~\ref{lem:genericarc}), the number of these points is the local intersection number $(\nu \cdot D_j)_0$. Hence
\[
n_j = (\nu \cdot D_j)_0 .
\]

\emph{Step 2 (the equation of $D_j$).} By Proposition~\ref{prop:BKL}(1), $D_j = \pi\big(\bigcup_{\alpha_\bul = j} \alpha^\perp\big) = \{f_j = 0\}$ as sets, the $D_j$ of \S\ref{ss:lie}. By Proposition~\ref{prop:support}, $f_j$ is a reduced equation for $D_j$. Here $(\nu \cdot D_j)_0$ is the intersection number of the parametrised arc, the length of $k[[w]]/(\phi_j(\nu))$ for a reduced local equation $\phi_j$ of $D_j$, which is what Step~1 counts. Hence
\[
n_j = (\nu \cdot D_j)_0 = \ord_w f_j(\nu)
\]
for every polynomial arc $\nu$ that classifies a flop.

\emph{Step 3 (any arc).} For the classifying arc $\nu$ of Lemma~\ref{lem:arcs}, put $m_j = \ord_w f_j(\nu)$, which is finite by that lemma. By Lemma~\ref{lem:truncation}, for $r$ large the truncation $\nu_r$ classifies a flop with the same orders $m_j$, so by Steps~1 and~2 its GV invariants are $n_j' = m_j$. Taking also $r \ge \max_j j(m_j + 1)$, Lemma~\ref{lem:truncation} gives $n_j = n_j' = m_j$.
\end{proof}

\begin{cor}[{Toda's formula; recovering \cite[Thm~1.1]{T}; see also \cite[Cor~4.5]{HT}}]\label{cor:toda}
$\dim \Acon = \sum_j j^2 n_j$.
\end{cor}

\begin{proof}
By Corollary~\ref{cor:rootsum} and Proposition~\ref{prop:GV}, $\dim \Acon = \sum_j j^2 \ord_w f_j(\nu) = \sum_j j^2 n_j$.
\end{proof}

\subsection{Proof of Theorem~B(1) and~(2)}\label{ss:proofB}

The deformation in Theorem~B(2) comes from moving the arc into general position (Lemma~\ref{lem:genericarc}).

\begin{thm}[$=$ Theorem~B]\label{thm:B}
Let $f\colon Y \to X = \Spec R$ be a flop of length $l$, with GV invariants $n_1, \dots, n_l$.
\begin{thmlist}
\item $\dim Z(\Acon) = \dim \HH_0(\Acon) = n_1 + \dots + n_l$.
\item $\Acon$ deforms to $\SGV$ (\S\ref{ss:moduli}), through a family obtained by deforming a truncation of its classifying arc (Lemmas~\ref{lem:truncation} and~\ref{lem:genericarc}).
\item Let $O$ be a complete local noetherian $k$-algebra with residue field $k$, let $f_O\colon \mathcal Y \to \mathcal X = \Spec R_O$ be a flat deformation of $f$ over $O$: a projective morphism of schemes flat over $O$, with $R_O$ complete local noetherian, whose fibre over the closed point of $\Spec O$ is $f$. Let $\Acon(\mathcal X) = A_O/A_Oe_0A_O$ be its relative contraction algebra, where $A_O = \End_{\mathcal Y}(\cO \oplus \mathcal M_O)$ for the relative tilting bundle of Van den Bergh. Then $\Acon(\mathcal X)$, $\HH_0(\Acon(\mathcal X))$ and the centre $Z(\Acon(\mathcal X))$ are free $O$-modules, of ranks $\sum_j j^2 n_j$, $\sum_j n_j$ and $\sum_j n_j$, with fibres $\Acon$, $\HH_0(\Acon)$ and $Z(\Acon)$ at the closed point.
\end{thmlist}
\end{thm}

Part (3) is Hua and Toda's conjecture \cite[Conj~4.3]{HT}, and we prove it in \S\ref{ss:defflop}.

\begin{proof}[Proof of (1) and (2)]
(1) Corollary~\ref{cor:rootsum} and Proposition~\ref{prop:GV} give $\dim \HH_0(\Acon) = \sum_j \ord_w f_j(\nu) = \sum n_j$, and $Z(\Acon) \cong \HH_0(\Acon)^*$ since $\Acon$ is symmetric (\S\ref{ss:moduli}).

(2) By Lemma~\ref{lem:truncation}, for $r$ large the truncation of $\nu$ modulo $w^r$, a polynomial map, classifies a flop with the same contraction algebra $\Acon$ and the same orders of vanishing of the $f_j$. Hence we may assume that $\nu$ is the germ at $0$ of a curve $\mathbb{A}^1 \to \Spec H_l$. For $\xi \in \Spec H_l$ put $\hat\nu(w, s) = \nu(w) + s\xi$ as in Lemma~\ref{lem:genericarc}, and set $\cC_{\hat\nu} = \cC_l \otimes_{H_l} k[[w, s]]$, via $\hat\nu$.

\emph{Flatness.} Pulling back the presentation of $\cC_l$ of Theorem~A(3) gives a square matrix over the two-dimensional regular ring $k[[w, s]]$ with nonzero determinant and cokernel $\cC_{\hat\nu}$, so $\cC_{\hat\nu}$ is Cohen--Macaulay of dimension $1$ (Lemma~\ref{lem:freeness}(3)). Modulo $s$ the determinant is that of the presentation pulled back along $\nu$, which is nonzero, so by Weierstrass preparation it is a unit times a monic polynomial $q \in k[[s]][w]$. By the adjugate identity $q$ annihilates $\cC_{\hat\nu}$, so $\cC_{\hat\nu}$ is finite over $k[[w, s]]/(q)$, which is free of finite rank over $k[[s]]$. Hence $\cC_{\hat\nu}$ is finite over $k[[s]]$, so free over it (Lemma~\ref{lem:freeness}(2)): a flat family of algebras over $k[[s]]$ with special fibre $\Acon$.

\emph{The generic fibre.} By Lemma~\ref{lem:genericarc}(1), with $m_j = n_j$ by Proposition~\ref{prop:GV}, for general $\xi$ the arc $\hat\nu$ meets each $D_j$ over $\bar K$ in exactly $n_j$ points near the origin, at each of which $f_j$ vanishes to order one. These points are defined over a finite extension $K' = k((s^{1/N}))$ of $k((s))$, and we argue over $K'$ before extending scalars to $\bar K$. The ring $k[[w, s]] \otimes_{k[[s]]} K' = k[[w, s^{1/N}]][1/s]$ is regular of dimension one, so its local ring at each such point $p$ is a discrete valuation ring. At $p$ only the factor $f_j^{j^2}$ of the determinant of Theorem~A(3) vanishes, so by Lemma~\ref{lem:dvr}(2) the factor of $\cC_{\hat\nu} \otimes_{k[[s]]} K'$ at $p$, a module over that discrete valuation ring, has length $j^2$. Its fibre $\cC_l \otimes_{H_l} k(p)$ also has dimension $j^2$ over $k(p) = K'$, as it becomes $M_j(\bar K)$ over $\bar K$ (Lemma~\ref{lem:genericarc}(1)), so a uniformiser at $p$ acts as zero (Lemma~\ref{lem:dvr}(1)), and this factor is $\cC_l \otimes_{H_l} k(p)$. Hence $\cC_{\hat\nu} \otimes \bar K \cong \prod_j M_j(\bar K)^{n_j} = \SGV \otimes \bar K$, so $\Acon$ deforms to $\SGV$ (\S\ref{ss:moduli}).
\end{proof}

\subsection{Proof of Theorem~B(3)}\label{ss:defflop}

The same method works over a complete local base and proves Theorem~B(3). Remark~\ref{rem:HT} then deduces Hua and Toda's conjecture in their own setting, a deformation over an open neighbourhood of $0$ in $\A^1$.

\begin{proof}[Proof of Theorem~B(3)]
The proof has three steps. We classify the deformation by a map to $\Spec H_l$ and apply Lemma~\ref{lem:basechange} over its base. We then show that the pulled-back presentations of Theorem~A(3) stay injective, so that $\Acon(\mathcal X)$ and its $\HH_0$ are free. Lemma~\ref{lem:symfamily} then gives the centre.

\emph{Step 1 (base change).} We show that $\mathcal Y \to \mathcal X$ is pulled back from the universal flop along a map $\Spec B \to \Spec H_l$, with $B$ as below. First, $Rf_{O*}\cO_{\mathcal Y} = \cO_{\mathcal X}$ by Lemma~\ref{lem:flatfamilies}(2) over $O$: $\mathcal Y$ and $\mathcal X$ are flat over $O$, the closed point of $\mathcal X = \Spec R_O$ lies over the closed point of $\Spec O$, as $R_O/\mathfrak m_OR_O = R$ is nonzero, and the fibre of $f_O$ over it is $f$, with $Rf_*\cO_Y = \cO_X$. In particular $\mathcal X$ is the affinisation of $\mathcal Y$, so identifying $\mathcal Y$ with a pullback of $\cY_l$ below also identifies $\mathcal X$ with the pullback of $\cX_l$. Choose $w \in R_O$ restricting on $X$ to the equation of a generic hyperplane section through the singular point, so that the corresponding section of $Y$ is the partial resolution $Y_0$ of the proof of Lemma~\ref{lem:arcs}. Put $B = k[[w]] \mathbin{\hat\otimes} O$.

\emph{Flatness over $B$.} We show that $\mathcal Y$ and $\mathcal X$ are flat over $B$, using only that they are flat over $O$ and that $w$ is a nonzerodivisor on their fibres $Y$ and $X$ over the closed point of $O$. Let $y$ be a point of $Y_0$, the fibre of $\mathcal Y$ over the closed point of $B$, and put $M = \cO_{\mathcal Y, y}$, a module over the local ring $\cO_{\mathcal Y, y}$ to which $O \to B \to \cO_{\mathcal Y, y}$ are local homomorphisms. As $M$ is flat over $O$, each $\mathfrak m_O^i M/\mathfrak m_O^{i+1} M \cong (\mathfrak m_O^i/\mathfrak m_O^{i+1}) \otimes_k M/\mathfrak m_O M$ is a direct sum of copies of $M/\mathfrak m_O M = \cO_{Y, y}$, on which $w$ is injective; so $w$ is injective on $M$, which is $\mathfrak m_O$-adically separated. Tensoring $0 \to M \xrightarrow{w} M \to M/wM \to 0$ with $k$ over $O$ gives $\Tor_1^O(M/wM, k) = \ker(w \text{ on } M/\mathfrak m_O M) = 0$, so $M/wM$ is flat over $O$ by the local criterion of flatness. The resolution $0 \to B \xrightarrow{w} B \to O \to 0$ of $O = B/wB$ gives $\Tor_q^B(M, O) = 0$ for $q \ge 1$, so the change-of-rings spectral sequence gives $\Tor_1^B(M, k) \cong \Tor_1^O(M/wM, k) = 0$, and $M$ is flat over $B$, again by the local criterion. Hence $\mathcal Y$ is flat over $B$: flatness at a point gives flatness at its generisations, and every point of $\mathcal Y$ specialises to a closed point, which lies in $Y_0$ as $\mathcal Y$ is proper over the local $\mathcal X$. The same argument, with $R_O$ in place of $M$, shows that $\mathcal X$ is flat over $B$, as $w$ is a nonzerodivisor on the domain $R$.

\emph{The classifying map.} Thus $\mathcal Y \to \mathcal X$ is a deformation of $Y_0$ over $B$ as in Lemma~\ref{lem:classifying}, so there is a map $\nu_O\colon \Spec B \to \Spec H_l$ along which $\mathcal Y \to \mathcal X$ is, after completion at the singular point, the pullback of the universal flop. By Lemma~\ref{lem:basechange}, the pullback of the universal tilting bundle is $\cO \oplus N_1 \cong \cO \oplus \mathcal M_O$, a projective generator of ${}^0\Per(\mathcal Y/\mathcal X)$, and $\Acon(\mathcal X) \cong \cC_{\nu_O} = \cC_l \otimes_{H_l} B$. As $\HH_0$ commutes with base change, $\HH_0(\Acon(\mathcal X)) \cong \HH_0(\cC_l) \otimes_{H_l} B$.

\emph{Step 2 (freeness).} By Theorem~A(3), $\cC_l$ is the cokernel of a square matrix over $H_l$; let $\varphi$ be its pullback to $B$, so that $\cC_{\nu_O} = \operatorname{coker}\varphi$. Modulo the maximal ideal $\mathfrak m_O$ of $O$, $\varphi$ is the pullback of that matrix along $\nu = \nu_O \bmod \mathfrak m_O$, a classifying arc of $f$, whose determinant $c\prod_j f_j(\nu)^{j^2}$ is nonzero (proof of Corollary~\ref{cor:rootsum}). Hence $\det\varphi$ has a coefficient that is a unit of $O$, and by Weierstrass preparation it is a unit times a monic polynomial $q$ in $w$, a nonzerodivisor, so $\varphi$ is injective. By the adjugate identity $q$ annihilates $\cC_{\nu_O}$, and $B/(q)$ is free of finite rank over $O$, so $\cC_{\nu_O}$ is finite over $O$. Hence $0 \to F \to F \to \cC_{\nu_O} \to 0$, with $F$ free over $B$, is a resolution by flat $O$-modules, and $\Tor_1^O(\cC_{\nu_O}, k) = \ker(\varphi \otimes_O k) = 0$. As $\cC_{\nu_O}$ is finite over $O$, it is free. The same argument applies to $\HH_0(\cC_l) \otimes_{H_l} B = \HH_0(\Acon(\mathcal X))$, with determinant $c'(f_1 \cdots f_l)(\nu)$.

\emph{Step 3 (the centre).} $\Acon(\mathcal X)$ and its $\HH_0$ are free $O$-modules of finite rank by what has just been shown, and the fibre at the closed point is the symmetric algebra $\Acon$ (Proposition~\ref{prop:Acon}(2)). Lemma~\ref{lem:symfamily} then shows that $Z(\Acon(\mathcal X))$ is free, of the same rank as $\HH_0$, with fibre $Z(\Acon)$. The ranks are $\dim \Acon = \sum_j j^2 n_j$ (Corollary~\ref{cor:toda}) and $\dim \HH_0(\Acon) = \sum_j n_j$ (Theorem~B(1)).
\end{proof}

\begin{rem}[Hua and Toda's formulation]\label{rem:HT}
Hua and Toda's conjecture concerns a deformation over an open subset $U \ni 0$ of $\A^1$, near the singular point $x$ of the central fibre, with contraction algebras defined through their own tilting bundle $\mathcal U$. Theorem~B(3) covers it. First, flatness at $x$ of a coherent sheaf can be checked after completing at $x$, and the relative tilting bundle, its endomorphism algebra, $\HH_0$ and the centre commute with this flat base change. Hence Theorem~B(3) over $O = \hat\cO_{U,0}$ gives flatness at $x$, and so near $x$, as the flat locus is open. Second, their algebra is a matrix algebra over ours. The argument of \cite[Prop~2.3]{HT}, through Orlov's description of the singularity category \cite[Prop~1.21]{Orl}, identifies their relative contraction algebra with $A'/A'e_0A'$, where $A' = \End(\mathcal U)$. That argument applies, as their total space is flat over a smooth curve with Gorenstein fibres and the pushforward $N'$ of $\mathcal U$ is flat over it with maximal Cohen--Macaulay fibres. Their conjecture is stated for $\mathbf{R}\Hom(N', N')$ in the singularity category, which near $x$ is $2$-periodic with $\Ext^1(N', N') = 0$ \cite[Prop~3.2, Rem~3.3]{HT}, so only this algebra matters. After completion at $x$ the line bundle from which Hua and Toda build $\mathcal U$ is ample, so $\mathcal U$ is a projective generator of ${}^0\Per$ \cite[Lemma~3.4.5 and the proof of Prop~3.2.5]{VdB}, hence $\cO^{a} \oplus \mathcal M^{b}$ with $a, b \ge 1$ (Theorem~\ref{thm:VdB}(3)), and $A'/A'e_0A' \cong M_b(\Acon(\mathcal X))$. Flatness of the algebra, of its $\HH_0$ and of its centre is invariant under this Morita equivalence. Over the closed point their algebra is the basic contraction algebra \cite[\S4]{HT}, of dimension $\sum_j j^2 n_j$ \cite[Thm~4.1]{HT}, so $b = 1$ and Hua and Toda's contraction algebra is isomorphic to $\Acon$.
\end{rem}

\subsection{Proof of Theorem~C}\label{ss:proofC}

We now combine Theorem~B with the reduction of Section~\ref{sec:reduction}.

\begin{thm}[$=$ Theorem~C]\label{thm:C}
For every flop of length $l$, with $d = \dim \Acon$, $[\Acon]$ is a smooth point of $\Ass_d$, and the closure of $G\cdot[\SGV]$ is the only irreducible component of $\Ass_d$ through it. Hence $\SGV$ is the only semisimple algebra to which $\Acon$ deforms, and $\dim Z$ is constant on a neighbourhood of $[\Acon]$.
\end{thm}

\begin{proof}
By Theorem~B(2), $\Acon$ deforms to $\SGV$; by Theorem~B(1), $\dim Z(\Acon) = \sum n_j$, which is $\dim Z(\SGV)$. Apply Proposition~\ref{prop:criterion}.
\end{proof}

\section{\texorpdfstring{The centre map (Proposition~D)}{The centre map (Proposition D)}}\label{sec:centremap}

In this section we show that the second half of \cite[Conj~1.4(2)]{JS}, smoothness of the centre map, fails for the Laufer flop when its target is a moduli space of algebras (as a map to a Grassmannian it is smooth for every flop, \S\ref{ss:centremap}). The first half, constancy of $\dim Z$ near $[\Acon]$, holds for every contraction algebra by Theorem~C. The argument rests on the fact that a smooth morphism maps a smooth point to a smooth point of the target. In Lemmas~\ref{lem:laufer}--\ref{lem:alg6} we show that $[L_1]$ is a smooth point of $\Ass_9$, as Theorem~C predicts, but that its centre, a $6$-dimensional local algebra needing three generators, is a singular point of the scheme of $6$-dimensional commutative algebra structures, which here is smoothly equivalent to $\Hilb^6(\A^3)$. Hence the map $A \mapsto Z(A)$ cannot be smooth at $[L_1]$ as a map to a moduli space of algebras, and the proof of Proposition~D makes this precise for three such targets.

Throughout this section $L_1 = k\langle x,y\rangle/(xy + yx,\ x^3 + y^2)$ is the contraction algebra of the Laufer flop $\mathsf{D}_{2,\infty}$ \cite[Example~1.3]{DW}, \cite[\S\S4.2, 4.4]{AM}, \cite[Example~5.2]{K}. It is graded with $\deg x = 2$ and $\deg y = 3$. Proposition~D does not use Sections~\ref{sec:cartan}--\ref{sec:flops}: we deduce smoothness of $[L_1]$ from Section~\ref{sec:reduction} and the known deformation of $L_1$ to $k^5 \times M_2(k)$, and use the irreducibility of $\Hilb^6(\A^3)$ for its centre.

\subsection{The algebra and its centre}

We first compute a basis of $L_1$ and its centre.

\begin{lem}\label{lem:laufer}
$L_1$ has basis $1, x, x^2, x^3, x^4, x^5, y, xy, x^2y$, in degrees $0, 2, 4, 6, 8, 10, 3, 5, 7$, so $\dim L_1 = 9$. Its centre has basis $1, x^2, x^3, x^4, x^5, x^2y$, and
\[
Z(L_1) \;\cong\; k[z_1,z_2,z_3]/(z_1z_3,\ z_2z_3,\ z_3^2,\ z_2^2,\ z_1^3,\ z_1^2z_2), \qquad z_1 = x^2,\quad z_2 = x^3,\quad z_3 = x^2y .
\]
This agrees with \cite[Prop~3.3]{JS}.
\end{lem}

\begin{proof}
The relations give $yx = -xy$ and $y^2 = -x^3$, so $x^2$ is central and every word reduces to $x^a$ or $x^a y$. Computing $y^3$ in two ways gives $y\cdot y^2 = -yx^3 = x^3y$ and $y^2\cdot y = -x^3y$, so $x^3y = 0$, and then $x^6 = y^4 = y\cdot y^3 = 0$. Hence the nine listed monomials span. They are independent by the diamond lemma for the rewriting system $yx \to -xy$, $y^2 \to -x^3$, $x^3y \to 0$, $x^6 \to 0$. An element of $L_1$ is central if and only if it commutes with $x$ and $y$. Since $L_1$ is graded and each degree is at most one-dimensional, this can be checked monomial by monomial: $x$, $y$ and $xy$ are not central (for instance $x\cdot xy = x^2y$ but $xy\cdot x = -x^2y \neq 0$), and the other six are. The monomials $1, z_1, z_2, z_3, z_1^2, z_1z_2$ have distinct degrees $0, 4, 6, 7, 8, 10$ and are nonzero, while $z_1z_3 = x^4y$, $z_2z_3 = x^5y$, $z_3^2 = -x^7$, $z_2^2 = x^6$, $z_1^3 = x^6$ and $z_1^2z_2 = x^7$ all vanish. This gives the presentation.
\end{proof}

\subsection{\texorpdfstring{$[L_1]$ is a smooth point of $\Ass_9$}{[L1] is a smooth point of Ass9}}

We next compute the derivations of $L_1$ degree by degree, and apply the argument of Proposition~\ref{prop:criterion}.

\begin{lem}\label{lem:lauferSmooth}
$h^1(L_1) = h^2(L_1) = 5$, and $[L_1]$ is a smooth point of $\Ass_9$, near which $\dim Z$ is constant. The closure of the orbit of $\SGV = k^5 \times M_2(k)$ is the only irreducible component of $\Ass_9$ through $[L_1]$.
\end{lem}

\begin{proof}
As $L_1$ is graded, so is $\Der(L_1)$, and each degree of $L_1$ is at most one-dimensional (Lemma~\ref{lem:laufer}). A derivation of degree $d$ sends $x$ and $y$ to multiples of the basis monomials of degrees $d + 2$ and $d + 3$, and the two relations give at most two linear conditions on the two coefficients. Checking $-3 \le d \le 8$ (outside this range both targets vanish) gives, up to scalars, exactly one derivation in each of the degrees $0, 2, 3, 4, 5, 6, 7, 8$ and none in the others: the Euler derivation $E$ ($x \mapsto 2x$, $y \mapsto 3y$), $\operatorname{ad} x$, $\operatorname{ad} y$, $x^2E$, $\operatorname{ad}(xy)$, $x^3E$, $x^2yE$ and $x^4E$. For instance, in degree $1$ the derivation $x \mapsto ay$, $y \mapsto bx^2$ must satisfy $2(b - a)x^3 = 0$ and $(a + 2b)x^2y = 0$, so it vanishes. Hence $\dim \Der(L_1) = 8$. Inner derivations have dimension $\dim L_1 - \dim Z(L_1) = 3$, so $h^1(L_1) = 5$, and $\HH^1(L_1) = Z(L_1)\cdot E$ (here $x^5E = 0$, while $x^2yE \neq 0$). Also $h^2(L_1) = 5$ by Lemma~\ref{lem:pan}. The algebra $L_1$ deforms to $S = k^5 \times M_2(k)$ \cite[Thm~4.2]{BPS}, \cite[Thm~3.13, Rem~3.4]{JS}, and $\dim Z(S) = 6 = \dim Z(L_1)$. Hence Lemma~\ref{lem:tangent} gives $\dim T_{[L_1]}\Ass_9 = \dim (G\cdot[S])$, and the proof of Proposition~\ref{prop:criterion}, with $S$ in place of $\SGV$, shows that $[L_1]$ is a smooth point, that the closure of $G\cdot[S]$ is the only irreducible component of $\Ass_9$ through it, and that $\dim Z$ is constant near it (see also Remark~\ref{rem:AD}). The smoothness does not use that $S = \SGV$. The equality holds because the abelianisation $k[x,y]/(xy, x^3 + y^2)$ has basis $1, x, x^2, x^3, y$, so $n_1 = 5$ (\S\ref{ss:flops}), and Toda's formula, $9 = n_1 + 4n_2$, gives $n_2 = 1$. Toda's formula applies because this flop is the completion of the algebraic flop of \cite[Example~1.3]{T} (its first case), where Toda finds these values.
\end{proof}

\subsection{\texorpdfstring{$Z(L_1)$ is a singular point of $\Hilb^6(\A^3)$}{Z(L1) is a singular point of Hilb6(A3)}}

We now show that the centre of $L_1$, presented as a quotient of a polynomial ring in three variables, is a singular point of $\Hilb^6(\A^3)$, by computing the tangent space there. Let $\Lambda = k[z_1,z_2,z_3]$ and $I = (z_1z_3,\allowbreak z_2z_3,\allowbreak z_3^2,\allowbreak z_2^2,\allowbreak z_1^3,\allowbreak z_1^2z_2)$. Then $\Lambda/I \cong Z(L_1)$ has dimension $6$ (Lemma~\ref{lem:laufer}), so $[I]$ is a point of $\Hilb^6(\A^3)$.

\begin{lem}\label{lem:hilbTangent}
$[I]$ is a singular point of $\Hilb^6(\A^3)$: the tangent space $T_{[I]} \Hilb^6(\A^3)$ has dimension at least $19$, while $\Hilb^6(\A^3)$ has dimension $18$.
\end{lem}

\begin{proof}
Since $\Hilb^d(\A^3)$ is irreducible for $d \le 7$ \cite[Cor~4]{Maz}, \cite[\S1]{DJNT}, $\Hilb^6(\A^3)$ is the closure of the locus of six distinct points and has dimension $18$, and the tangent space at $[I]$ is $\Hom_\Lambda(I, \Lambda/I)$. It contains the three translations, sending $g \in I$ to $\partial g/\partial z_i$ modulo $I$, which are $\Lambda$-linear as $\partial/\partial z_i$ is a derivation. It also contains the following sixteen homomorphisms, each sending the generators of $I$ shown to the standard monomials shown, and the other generators to $0$:
\[
\begin{gathered}
z_1^3 \mapsto z_3,\quad z_1^3 \mapsto z_2,\quad z_1^3 \mapsto z_1z_2,\quad z_1^2z_2 \mapsto z_3,\quad (z_1^3, z_1^2z_2) \mapsto (z_1, z_2), \\
z_1z_3 \mapsto z_2,\quad z_1z_3 \mapsto z_1z_2,\quad z_1z_3 \mapsto z_1^2,\quad z_2z_3 \mapsto z_1z_2,\quad z_2z_3 \mapsto z_1^2, \\
z_2^2 \mapsto z_3,\quad z_2^2 \mapsto z_1,\quad z_2^2 \mapsto z_1z_2,\quad z_2^2 \mapsto z_1^2,\quad z_3^2 \mapsto z_1z_2,\quad z_3^2 \mapsto z_1^2.
\end{gathered}
\]
Each is well defined, as a direct check shows that it respects the fifteen relations $z^{\mathrm{lcm}(a, b) - a} z^a = z^{\mathrm{lcm}(a, b) - b} z^b$ between pairs of generators $z^a$, $z^b$ of $I$, which generate all relations among the generators of a monomial ideal. The torus $(k^\times)^3$ acts on $\Hom_\Lambda(I, \Lambda/I)$, as $I$ is a monomial ideal, and these nineteen homomorphisms are weight vectors. Their weights are distinct, except for $z_1z_3 \mapsto z_1^2$ and $z_2z_3 \mapsto z_1z_2$, which are nonzero on different generators, so the nineteen are linearly independent. Hence $\dim \Hom_\Lambda(I, \Lambda/I) \ge 19 > 18$. A computer calculation, which the lemma does not use, gives the full dimension, $24$.
\end{proof}

\subsection{The centre map}\label{ss:centremap}

We now define the centre map and prove Proposition~D for three choices of its target. Part (2) of \cite[Conj~1.4]{JS} reads: ``the dimension of the centre is constant in an open neighbourhood of $[A]$ and the associated map $A \mapsto Z(A)$ is smooth'' (here $A = L_1$). It does not name the target of the map. By Lemma~\ref{lem:lauferSmooth}, $[L_1]$ has a smooth open neighbourhood $U \subset \Ass_9$ on which $\dim Z = 6$, and we may take $U$ to be $G$-stable, replacing it by the union of its translates. Hence the centres $Z(V, \mu)$, for $\mu \in U$, form a rank-$6$ subbundle of $V \times U$ (here $\dim V = 9$), which defines a map $\varphi$ from $U$ to the Grassmannian $\mathrm{Gr}$ of $6$-dimensional subspaces of $V$ containing $1$. This map is smooth. Indeed, $\varphi$ is $G$-equivariant, and $G$ acts transitively on $\mathrm{Gr}$. Then for $u \in U$ the orbit map $G \to \mathrm{Gr}$, $g \mapsto g\cdot\varphi(u) = \varphi(g\cdot u)$, is smooth and factors through $\varphi$. Hence the differential of $\varphi$ is onto at every point of $U$, and $\varphi$ is smooth, as $U$ and $\mathrm{Gr}$ are. The same argument, with Theorem~C (a smooth neighbourhood on which $\dim Z$ is constant) and Theorem~B(1) ($\dim Z(\Acon) = n_1 + \dots + n_l$) in place of Lemma~\ref{lem:lauferSmooth}, shows that for every contraction algebra the map to the Grassmannian is smooth near $[\Acon]$.

Jelisiejew and Smoktunowicz motivate the conjecture by the deformations of the commutative algebra $Z(A)$ \cite[\S1]{JS}, so we record $Z(A)$ as a point of a moduli space of algebras. Let $\Alg_6 \subset \Ass_6$ be the closed subscheme of commutative structures (\S\ref{ss:moduli}, with $\dim V = 6$), and $\sAlg_6$ the stack of $6$-dimensional commutative algebras, the quotient of $\Alg_6$ by $G$, or equivalently the quotient by $GL_6$ of the commutative algebra structures on $k^6$ with any unit \cite[Prop~3.1, Def~3.2]{Poo}. The subbundle of centres is a family of $6$-dimensional commutative algebras over $U$, and so defines the centre map $\mathcal Z\colon U \to \sAlg_6$. Near $[L_1]$ a local trivialisation of the subbundle that preserves $1$, a \emph{frame}, lifts $\mathcal Z$ to a map to $\Alg_6 \subset \Ass_6$, the \emph{framed} map, which we also denote by $\mathcal Z$. We take the target to be $\sAlg_6$, and the proof of Proposition~D also covers $\Alg_6$ and the scheme $\Ass_6$ of all $6$-dimensional unital associative algebras. The scheme $\Ass_6$ is their moduli space $M_6$, and part~(1) of the conjecture is stated in their spaces $M_d$, with $d = \dim A$.

\begin{lem}\label{lem:alg6}
$Z(L_1)$ is a singular point of $\Alg_6$.
\end{lem}

\begin{proof}
Let $\mathcal G$ be the scheme of pairs $(B, (b_1, b_2, b_3))$ with $B \in \Alg_6$ and $b_1, b_2, b_3$ generating $B$. $\mathcal G$ is an open subscheme of the trivial bundle $\Alg_6 \times V^3$, so the forgetful map $\mathcal G \to \Alg_6$ is smooth. The map $\mathcal G \to \Hilb^6(\A^3)$ sending $(B, b)$ to the kernel of $k[z_1,z_2,z_3] \to B$, $z_1, z_2, z_3 \mapsto b_1, b_2, b_3$, is a $G$-torsor \cite[Def~1.4, \S\S3--4]{Poo}, with the unit fixed: its fibre over $[I]$ is the set of linear isomorphisms $\Lambda/I \to V$ sending $1$ to the unit. The point $(Z(L_1), (z_1, z_2, z_3))$ of $\mathcal G$ lies over $Z(L_1)$ and over $[I]$, and both maps are smooth, so $\Alg_6$ is regular at $Z(L_1)$ if and only if $\Hilb^6(\A^3)$ is regular at $[I]$, and likewise for the reduced structures. Now apply Lemma~\ref{lem:hilbTangent}.
\end{proof}

\begin{proof}[Proof of Proposition~D]
By Lemma~\ref{lem:lauferSmooth}, $[L_1]$ is a smooth point of $\Ass_9$. If $\mathcal Z$ were smooth at $[L_1]$, then $Z(L_1)$ would be a smooth point of $\sAlg_6$, since a smooth morphism is flat and regularity descends along flat local maps. That contradicts Lemma~\ref{lem:alg6}, as the quotient map $\Alg_6 \to \sAlg_6$ is smooth. The same argument applies with $\Alg_6$ in place of $\sAlg_6$ after choosing a frame. For $\Ass_6$, the framed map $\mathcal Z$, defined near $[L_1]$, factors through the closed subscheme $\Alg_6$. Hence the local homomorphism $\cO_{\Ass_6, Z(L_1)} \to \cO_{\Ass_9, [L_1]}$ kills the ideal of $\Alg_6$. If $\mathcal Z$ were smooth at $[L_1]$, this homomorphism would be flat and local, hence injective. Then $\Alg_6$ would equal $\Ass_6$ near $Z(L_1)$, and $\mathcal Z$ would be smooth at $[L_1]$ as a map to $\Alg_6$, which the previous case excludes.
\end{proof}

\begin{rem}[Reduced structures]\label{rem:reduced}
Proposition~D also holds with reduced structures, as $[I]$ is a singular point of $\Hilb^6(\A^3)_{\mathrm{red}}$ too. Indeed, for each of the sixteen homomorphisms $\varphi$ listed in the proof of Lemma~\ref{lem:hilbTangent}, take a monomial order refining a positive grading of $z_1, z_2, z_3$ in which $\varphi(z^a)$ has lower degree than $z^a$ for each generator $z^a$. Then, with a parameter $s$ in degree $0$, the polynomials $z^a + s\varphi(z^a)$ have leading terms $z^a$ and form a Gr\"obner basis, as is straightforward to check with computer algebra, so they define a flat family over $k[s]$, a curve through $[I]$ in $\Hilb^6(\A^3)$ with tangent vector $\varphi$. With the three translations, the reduced tangent space at $[I]$ therefore has dimension at least $16 + 3 = 19 > 18$.
\end{rem}

\section{The root formula and its applications}\label{sec:examples}

In this section we apply the results of Section~\ref{sec:flops} to explicit arcs, with the aim of determining all the GV invariants that can occur. Lifting an arc to $\h$ gives the root formula of \S\ref{ss:rootformula}, which writes each GV invariant, and so each dimension in Theorem~B, as a sum over roots. Lemma~\ref{lem:flopcrit} puts the condition of Lemma~\ref{lem:arcflop} for an arc to classify a flop in terms of the parameters of \cite{K}. In length~$2$ we continue Example~\ref{ex:running} in \S\ref{ss:rootformula2}, and in Lemma~\ref{lem:length2} we recover the GV invariants found by Brown and Wemyss \cite[Cor~8.13]{BW}, realising each pair by an explicit arc. In Proposition~\ref{prop:typeE} we then determine exactly which GV invariants occur in lengths $3$--$6$, and realise each of them by an explicit arc. In \S\ref{ss:typeE} we also revisit the arcs of \cite[Examples~5.3--5.5]{K} and a second arc in length~$3$. Finally, in \S\ref{ss:laufer} we treat the Laufer flops, realising the first by an explicit arc. Throughout, the presentations and parameters of \cite{K} are taken with the corrections of Appendix~\ref{app:corrections}. This section uses Sections~\ref{sec:cartan} and~\ref{sec:squeeze} only through Section~\ref{sec:flops}, and Section~\ref{sec:centremap} for the Laufer flop.

\subsection{The root formula}\label{ss:rootformula}

In this subsection we lift the classifying arc to $\h$, where each $f_j$ is a product of roots, and write each GV invariant as a sum over roots (Corollary~\ref{prop:rootformula}). Remark~\ref{rem:Kpresentation} first makes the presentation of \cite{K} explicit, for Lemma~\ref{lem:flopcrit} and the examples.

\begin{rem}[The explicit presentation]\label{rem:Kpresentation}
The identification $A_l \otimes_{H_l} k[\h] \cong e\cA_\Gamma e$ of \S\ref{ss:universal} is explicit, following \cite[\S\S2.2, 4.5, Examples~4.24--4.27]{K}. In its notation, let $H_\Gamma = k[t_i : i \in \tGamma_0]/(\sum \omega_i t_i)$, where $\omega_i$ are the coefficients of $\delta$ (so $\omega_0 = 1$), identified with $k[\h]$ via $t_i = \lambda_i$. There $e\cA_\Gamma e$ is presented over $H_\Gamma$, and a change of variables puts all coefficients in $H_l = H_\Gamma^{W_J}$; the resulting $H_l$-presentation is its Theorem~1.3 (see Appendix~\ref{app:corrections} for the wording of Prop~4.21 there). In the arm coordinates $\tau_1, \dots, \tau_m$ of \S\ref{ss:universal} (``Generators of $H_l$''), the parameters $T$ of that presentation are, up to sign, the coefficients of the arm polynomial $\prod_k (x + \tau_k)$ of each arm, with $T_i$ the coefficient of $x^i$ (after the corrections in Appendix~\ref{app:corrections}); for instance, in length $3$ the arms $\beta$ and $\gamma$ (Figure~\ref{fig:pairs}) have arm polynomials $x^3 - T_1^\beta x - T_0^\beta$ and $x^3 - T_1^\gamma x - T_0^\gamma$. Besides this section, only Remark~\ref{rem:noER} uses this.
\end{rem}

We now lift the arc.

\begin{lem}\label{lem:lift}
Let $\nu\colon \Spec k[[w]] \to \Spec H_l$ send the closed point to the origin, for instance the classifying arc of a flop (Lemma~\ref{lem:arcs}). There are $N \ge 1$ and a map $\tnu\colon \Spec k[[v]] \to \h$ with $\pi \circ \tnu = \nu \circ (w \mapsto v^N)$.
\end{lem}

\begin{proof}
$k[[w]] \otimes_{H_l} k[\h]$ is finite and free of rank $|W_J|$ over $k[[w]]$, and it is a complete local ring, as $\pi^{-1}(0) = \{0\}$. Choose a minimal prime of it. The quotient is a complete local domain of dimension one, finite over $k[[w]]$, and its normalisation is $k[[v]]$ with $w = c v^N$ for a unit $c$. Since $k$ is algebraically closed of characteristic zero, we may replace $v$ by $c^{1/N}v$. This gives $\tnu$.
\end{proof}

For a linear form $\alpha$ on $\h$, write $\ord_w \alpha(\tnu) := (1/N)\cdot\ord_v \alpha(\tnu(v)) \in (1/N)\Z$.

\begin{cor}[Root formula]\label{prop:rootformula}
For every flop of length $l$, with a classifying arc $\nu$ and a lift $\tnu$ as in Lemma~\ref{lem:lift}, $n_j = \ord_w f_j(\nu) = \sum_{\alpha_\bul = j} \ord_w \alpha(\tnu)$. Hence
\[
\dim \HH_0(\Acon) = \sum_{\alpha \in \Phi^+_\bul} \ord_w \alpha(\tnu), \qquad \dim \Acon = \sum_{\alpha \in \Phi^+_\bul} \alpha_\bul^2 \, \ord_w \alpha(\tnu).
\]
\end{cor}

\begin{proof}
The first equality is Proposition~\ref{prop:GV}. Since $\pi^* f_j = \prod_{\alpha_\bul = j} \alpha$ (\S\ref{ss:lie}) and $\pi \circ \tnu = \nu \circ (w \mapsto v^N)$, $\ord_w f_j(\nu) = \frac1N \ord_v f_j(\pi(\tnu(v))) = \sum_{\alpha_\bul = j} \ord_w \alpha(\tnu)$. The dimensions follow from Corollary~\ref{cor:rootsum}.
\end{proof}

The next lemma puts the smoothness condition of Lemma~\ref{lem:arcflop}(1) in terms of the parameters $T$, and gives the leading terms of the coordinates of each arm along the arc.

\begin{lem}[The $T_0$-criterion]\label{lem:flopcrit}
Let $\nu$ and $Y$ be as in Lemma~\ref{lem:arcflop}.
\begin{thmlist}
\item $Y$ is smooth along $C$ if and only if, on each arm of $J$, the constant term of the arm polynomial (\S\ref{ss:universal}), which is $T_0$ up to sign (Remark~\ref{rem:Kpresentation}), has order~$1$ along $\nu$.
\item In that case the roots of the arm polynomial of an arm of size $m$ lie in $k[[w^{1/m}]]$, with leading terms $\varrho\,\zeta^i w^{1/m}$ ($1 \le i \le m$), for some $\varrho \neq 0$ and a primitive $m$th root of unity $\zeta$.
\end{thmlist}
\end{lem}

\begin{proof}
(1) By Lemma~\ref{lem:arcflop}(1), $Y$ is smooth along $C$ if and only if $\ord_w u_0(\nu) = 1$ at each singular point of $Y_0$, where $u_0$ is the constant term of the local model there (``Arcs that classify flops'', \S\ref{ss:basechange}). Write $T_0$ for the constant term of the arm polynomial, up to sign. Fix a singular point $p$ of $Y_0$, on an arm $a$ of size $m$. We compare $u_0$ with $T_0$ through the discriminant $\Delta_u$ of $z^m + \dots + u_0$, regarded, like $u_0$, as a function on the base, and the discriminant $\Delta_T$ of the arm polynomial, in $H_l$. The local model determines $u_0$ only up to the formally smooth map of bases, which is why Steps~1--4 are needed: they show that the two discriminants differ by a unit in the completion of $H_l$ at the origin, and Step~5 that each has order $m - 1$ along $\nu$ exactly when its constant term has order $1$.

\emph{Step 1 (irreducible germs).} $\Delta_u$ and $\Delta_T$ are both reduced with irreducible zero germ at the origin. Indeed, the argument for $d_W$ in \S\ref{ss:kleinian}, with the symmetric group of the arm acting transitively on the roots of the arm, shows that the discriminant of the arm polynomial is irreducible in $H_l$, and that of $z^m + \dots + u_0$ in the polynomial ring of the $u_i$. Both are weighted homogeneous in these polynomial rings, for positive weights, so they stay irreducible in the completions at the origin (the lowest-weight terms of a factorisation there would factor them in the polynomial rings), and $\Delta_u$ stays so when pulled back along the formally smooth map of bases.

\emph{Step 2 (singular fibres over closed points).} Write $\widehat{\mathcal O}$ for the completion of $H_l$ at the origin and $p_b$ for the singular point on an arm $b$. We show that over a general point of the image of a root hyperplane of an arm $b$ the fibre has a single singular point, whose $k^\times$-orbit closure meets $Y_0$ only in $p_b$, and that this image is the zero set of the discriminant $\Delta_T^b$ of the arm polynomial of $b$ (so $\Delta_T = \Delta_T^a$). Over a general point $\lambda$ of a hyperplane $\alpha^\perp$ the only complete curve of $\chi^{-1}(\lambda)$ is a $(-2)$-curve $E$ of class $\alpha$ (Proposition~\ref{prop:simres}(3)). The blow-down $\cY_\Gamma \to \cY_l \times_{H_l} \h$ of Proposition~\ref{prop:simres}(2) is the minimal resolution of the fibre $Y_{\pi(\lambda)}$ of $\cY_l$. Hence $Y_{\pi(\lambda)}$ is singular exactly when the blow-down contracts $E$, and then at a single $\mathsf{A}_1$ point $q(\lambda)$, the image of $E$. Let $\mathcal L$ be a line bundle on $\cY_l$ that is ample over $\cX_l$. Its pullback to $\cY_\Gamma$ has degree $0$ on the curves $E_j$, $j \in J$, which the blow-down contracts, and positive degree $d$ on $E_\bul$, which it maps birationally onto $C$; as $E$ has class $\alpha$, the pullback has degree $\alpha_\bul d$ on $E$ (Proposition~\ref{prop:simres}(3)). Hence $E$ is contracted exactly when $\alpha_\bul = 0$. The support of such an $\alpha$ is a connected set of vertices of $J$, so it lies in one arm $b$. As $s \to 0$ in $k^\times$ the curve $s\cdot E$ specialises into a connected union of curves of that arm in $\chi^{-1}(0)$, contracted to $p_b$. Hence $s\cdot q(\lambda) = q(s\lambda) \to p_b$, and the closure $\overline Q$ of the orbit $k^\times\cdot q(\lambda)$, the image of $\A^1 \to \cY_l$, $s \mapsto s\cdot q(\lambda)$, meets $Y_0$ only in $p_b$. For a root $\alpha$ of the arm $b$, the image of $\alpha^\perp$ in $\Spec H_l$ is $V(\Delta_T^b)$: the roots of the arm are the differences $\tau_i - \tau_{i'}$ of its coordinates (``Generators of $H_l$'', \S\ref{ss:universal}), and $\Delta_T^b$ vanishes where two of them coincide.

\emph{Step 3 ($\Delta_u$ is an associate of some $\Delta_T^b$).} Let $\mathscr S \subset \cY_l$ be the relative singular locus of $\cY_l \to \Spec H_l$, the closed subscheme defined by the Fitting ideal $\operatorname{Fitt}_2(\Omega_{\cY_l/H_l})$ (the fibres are surfaces). It is $k^\times$-stable, its closed points are the singular points of the fibres over closed points, and $\mathscr S \to \Spec H_l$ is quasi-finite, as its fibres over closed points are finite and the quasi-finite locus is open (Zariski's Main Theorem) and contains every closed point of $\mathscr S$, so it is all of $\mathscr S$. Fitting ideals commute with base change, so the completion of $\cO_{\mathscr S,p}$ is the relative singular locus of the completed family at $p$, that is, the pullback along the formally smooth map of bases of the critical locus $\{x = y = 0,\ \Psi = \Psi' = 0\}$ of the local model, where $\Psi = z^m + u_{m-2}z^{m-2} + \dots + u_0$. That critical locus, $\{\Psi = \Psi' = 0\}$ in the $(z, u)$-space, is a regular scheme of codimension one in the base (its equations are solved for $u_0$ and $u_1$, or for $u_0$ and $z$ if $m = 2$), finite over the base with image $V(\Delta_u)$, as $\Delta_u$ vanishes exactly where $\Psi$ has a repeated root. Its closed fibre is the single point $z = 0$. Hence $\widehat{\cO}_{\mathscr S,p}$ is a regular local domain of dimension $\dim \widehat{\mathcal O} - 1$, finite over $\widehat{\mathcal O}$, and the image of $\Spec \widehat{\cO}_{\mathscr S,p}$ in $\Spec \widehat{\mathcal O}$ is $V(\Delta_u)$. In particular $\cO_{\mathscr S,p}$ is a domain, so a unique irreducible component $\mathscr S_1$ of $\mathscr S$ passes through $p$. The closure of its image in $\Spec H_l$ is $V(\mathfrak p)$ for the prime $\mathfrak p = \ker(H_l \to \widehat{\cO}_{\mathscr S,p})$, of dimension $\dim \mathscr S_1 = \dim H_l - 1$. Thus $\mathfrak p = (\Delta)$ for an irreducible $\Delta \in H_l$. The fibres over the image of $\mathscr S_1$ are singular, so $\Delta$ divides the discriminant $d_W \in k[\h]^W \subseteq H_l$ of \S\ref{ss:kleinian}. The irreducible factors of $d_W$ in $H_l$ are, up to units, the $f_j$ and the $\Delta_T^b$, as every positive root either has $\alpha_\bul \neq 0$ or lies in one arm. The fibre over a general point of $D_j = V(f_j)$ is smooth, by Step~2, as $\alpha_\bul = j \neq 0$ there; so $\Delta = \Delta_T^b$ for some arm $b$. Then $\Delta_T^b \in \mathfrak p$ vanishes on $\Spec \widehat{\cO}_{\mathscr S,p}$, so $\Delta_T^b \in \sqrt{(\Delta_u)} = (\Delta_u)$, and as $\Delta_T^b$ is irreducible in $\widehat{\mathcal O}$ and $\Delta_u$ is not a unit, $\Delta_u$ and $\Delta_T^b$ are associates.

\emph{Step 4 ($b = a$).} Let $\lambda$ be a general point of $\alpha^\perp$ for a root $\alpha$ of the arm $b$. The arc $s \mapsto \pi(s\lambda)$ in $\Spec H_l$, centred at the origin, lies in $V(\Delta_T^b)$. Hence the point of $\Spec \widehat{\mathcal O}$ that its generic point defines lies in $V(\Delta_T^b) = V(\Delta_u)$, the image of the finite morphism $\Spec \widehat{\cO}_{\mathscr S,p} \to \Spec \widehat{\mathcal O}$. Lifting that point to $\Spec \widehat{\cO}_{\mathscr S,p}$ and extending over the valuation ring gives, for some $e \ge 1$, an arc $\sigma\colon \Spec k[[s]] \to \mathscr S$ with $\sigma(0) = p$ lying over the arc $s \mapsto \pi(s^e\lambda)$. The points of $\mathscr S$ over the punctured curve $\{\pi(s\lambda) : s \neq 0\}$ are the singular points $q(s\lambda) = s\cdot q(\lambda)$ of its fibres, so the generic point of $\sigma$ lies in $\overline Q$, and $p = \sigma(0) \in \overline Q \cap Y_0 = \{p_b\}$. Hence $b = a$, and $\Delta_u$ and $\Delta_T$ differ by a unit.

\emph{Step 5 (orders).} If $u_{m-2}, \dots, u_0 \in wk[[w]]$, the discriminant of $z^m + \dots + u_0$ has order at least $m - 1$, with equality if and only if $\ord_w u_0 = 1$: it is weighted homogeneous of weight $m(m-1)$ when $u_i$ has weight $m - i$, so $u_0^{m-1}$, with coefficient $\pm m^m$, is its only monomial of degree $m - 1$, and all others have higher degree. Applying this to both discriminants proves the criterion.

(2) The coefficients of the arm polynomial are, up to sign, elementary symmetric functions of positive degree in the coordinates of the arm, so they lie in $wk[[w]]$ along $\nu$; as $T_0$ has order $1$, the arm polynomial is an Eisenstein polynomial over $k[[w]]$. Hence it is irreducible over $k((w))$, and its roots generate the extension of degree $m$, which is $k((w^{1/m}))$. They are integral over $k[[w]]$, so lie in $k[[w^{1/m}]]$, and by the Newton polygon each has order $1/m$, with leading coefficient an $m$th root of the leading coefficient of $\pm T_0$. As $w^{1/m} \mapsto \zeta w^{1/m}$ permutes the roots, these leading coefficients are the $\varrho\,\zeta^i$.
\end{proof}

\subsection{The root formula in length 2}\label{ss:rootformula2}

We now work out the root formula for Example~\ref{ex:running}, and determine the GV invariants of the flops of length $2$. Take $\Gamma = \mathsf{D}_4$ with $\bul$ the central vertex, $J = \{1, 3, 4\}$ and $W_J = (\Z/2)^3$. To compare with \cite[Example~4.23]{K}, we label the arms of $J$ so that the loops $\beta$, $\gamma$, $\delta$ at $\bul$ belong to the vertices $1$, $3$, $4$ (Figures~\ref{fig:pairs} and~\ref{fig:laufer}). Their parameters are $(T_0^\beta, T_0^\gamma, T_0^\delta) = \frac14(t_1^2, t_3^2, t_4^2)$, and on the arm of vertex $i$ the arm coordinates of Remark~\ref{rem:Kpresentation} are $\pm\frac12 t_i$. In the coordinates $t_i$ of \S\ref{ss:lie}, the reflection $s_i$ ($i \in J$) sends $t_i$ to $-t_i$ and $t_\bul$ to $t_\bul + t_i$, and fixes the other coordinates. The highest root $\theta$ has coefficient $2$ at $\bul$ and $1$ at each vertex of $J$, and it is orthogonal to the simple roots of $J$, so
\[
t := \theta(\lambda) = 2t_\bul + t_1 + t_3 + t_4, \qquad H_2 = k[t] \otimes P, \qquad P = k[t_1^2,\, t_3^2,\, t_4^2].
\]
This $t$ is $-\lambda_0$, because the constraint of \S\ref{ss:simres} that $\lambda$ vanish on the minimal imaginary root reads $\lambda_0 + \theta(\lambda) = 0$, while the invariant form called $t$ in \cite{K} is $\lambda_0$. (Only in lengths $1$ and $2$ is $t$ a multiple of $\theta(\lambda)$: in lengths $3$--$6$ the highest root is not orthogonal to the simple roots of $J$, and $t$ is the form of \S\ref{ss:lie}.) The generators of $P$ have degrees $2, 2, 2$, with product $8 = |W_J|$, as in \S\ref{ss:universal}. The positive roots with $\alpha_\bul \neq 0$ are $\theta$, with $\alpha_\bul = 2$, and the eight roots $\alpha_U$ of Example~\ref{ex:running}, with $\alpha_\bul = 1$. Since
\[
\alpha_U(\lambda) = \tfrac12\,\bigl(t + \varepsilon_1 t_1 + \varepsilon_3 t_3 + \varepsilon_4 t_4\bigr), \qquad \varepsilon_i = +1 \text{ if } i \in U,\ -1 \text{ otherwise},
\]
the two components of $D$ in $\Spec H_2$ are $D_2 = \{t = 0\}$ and $D_1 = \{f_1 = 0\}$, where
\[
f_1 = \prod_{\varepsilon \in \{\pm 1\}^3} \tfrac12\,\bigl(t + \varepsilon_1 t_1 + \varepsilon_3 t_3 + \varepsilon_4 t_4\bigr) \;\in\; k[t,\, t_1^2,\, t_3^2,\, t_4^2]
\]
is an octic in $H_2$ that splits into the eight root forms on $\h$. For an arc $\nu(w) = (t(w), T_0^\beta(w),$ $T_0^\gamma(w), T_0^\delta(w))$ through the origin, Corollary~\ref{prop:rootformula} gives
\[
n_2 = \ord_w t(\nu), \qquad n_1 = \ord_w f_1(\nu) = \sum_{\varepsilon} \ord_w\, \tfrac12\bigl(t + \varepsilon \cdot \mathbf t\bigr)(\tnu),
\]
where $\mathbf t = (t_1, t_3, t_4)$ is evaluated on the lift $\tnu$ of Lemma~\ref{lem:lift}. Here $N = 2$ suffices: along any arc each $T_0$ is zero or a power of $w$ times a unit, so the corresponding $t_i = \pm 2T_0^{1/2}$ lies in $k[[w^{1/2}]]$. Individual terms may lie in $\frac12\Z$, but their sum over a $W_J$-orbit is $\ord_w$ of an element of $H_2$, and so an integer.

\begin{lem}[{GV invariants in lengths 1 and 2; recovering \cite[Cor~8.13]{BW}}]\label{lem:length2}
Let $f$ be a flop of length $l \le 2$.
\begin{thmlist}
\item If $l = 1$, its GV invariant $n_1$ is a positive integer, and each $n_1 \ge 1$ occurs.
\item If $l = 2$, its GV invariants are $(2m + 3, m)$ with $m \ge 1$, or $(2i, m)$ with $i \ge 2$ and $m \ge i - 1$, and each of these is realised by an explicit arc.
\end{thmlist}
\end{lem}

\begin{proof}
(1) In length $1$, $J$ is empty and $f_1 = t$. A flop of length $1$ has $n_1 = \ord_w t(\nu)$ by Proposition~\ref{prop:GV}, which is at least $1$ as its classifying arc $\nu$ sends the closed point to the origin. Conversely, for $m \ge 1$ the arc $t = w^m$ classifies a flop with $n_1 = m$, by Lemma~\ref{lem:arcflop}(2), whose smoothness condition is empty here, and Proposition~\ref{prop:GV}.

(2) By Lemma~\ref{lem:flopcrit}(1), the pulled-back $Y$ is smooth if and only if $\ord_w T_0(\nu) = 1$ on each arm. For such arcs $t_i/w^{1/2} \in k[[w]]$ with nonzero constant term $c_i$, and, writing $\mathbf c = (c_1, c_3, c_4)$, $\varepsilon\cdot\mathbf c = 0$ for at most one pair $\pm\varepsilon$. Write $n_2 = \ord_w t$. If $n_2 < \infty$, then $f_2(\nu) = t(\nu)$ is nonzero, and so is $f_1(\nu)$: each $\alpha_U = \frac12(t + \varepsilon\cdot\mathbf t)$ is nonzero along the lifted arc, since $t \in k[[w]]$ is nonzero while $\varepsilon\cdot\mathbf t \in w^{1/2}k[[w]]$. Hence $\nu$ classifies a flop with GV invariants $(n_1, n_2)$, where $n_1 = \ord_w f_1(\nu)$ (Lemma~\ref{lem:arcflop}(2) and Proposition~\ref{prop:GV}). If no pair vanishes, all eight root forms have order $\frac12$, so $n_1 = 4$. If a pair vanishes, let $r \in \{1, 2, \dots\} \cup \{\infty\}$ be $\ord_w$ of $(\varepsilon\cdot\mathbf t)/w^{1/2}$, with $r = \infty$ if this vanishes identically, as for the arc $T_0^\delta = 4w$ below. Then the pair contributes $2\cdot\min(n_2, r + \frac12)$ and the other six $\frac12$ each. This gives $n_1 = 2n_2 + 3$ when $r \ge n_2$, and $n_1 = 2r + 4$ when $r < n_2$. The classifying arc of every flop of length $2$ is such an arc, with $n_2$ finite (Lemmas~\ref{lem:arcs} and~\ref{lem:flopcrit}(1)), so its GV invariants are among the pairs $(4, m)$, $(2m + 3, m)$ and $(2r + 4, m)$ with $1 \le r < m$, which are those of the statement, with $i = 2$ and $i = r + 2$. Each pair is realised by an explicit arc. With $t = w^m$ and $T_0^\beta = T_0^\gamma = w$, taking $T_0^\delta = w$ gives no vanishing pair and $(4, m)$, taking $T_0^\delta = 4w$ gives $(2m + 3, m)$, and taking $T_0^\delta = 4w(1 + w^r)$ with $r < m$ gives $(2r + 4, m)$, since then $(\varepsilon\cdot\mathbf t)/w^{1/2} = \mp 2w^r + \cdots$ for $\varepsilon = \pm(1, 1, -1)$. For instance, a general linear arc $\nu(w) = (aw, a_1w, a_3w, a_4w)$ has $n_2 = 1$ and no vanishing pair, so $(n_1, n_2) = (4, 1)$.
\end{proof}

\subsection{GV invariants of flops}\label{ss:typeE}

Lemmas~\ref{lem:arcflop} and~\ref{lem:flopcrit} and the root formula determine which GV invariants occur in every length, and explicit arcs realise them. Lengths $1$ and $2$ are Lemma~\ref{lem:length2}, and we now treat lengths $3$--$6$, where the work is.

\begin{prop}[GV invariants in lengths 3--6; $=$ Proposition~E]\label{prop:typeE}
Let $f$ be a flop of length $l \ge 3$.
\begin{thmlist}
\item If $l = 3$, its GV invariants are $(6, 3, 1)$ or $(7, n_2, 1)$ with $n_2 \ge 4$, and each of these occurs.
\item If $l = 4$, $5$ or $6$, its GV invariants are $(6, 5, 2, 1)$, $(8, 6, 4, 2, 1)$ or $(6, 6, 4, 3, 2, 1)$ respectively, and these occur. In particular $\dim \Acon$ is $60$, $125$ or $200$, and $\dim Z(\Acon)$ is $14$, $21$ or $22$.
\end{thmlist}
\end{prop}

For one-parameter families linear in the arc parameter, \cite[\S\S4.3--4.6]{CDSV} found one tuple in each of the lengths $3$--$6$ by a count of zero modes (compare Corollary~\ref{prop:rootformula}). These are the GV invariants in the first four rows of Table~\ref{tab:arcs}, and in length $3$ they come from a single family at general values of its parameters. Proposition~\ref{prop:typeE} shows that every flop of length $4$, $5$ or $6$ has the tuple of its length, and that in length $3$ the other tuples are the $(7, n_2, 1)$ with $n_2 \ge 4$. The proof describes the roots with $\alpha_\bul \neq 0$ through their components on the arms of $J$, which Table~\ref{tab:levi} lists. Brown and Wemyss obtain the same list from their classification of contraction algebras, lower bounds on the $n_j$, upper bounds on $\dim \Acon$ and Toda's formula \cite[Cor~1.6]{BW26}.

\begin{table}[ht]
\caption{The roots with $\alpha_\bul = j$ in lengths $3$--$6$. On an arm of size $m$ the arm term of such a root is $\tau_S$, with $|S|$ depending only on $j$ and the arm (proof of Proposition~\ref{prop:typeE}); each entry gives $|S|$ for the arms, in the order of the second row, and the number $\prod \binom{m}{|S|}$ of these roots.}\label{tab:levi}
\begin{tabular}{@{}ccccc@{}}
\toprule
 & $l = 3$ & $l = 4$ & $l = 5$ & $l = 6$ \\
Arm sizes $m$ & $3, 3, 2$ & $4, 3, 2$ & $5, 4$ & $5, 3, 2$ \\
\midrule
$j = 1$ & $(2, 2, 1)$; $18$ & $(3, 2, 1)$; $24$ & $(3, 3)$; $40$ & $(4, 2, 1)$; $30$ \\
$j = 2$ & $(1, 1, 0)$; $9$ & $(2, 1, 0)$; $18$ & $(1, 2)$; $30$ & $(3, 1, 0)$; $30$ \\
$j = 3$ & $(0, 0, 1)$; $2$ & $(1, 0, 1)$; $8$ & $(4, 1)$; $20$ & $(2, 0, 1)$; $20$ \\
$j = 4$ & & $(0, 2, 0)$; $3$ & $(2, 0)$; $10$ & $(1, 2, 0)$; $15$ \\
$j = 5$ & & & $(0, 3)$; $4$ & $(0, 1, 1)$; $6$ \\
$j = 6$ & & & & $(4, 0, 0)$; $5$ \\
\bottomrule
\end{tabular}
\end{table}

\begin{proof}
\emph{Outline.} By the root formula (Corollary~\ref{prop:rootformula}), $n_j$ is the sum of the orders along $\tnu$ of the roots with $\alpha_\bul = j$. We write each root through its terms on the arms of $J$ (\emph{Arm terms} below). On an arm of size $m$ the coordinates have order $\frac1m$ along $\tnu$ (Lemma~\ref{lem:flopcrit}(2)), so a root has order $\frac1m$ for the largest size $m$ of an arm on which its term is nonzero, unless leading terms cancel. They cancel only for six of the roots with $\alpha_\bul = 2$ in length $4$ (\emph{Orders} below), and in length $3$ exactly when $T_0^\beta$ and $T_0^\gamma$ have opposite leading coefficients (\emph{Length~$3$} below), which gives the two cases of~(1).

\emph{Arm terms.} The linear forms on $\h$ vanishing on $\varpi_\bul^\vee$ are spanned by the simple roots $t_i$, $i \in J$, and $t$ is orthogonal to them for the invariant form, being $W_J$-invariant. Hence a positive root $\alpha$ is $\alpha_\bul\, t/t(\varpi_\bul^\vee)$ plus the sum, over the arms of $J$, of its orthogonal projections to the spans of the roots of the arms, its \emph{arm terms}. On an arm of type $\mathsf{A}_{m-1}$ the coordinates $\tau_1, \dots, \tau_m$ of Remark~\ref{rem:Kpresentation} satisfy $\sum_i \tau_i = 0$, and after changing their common sign if necessary they can be numbered so that the simple roots of the arm are the $\tau_i - \tau_{i+1}$, starting from the end of the arm nearer to the vertex joined to $\bul$ \cite[Examples~4.24--4.27]{K}. Let $i_0$ be the position of that vertex, so $i_0 = 1$, except on the arm of type $\mathsf{A}_4$ in length $5$, where $i_0 = 2$. For $S \subseteq \{1, \dots, m\}$ put $\tau_S = \sum_{i \in S}\tau_i$. We claim that if $\alpha_\bul = j > 0$, the arm term of $\alpha$ is $\tau_S$ with $|S| \equiv -j\,i_0 \pmod m$. Write it as $\sum_i c_i\tau_i$. For $i \neq i'$, $\tau_i - \tau_{i'}$ is a root of $J$, so it is not $\pm\alpha$, and as $\Gamma$ is simply laced, $(\alpha, \tau_i - \tau_{i'}) \in \{-1, 0, 1\}$. This pairing is $c_i - c_{i'}$, since $(\tau_a, \tau_b) = \delta_{ab} - 1/m$ and the other summands of $\alpha$ are orthogonal to the arm. Thus the $c_i$ take at most two values, which differ by $1$, and as $\sum_i \tau_i = 0$ the arm term is $\tau_S$, with $S$ the set where $c_i$ takes the larger value. The arm term of $t_\bul$ is $-(\tau_1 + \dots + \tau_{i_0})$, as its pairing with $\tau_i - \tau_{i+1}$ is $-1$ for $i = i_0$ and $0$ otherwise. Hence, writing $\alpha = j\,t_\bul + \sum_{i \in J} a_i t_i$ with integers $a_i$, the arm term of $\alpha$ is $-j(\tau_1 + \dots + \tau_{i_0})$ plus an integral combination of the $\tau_i - \tau_{i+1}$. An integral combination of the $\tau_i$ representing a given form is unique up to adding the same integer to every coefficient, so the sum of the coefficients is well defined modulo $m$; it is $-j\,i_0$ here, and $|S|$ for $\tau_S$. This proves the claim.

Hence a root with $\alpha_\bul = j$ is determined by a subset $S$ for each arm, with $|S| \in \{0, \dots, m - 1\}$ fixed by the claim (and $\tau_S = 0$ when $S = \emptyset$). The group $W_J$ permutes the coordinates of each arm, so it acts transitively on these tuples of subsets, and it preserves the nonempty set of roots with $\alpha_\bul = j$. Hence these roots correspond to the tuples, and their number is $\prod\binom{m}{|S|}$, the product over the arms, as in Table~\ref{tab:roots}. Table~\ref{tab:levi} lists the values of $|S|$.

\emph{Orders.} Let $\nu$ be the classifying arc of $f$ and $\tnu$ its lift (Lemma~\ref{lem:lift}). Along $\tnu$, $t$ has order at least $1$ (possibly infinite). The coordinates of an arm of size $m$ are, up to a common sign, the roots of its arm polynomial, so by Lemma~\ref{lem:flopcrit}(2) they have leading terms $\pm\varrho\,\zeta^i w^{1/m}$, $1 \le i \le m$, in some order. Hence, for $0 < |S| < m$, $\tau_S$ has order $1/m$ unless the corresponding sum of $|S|$ distinct $m$th roots of unity vanishes. For $m \le 5$ this happens only when $m = 4$ and the two roots of unity are opposite, for two of the six subsets with $|S| = 2$ (for $m$ prime no proper nonempty set of $m$th roots of unity has sum zero); then $\tau_S$ lies in $k[[w^{1/4}]]$, so it has order at least $\frac12$.

\emph{Lengths $4$, $5$ and $6$.} Here the arms have different sizes, listed in decreasing order in Table~\ref{tab:levi}, and in each entry the first nonzero value of $|S|$ belongs to an arm of size $m \ge 3$. Its arm term $\tau_S$ has order exactly $1/m$, unless $l = 4$, $j = 2$ and $\tau_S$ cancels. Every other arm term is zero, or has order $1/m'$ with $m' < m$, or order at least $\frac12$, and the multiple of $t$ has order at least $1$, all larger than $1/m$ as $m \ge 3$. Hence the root form has order $1/m$. In length $4$ with $j = 2$, the arm term on the arm of size $4$ cancels for $6$ of the $18$ roots, which then have order $\frac13$, from the arm of size $3$. Summing,
\[
(n_1, \dots, n_l) = \bigl(\tfrac{24}{4}, \tfrac{12}{4} + \tfrac63, \tfrac84, \tfrac33\bigr), \quad \bigl(\tfrac{40}{5}, \tfrac{30}{5}, \tfrac{20}{5}, \tfrac{10}{5}, \tfrac44\bigr), \quad \bigl(\tfrac{30}{5}, \tfrac{30}{5}, \tfrac{20}{5}, \tfrac{15}{5}, \tfrac63, \tfrac55\bigr)
\]
in lengths $4$, $5$ and $6$, which are the tuples in~(2). The same computation applies to any arc satisfying the condition of Lemma~\ref{lem:flopcrit}(1), such as the arcs of \cite[Examples~5.4, 5.5]{K} below, and shows that every root form has finite order along it. Hence such an arc classifies a flop (Lemma~\ref{lem:arcflop}(2)), and these tuples occur. The dimensions follow from Corollary~\ref{cor:toda} and Theorem~B(1).

\emph{Length $3$.} Here the arms have sizes $3$, $3$ and $2$, with coordinates $\tau^\beta$, $\tau^\gamma$ and $\tau^\delta$ (Figure~\ref{fig:pairs}), of orders $\frac13$, $\frac13$ and $\frac12$. By Table~\ref{tab:levi}, as $\tau_S = -\tau_{S^c}$, the root forms with $\alpha_\bul = 3$, $2$ and $1$ are, up to multiples of $t$, $\pm\tau^\delta$, $\tau^\beta + \tau^\gamma$ and $-(\tau^\beta + \tau^\gamma) \pm \tau^\delta$, where $\tau^\beta$, $\tau^\gamma$, $\tau^\delta$ run over the coordinates of the three arms. Hence $n_3 = 2\cdot\frac12 = 1$. In the coordinates of \cite[Example~4.24]{K}, which are already numbered as above, the $-\tau^\beta$ are the roots of $x^3 - T_1^\beta x - T_0^\beta$, and the $-\tau^\gamma$ those of $x^3 - T_1^\gamma x - T_0^\gamma$ (Remark~\ref{rem:Kpresentation}). Therefore the leading coefficient of $\tau^\beta + \tau^\gamma$ vanishes for some of the nine pairs exactly when $T_0^\beta$ and $T_0^\gamma$ have opposite leading coefficients, and then for exactly three pairs. Otherwise every root form with $\alpha_\bul = 1$ or $2$ has order $\frac13$, and the GV invariants are $(6, 3, 1)$. If they cancel, the three cancelling sums lie in $k[[w^{1/3}]]$, so they have order at least $\frac23$. The six root forms with $\alpha_\bul = 1$ containing them then have order exactly $\frac12$, because of $\pm\tau^\delta$, and the other twelve have order $\frac13$, so $n_1 = 4 + 3 = 7$. Also $n_2 \ge 6\cdot\frac13 + 3\cdot\frac23 = 4$.

Each case occurs. The first occurs for the arc of \cite[Example~5.3]{K}, which has $t = T_1^\beta = T_1^\gamma = 0$ and $T_0^\beta = T_0^\gamma = T_0^\delta = w$. It satisfies the condition of Lemma~\ref{lem:flopcrit}(1), and as $T_0^\beta$ and $T_0^\gamma$ have the same leading coefficient, every root form has order $\frac13$ or $\frac12$ along it, by the computation above. Hence it classifies a flop (Lemma~\ref{lem:arcflop}(2)) with GV invariants $(6, 3, 1)$. For the second, take $T_0^\beta = T_0^\delta = w$ and $T_0^\gamma = -w$, with the other parameters zero. The arm polynomials are $x^3 - w$, $x^3 + w$ and $x^2 - w$, with roots $\zeta^i w^{1/3}$, $-\zeta^i w^{1/3}$ and $\pm w^{1/2}$, and the three cancelling sums vanish identically. Now change the arc in one of three ways, keeping $t = 0$ in (i) and (iii); in (ii) the three sums $\tau^\beta + \tau^\gamma$ still vanish identically.
\begin{enumerate}[label=(\roman*)]
\item Set $T_1^\gamma = w^a$, with $a \ge 1$. The roots of $x^3 - w^a x + w$ are $-\zeta^i w^{1/3} - \frac13\zeta^{2i} w^{a - 1/3} + \cdots$, so the three sums acquire order $a - \frac13$.
\item Set $t = w^a$, with $a \ge 1$. The roots do not change, and the three cancelling root forms with $\alpha_\bul = 2$ become nonzero multiples of $w^a$, of order $a$.
\item Replace $T_0^\gamma$ by $-w - w^a$, with $a \ge 2$. The roots of $x^3 + w + w^a$ are $-\zeta^i w^{1/3}(1 + w^{a-1})^{1/3}$, so the three sums acquire order $a - \frac23$.
\end{enumerate}
In each case every root form has finite order along the arc, so it classifies a flop by Lemmas~\ref{lem:flopcrit}(1) and~\ref{lem:arcflop}(2), and $n_2 = 6\cdot\frac13 + 3\cdot(\text{that order})$ is $3a + 1$, $3a + 2$ and $3a$ respectively. Together these take every value $n_2 \ge 4$.
\end{proof}

\paragraph{Explicit arcs.} The examples of \cite[Examples~5.3--5.5]{K} are arcs in $\Spec H_l$, given by ring maps $H_l \to k[w]$ with $t \mapsto 0$ (the arc parameter is called $T$ there), so the root formula applies to them directly. The coordinates $\tau$ of each arm (Figure~\ref{fig:pairs}) are, up to a common sign, the roots of its arm polynomial (Remark~\ref{rem:Kpresentation}), whose coefficients are read off from the arc. For $l = 3, 4, 6$ each arm polynomial is $x^m \pm w$, and for $l = 5$ so is that of the arm of size $4$, while the arm of size $5$ has $T_1 = T_0 = w$ and $T_2 = T_3 = 0$. Killing $e_0$ in the relations of the presentation along these arcs gives algebras whose local blocks at the origin are the contraction algebras, where $b$ and $c$ are the images of the loops $\beta$ and $\gamma$, up to sign:
\begin{itemize}
\item $b^3 = c^3 = (b + c)^2$ in length $3$;
\item $b^2 = c^4 = -(b + c)^3$ in length $4$;
\item $cbc + c^2b + cb^3 = (c + b^2)^2 + bcb = -b^4$ in length $5$ (where $-Te_4$ becomes $-b^4$, as killing $e_0$ kills $a^*a$ in the relation $a^*a = \delta^4 - T_0^\delta e_4$ of \cite[Thm~1.3(5)]{K}, and $\beta = \delta$ along this arc);
\item $b^2 = c^3 = -(b + c)^5$ in length $6$.
\end{itemize}
As stated there, these arcs classify flops (proof of Proposition~\ref{prop:typeE}). In length~$5$ these relations differ from those displayed in its Example~5.5, which belong to a different arc (Appendix~\ref{app:corrections}).

Table~\ref{tab:arcs} lists the GV invariants of these arcs, given by the proof of Proposition~\ref{prop:typeE}, and the resulting dimensions of $\Acon$ and its centre (Corollary~\ref{cor:toda} and Theorem~B(1)). The GV invariants of the first row are also computed in \cite[Example~5.3]{K}.

\begin{table}[ht]
\caption{The arcs of \cite[Examples~5.3--5.5]{K}, and a second arc in length $3$.}\label{tab:arcs}
\begin{tabular}{@{}ccccc@{}}
\toprule
Length & Arc & $(n_1, \dots, n_l)$ & $\dim \Acon$ & $\dim Z(\Acon)$ \\
\midrule
3 ($\mathsf{E}_6$) & \cite[5.3]{K} & $(6, 3, 1)$ & $27$ & $10$ \\
4 ($\mathsf{E}_7$) & \cite[5.4]{K} & $(6, 5, 2, 1)$ & $60$ & $14$ \\
5 ($\mathsf{E}_8$) & \cite[5.5]{K} & $(8, 6, 4, 2, 1)$ & $125$ & $21$ \\
6 ($\mathsf{E}_8$) & \cite[5.4]{K} & $(6, 6, 4, 3, 2, 1)$ & $200$ & $22$ \\
3 ($\mathsf{E}_6$) & family~(i), $a = 1$ & $(7, 4, 1)$ & $32$ & $12$ \\
\bottomrule
\end{tabular}
\end{table}

The last row comes from the length-$3$ arc $t = 0$, $T_0^\beta = w$, $T_1^\beta = 0$, $T_0^\gamma = -w$, $T_1^\gamma = w$, $T_0^\delta = w$, whose arm polynomials are $x^3 - w$, $x^3 - wx + w$ and $x^2 - w$. It is the case $a = 1$ of family~(i) in the proof of Proposition~\ref{prop:typeE}, so it classifies a flop. Its contraction algebra is the local block at the origin of $k\langle b, c\rangle/\bigl(b^3 - (b + c)^2,\ c^3 - (b + c)^2c + (b + c)^2\bigr)$.

\paragraph{The flops $\mathsf{E}_{6,n}$.} Brown and Wemyss announced in \cite[\S1.5, Table~1]{BW}, and prove in \cite[Thms~1.2, 1.3]{BW26}, that their normal forms $\mathsf{E}_{6,n} = x^3 + xy^3 + y^n$, $n \ge 4$, are realised by flops of type $\mathsf{E}_6$ (see also \cite[Thm~1.7]{Wem}), that is, their Jacobi algebras are contraction algebras of flops of length $3$ (\S\ref{ss:universal}). For every $p$ with $\ord p = n \ge 4$, \cite[Cor~4.12]{Pan} constructs a flop whose contraction algebra is the Jacobi algebra of $x^3 + xy^3 + p(y)$, and Lemma~\ref{lem:length2}, Proposition~\ref{prop:typeE} and Toda's formula give its length and GV invariants. Indeed $n_1$ is the dimension of the abelianisation (\S\ref{ss:flops}), which for a Jacobi algebra is the commutative Jacobian algebra of the same potential. Hence $n_1$ is the Milnor number of $x^3 + xy^3 + p(y)$. This singularity is semi-quasihomogeneous, with principal part $x^3 + c\,y^4$ ($c \neq 0$ the coefficient of $y^4$ in $p$) for $n = 4$ and $x^3 + xy^3$ for $n \ge 5$, so it is the plane curve singularity $\mathsf{E}_6$, respectively $\mathsf{E}_7$, and $n_1 = 6$, respectively $7$. The Jacobi algebra of $x^3 + xy^3 + p(y)$ has dimension $27$ for $n = 4$ and $4(n + 3)$ for $n \ge 5$ \cite[Thm~3.1]{Pan} (for $p = y^n$, $n \ge 5$, see also \cite[\S1.3]{BW}), and by Toda's formula (Corollary~\ref{cor:toda}) this is $\sum_j j^2 n_j$. Among the GV invariants that Lemma~\ref{lem:length2} and Proposition~\ref{prop:typeE} allow in lengths $1$--$6$, only $(6, 3, 1)$ has $n_1 = 6$ and $\sum_j j^2 n_j = 27$, and only $(7, n - 1, 1)$ has $n_1 = 7$ and $\sum_j j^2 n_j = 4(n + 3)$. Hence the flop has length $3$, with these GV invariants (see also \cite[Thm~5.3]{Pan} and \cite[Cor~1.6]{BW26}). 

For polynomial $p$, Smoktunowicz and Smoktunowicz construct explicit flat deformations of this Jacobi algebra to $k^6 \times M_2(k)^3 \times M_3(k)$ if $n = 4$ and to $k^7 \times M_2(k)^{n-1} \times M_3(k)$ if $n \ge 5$ \cite[Thms~4.3, 7.1]{SS}, assuming that it is local of the dimension above, which holds as it is a contraction algebra of that dimension. These are the algebras $\SGV$ for the invariants above, so by Theorem~C the Jacobi algebra deforms to no other semisimple algebra. For these algebras this is \cite[Thm~3.1]{Pan}, proved there by the argument of Section~\ref{sec:reduction}: the centre has the dimension of the semisimple fibre, and $h^1 = h^2$, so the algebra is a smooth point.

\subsection{Laufer flops}\label{ss:laufer}

Finally we treat the Laufer flops, whose contraction algebras are known explicitly, and realise the first of them by an arc. The Laufer flops have length $2$ ($\Gamma = \mathsf{D}_4$), and by \cite[Example~3.14]{DW} their contraction algebras are
\[
L_m = k\langle x, y\rangle / (xy + yx,\ x^{2m+1} + y^2), \qquad m \ge 1.
\]
$L_1$ is the algebra of Section~\ref{sec:centremap}.

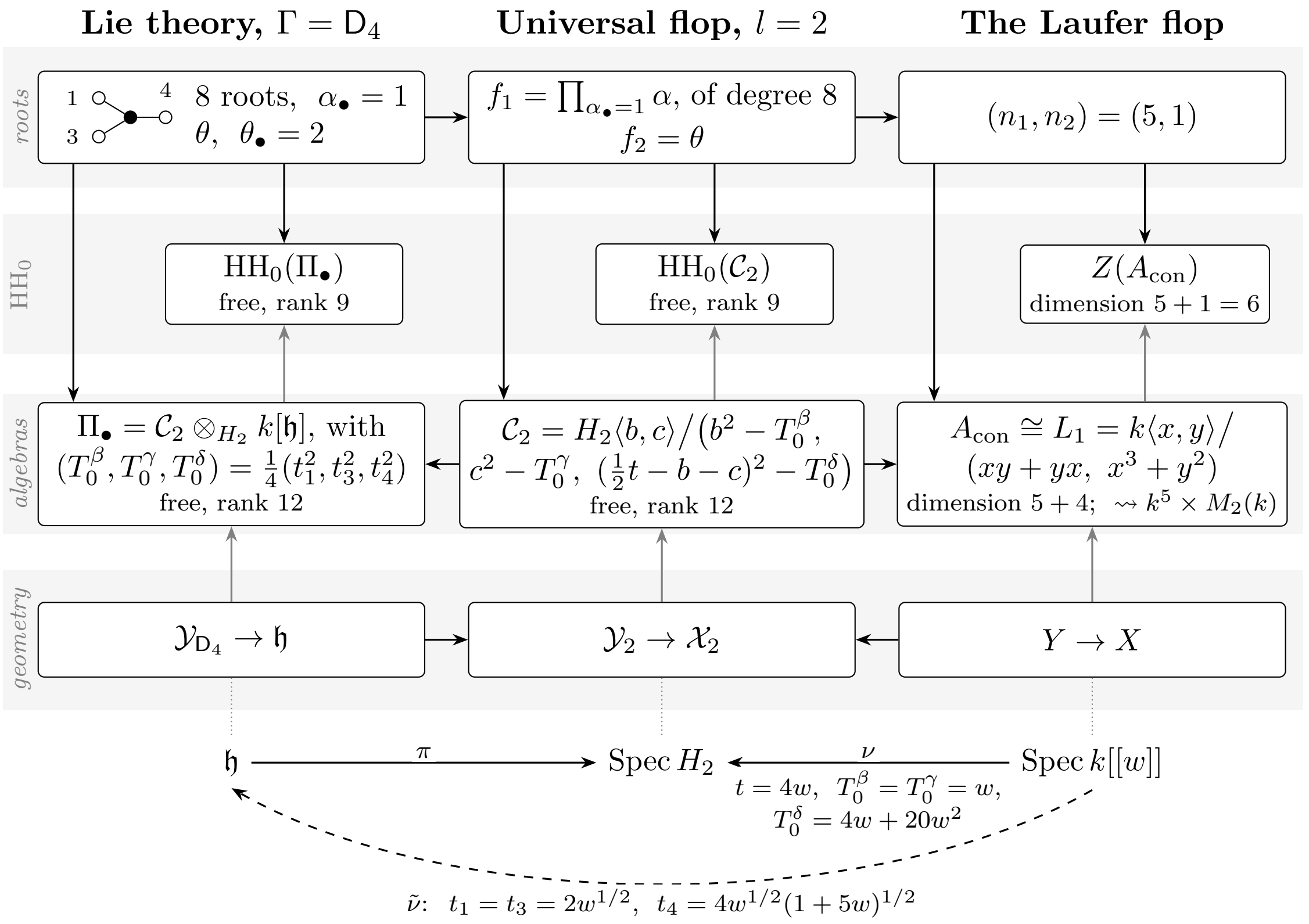
\begin{figure}[t]
\centering
\begin{tikzpicture}[
  font=\small,
  box/.style={draw, rounded corners=2pt, fill=white, align=center, inner sep=3pt,
              minimum height=0.85cm, minimum width=4.4cm},
  hbox/.style={box, minimum width=2.7cm},
  lab/.style={font=\scriptsize, align=center, fill=white, inner sep=1.2pt},
  arr/.style={-{Stealth[length=5pt]}, semithick},
  dyn/.style={circle, draw, inner sep=0pt, minimum size=4pt},
]
\def\xa{0} \def\xb{4.9} \def\xc{9.8}
\def\yr{0} \def\yh{-1.9} \def\ya{-3.95} \def\yg{-5.95} \def\yb{-7.35}
\def\hoff{0.6}   %
\def\lane{-1.8}  %

\node[font=\bfseries] at (\xa,1.05) {Lie theory, $\Gamma = \mathsf{D}_4$};
\node[font=\bfseries] at (\xb,1.05) {Universal flop, $l = 2$};
\node[font=\bfseries] at (\xc,1.05) {The Laufer flop};

\node[box, minimum height=1.05cm] (R1) at (\xa,\yr) {};
\begin{scope}[shift={(\xa-1.15,\yr)}]
  \node[dyn, fill=black] (d2) at (0,0) {};
  \node[dyn] (d1) at (-0.36,0.22) {};
  \node[dyn] (d3) at (-0.36,-0.22) {};
  \node[dyn] (d4) at (0.4,0) {};
  \foreach \v in {1,3,4} \draw (d2) -- (d\v);
  \node[font=\tiny, left=0.5pt of d1] {$1$};
  \node[font=\tiny, left=0.5pt of d3] {$3$};
  \node[font=\tiny, above=0.5pt of d4] {$4$};
\end{scope}
\node[align=left] at (\xa+0.8,\yr)
  {$8$ roots, \ $\alpha_\bul = 1$\\ $\theta$, \ $\theta_\bul = 2$};
\node[box, minimum height=1.05cm] (R2) at (\xb,\yr)
  {$f_1 = \prod\nolimits_{\alpha_\bul = 1} \alpha$, of degree $8$\\[2pt] $f_2 = \theta$};
\node[box, minimum height=1.05cm] (R3) at (\xc,\yr) {$(n_1, n_2) = (5, 1)$};
\node[hbox] (H1) at (\xa+\hoff,\yh) {$\HH_0(\Pi_\bul)$\\[-1pt] {\scriptsize free, rank $9$}};
\node[hbox] (H2) at (\xb+\hoff,\yh) {$\HH_0(\cC_2)$\\[-1pt] {\scriptsize free, rank $9$}};
\node[hbox] (H3) at (\xc+\hoff,\yh) {$Z(\Acon)$\\[-1pt] {\scriptsize dimension $5 + 1 = 6$}};
\node[box, minimum height=1.4cm] (A1) at (\xa,\ya) {$\Pi_\bul = \cC_2 \otimes_{H_2} k[\h]$, with\\
  $(T_0^\beta, T_0^\gamma, T_0^\delta) = \frac14(t_1^2, t_3^2, t_4^2)$\\[-1pt]
  {\scriptsize free, rank $12$}};
\node[box, minimum height=1.4cm] (A2) at (\xb,\ya) {$\cC_2 = H_2\langle b, c\rangle \big/ \bigl(b^2 - T_0^\beta,$\\
  $c^2 - T_0^\gamma,\ (\tfrac12 t - b - c)^2 - T_0^\delta\bigr)$\\[-1pt]
  {\scriptsize free, rank $12$}};
\node[box, minimum height=1.4cm] (A3) at (\xc,\ya) {$\Acon \cong L_1 = k\langle x, y\rangle \big/$\\
  $(xy + yx,\ x^3 + y^2)$\\[-1pt]
  {\scriptsize dimension $5 + 4$; \ $\rightsquigarrow k^5 \times M_2(k)$}};
\node[box] (G1) at (\xa,\yg) {$\cY_{\mathsf{D}_4} \to \h$};
\node[box] (G2) at (\xb,\yg) {$\cY_2 \to \cX_2$};
\node[box] (G3) at (\xc,\yg) {$Y \to X$};

\node (B1) at (\xa,\yb) {$\h$};
\node (B2) at (\xb,\yb) {$\Spec H_2$};
\node (B3) at (\xc,\yb) {$\Spec k[[w]]$};
\draw[arr] (B1) -- node[lab, above] {$\pi$} (B2);
\draw[arr] (B3) -- node[lab, above] {$\nu$}
  node[lab, below, yshift=-1pt] {$t = 4w$, \ $T_0^\beta = T_0^\gamma = w$,\\ $T_0^\delta = 4w + 20w^2$} (B2);
\draw[arr, dashed] (B3.south) to[out=-150, in=-30, looseness=0.75]
  node[lab, below] {$\tnu$: \ $t_1 = t_3 = 2w^{1/2}$, \ $t_4 = 4w^{1/2}(1 + 5w)^{1/2}$} (B1.south);
\foreach \c in {1,2,3} \draw[densely dotted, gray] (G\c.south) -- (B\c.north);

\draw[arr] (R1) -- (R2);
\draw[arr] (R2) -- (R3);
\draw[arr] (A2) -- (A1);
\draw[arr] (A2) -- (A3);
\draw[arr] (G1) -- (G2);
\draw[arr] (G3) -- (G2);

\foreach \c/\x in {1/\xa, 2/\xb, 3/\xc}{
  \draw[arr] (R\c.south -| {\x+\lane,0}) -- (A\c.north -| {\x+\lane,0});
  \draw[arr] (R\c.south -| H\c) -- (H\c);
  \draw[arr, gray] (A\c.north -| H\c) -- (H\c);
  \draw[arr, gray] (G\c) -- (A\c);}

\begin{scope}[on background layer]
  \foreach \y/\t in {\yr/roots, \yh/$\HH_0$, \ya/algebras, \yg/geometry}{
    \fill[black!4] (-2.6,\y-0.8) rectangle (\xc+2.4,\y+0.8);
    \node[font=\scriptsize\itshape, rotate=90, gray] at (-2.4,\y) {\t};}
\end{scope}
\end{tikzpicture}
\caption{The Laufer flop of Example~\ref{ex:lauferarc} in the positions of Figure~\ref{fig:story}, with the arc $\nu$ and its lift $\tnu$. Ranks are over $k[\h_J]$ (left) and $P$ (middle), and $\theta$ is the highest root.}\label{fig:laufer}
\end{figure}

\begin{ex}[The Laufer flop by an arc]\label{ex:lauferarc}
We realise the Laufer flop of Proposition~D by an explicit arc and follow it through Section~\ref{sec:flops}; Figure~\ref{fig:laufer} places it in Figure~\ref{fig:story}. The maps used for the Laufer flops in \cite[Example~5.2]{K} and \cite[\S5]{Kaw} are maps into the total space $\cX_2$, not arcs in $\Spec H_2$ as Lemma~\ref{lem:basechange} needs: they set parameters of the base equal to coordinates of the total space, so that $X$ is cut out of $\cX_2$ by three equations. For $l = 2$, \cite[Thm~1.3]{K} presents $A_2$ over $H_2 = k[t, T_0^\beta, T_0^\gamma, T_0^\delta]$ by the quiver
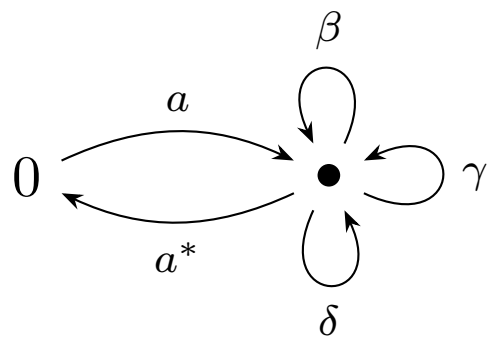
\[
\begin{tikzpicture}[baseline=(b.base), >={Stealth[length=4pt]}]
  \node (z) at (0,0) {$0$};
  \node (b) at (2,0) {$\bul$};
  \draw[->, bend left=25] (z) to node[above, font=\scriptsize] {$a$} (b);
  \draw[->, bend left=25] (b) to node[below, font=\scriptsize] {$a^*$} (z);
  \draw[->] (b) to[out=65, in=115, looseness=9] node[above, font=\scriptsize] {$\beta$} (b);
  \draw[->] (b) to[out=-25, in=25, looseness=9] node[right, font=\scriptsize] {$\gamma$} (b);
  \draw[->] (b) to[out=-115, in=-65, looseness=9] node[below, font=\scriptsize] {$\delta$} (b);
\end{tikzpicture}
\]
with relations $aa^* = \lambda_0e_0$, $\beta^2 = T_0^\beta e_\bul$, $\gamma^2 = T_0^\gamma e_\bul$, $\delta^2 = T_0^\delta e_\bul$ and $a^*a + \beta + \gamma + \delta = \frac12\lambda_0 e_\bul$; here $\lambda_0 = -t$ and $(T_0^\beta, T_0^\gamma, T_0^\delta) = \frac14(t_1^2, t_3^2, t_4^2)$ (\S\ref{ss:rootformula2}). Killing $e_0$ kills $a$ and $a^*$; writing $b = -\beta$, $c = -\gamma$ and eliminating $\delta$ gives
\[
\cC_2 = H_2\langle b, c\rangle \big/ \bigl(b^2 - T_0^\beta,\ c^2 - T_0^\gamma,\ (\tfrac12 t - b - c)^2 - T_0^\delta\bigr).
\]
Take the arc $\nu$ given by
\[
t = 4w, \qquad T_0^\beta = T_0^\gamma = w, \qquad T_0^\delta = 4w + 20w^2 .
\]
As $T_0^\delta = \frac14 t_4^2$, this gives $t_4 = \pm 4w^{1/2}(1 + 5w)^{1/2}$ on the lift, as in Figure~\ref{fig:laufer}.
Each $T_0$ has order $1$ and $t$ has finite order, so $\nu$ classifies a flop (Lemmas~\ref{lem:flopcrit}(1) and~\ref{lem:arcflop}(2), as in the proof of Lemma~\ref{lem:length2}). By Lemma~\ref{lem:basechange} its contraction algebra is $\Acon = k[[w]]\langle b, c\rangle/\bigl(b^2 - w,\ c^2 - w,\ (2w - b - c)^2 - 4w - 20w^2\bigr)$. We show that $\Acon \cong L_1$, so that $\nu$ classifies the Laufer flop.

\emph{Lemma~\ref{lem:basechange}: $\Acon \cong L_1$.} In $L_1 = k\langle x,y\rangle/(xy + yx,\ x^3 + y^2)$ put $u = x + \frac23 x^3$, $b = \frac12(u + y)$ and $c = \frac12(u - y)$. As $uy = -yu$, $b^2 = c^2 = \frac14(u^2 + y^2) =: w$, which is central, and nilpotent as $x$ and $y$ are. The third relation, $(2w - u)^2 = 4w + 20w^2$, is equivalent to $-y^2 = 4w(u + 4w)$, and it holds: by $y^2 = -x^3$ and $x^6 = 0$ (Lemma~\ref{lem:laufer}), $4w = x^2 - x^3 + \frac43 x^4$ and $4w(u + 4w) = x^3$. Hence $b, c \mapsto b, c$ defines a map $\Acon \to L_1$. It is onto, as $u = b + c$, $y = b - c$ and $x = u - \frac23 u^3 + \frac43 u^5$; both algebras have dimension $9$ (Lemma~\ref{lem:laufer}, and Corollary~\ref{cor:rootsum} as in the next paragraph), so it is an isomorphism. As the contraction algebra determines $R$ (the Donovan--Wemyss conjecture, proved in \cite{JKM,KLW}, building on \cite{HK}), $X$ is the Laufer singularity, and $Y \to X$ is its flopping contraction, up to exchanging $Y$ and the flopped variety $Y^+$. This exchange replaces $\Acon$ by its opposite algebra (\S\ref{ss:flops}). Hence $\nu$ classifies the Laufer flop.

\emph{Corollary~\ref{cor:rootsum} and Proposition~\ref{prop:GV}.} Along the lift in Figure~\ref{fig:laufer}, the pair $\varepsilon = \pm(1, 1, -1)$ vanishes, with $\varepsilon\cdot\mathbf t = \mp 10w^{3/2} + \cdots$, so its two root forms have order $1$ and the other six order $\frac12$. This is the case $r = 1 = n_2$ of the proof of Lemma~\ref{lem:length2}, so $(n_1, n_2) = (5, 1)$ and $\dim \Acon = 5 + 4 = 9$.

\emph{Theorems~B(2) and~C.} For general $\xi$, the moved arc $\nu + s\xi$ meets $D_1$ in five points and $D_2$ in one, so $\Acon$ deforms to $k^5 \times M_2(k)$. Hence $[L_1]$ is a smooth point of $\Ass_9$ and $k^5 \times M_2(k)$ is its only semisimple deformation; Section~\ref{sec:centremap} shows that the centre map is nevertheless not smooth at $[L_1]$ as a map to a moduli space of algebras.

The quadratic term of $T_0^\delta$ is needed to obtain the Laufer flop, and the arc without it shows that $h^1$, unlike $\dim Z$, is not determined by the GV invariants. With $T_0^\delta = 4w$ the root orders, and so the dimension and the GV invariants, are the same. However, the same substitution with $u = x - \frac7{12}x^3$ identifies the contraction algebra with the local block at the origin of $k\langle x,y\rangle/(xy + yx,\ y^2 + x^3 - \frac54 x^4)$. The derivations of this block form a space of dimension $7$, against $8$ for $L_1$ (straightforward to check with computer algebra), so $h^1 = 4$ against $5$. Thus that arc gives another of the flops with GV invariants $(5, 1)$ \cite[Cor~8.13]{BW}.
\end{ex}

\begin{ex}[The Laufer flops]\label{ex:laufer}
Grade $L_m$ by $\deg x = 2$ and $\deg y = 2m + 1$. As in Lemma~\ref{lem:laufer}, $x^2$ is central, $x^{2m+1}y = 0$ and $x^{4m+2} = 0$, and $L_m$ has basis $x^a$ ($0 \le a \le 4m + 1$) and $x^a y$ ($0 \le a \le 2m$). Hence $\dim L_m = 6m + 3$. The abelianisation $k[x,y]/(xy, x^{2m+1} + y^2)$ has basis $1, x, \dots, x^{2m+1}, y$, so $n_1 = 2m + 3$, and Toda's formula (Corollary~\ref{cor:toda}) gives $n_2 = m$, as in \cite[Example~1.3]{T}. Each degree of $L_m$ is at most one-dimensional, since $x^a$ has even degree and $x^a y$ odd degree, so the centre can be found monomial by monomial. As $yx^a = (-1)^a x^a y$, the power $x^a$ is central if and only if $a$ is even or $a \ge 2m + 1$. And $x^a y$ commutes with $x$ only if $x^{a+1}y = 0$, so the only central one is $x^{2m}y$. Hence
\[
\dim Z(L_m) = (2m+1) + (m+1) + 1 = 3m + 3 = n_1 + n_2 ,
\]
in line with Theorem~B(1). Both dimensions agree with \cite[Prop~3.3]{JS}, as $L_m$ is the Jacobi algebra of $\mathsf{D}_{m+1,\infty}$ \cite[Table~1]{BW}. By Theorem~C, each $[L_m]$ is a smooth point of $\Ass_{6m+3}$, and $k^{2m+3} \times M_2(k)^m$ is the only semisimple algebra to which $L_m$ deforms. The proof of Lemma~\ref{lem:length2} recovers these invariants from the root formula, as the case $n_2 = m$ and $r \ge m$; by Brown and Wemyss's classification \cite[Cor~8.13]{BW} several flops share them (\cite{BW18} gave the first two non-isomorphic flops with the same GV invariants).
\end{ex}

\appendix

\section{Corrections to \texorpdfstring{\cite{K}}{[K]}}\label{app:corrections}

In this appendix we correct some slips in \cite{K} in the parts used here, and one in \cite[Remark~4.14]{K}. None of them affects the main results of \cite{K}. Misprints in the preprojective relations of \cite[Examples~4.6, 4.7, 4.18--4.20]{K} do not affect what we use and are not listed. Here we use the notation of \cite{K}, including its vertex labels, its form $t$, which is minus ours (\S\ref{ss:lie}), and its arc parameter $T$, which is $w$ in \S\ref{ss:typeE}.

\paragraph{Statements.}
\begin{itemize}
\item Prop~1.7 (in the Introduction): the statement that \emph{every} deformation of $\Acon$ with a semisimple fibre has fibre $\bigoplus_i \mathrm{Mat}_{i \times i}(\mathbb{C})^{\oplus n_i}$ is attributed to \cite[Thm~4.2]{HT}, which treats only the deformation induced by a generic deformation of the flop. The general statement is Theorem~C of the present paper.
\item Prop~4.21: the conclusion should read $A_l \otimes_{H_l} H_\Gamma \cong e\cA_\Gamma e$, which is what the proof shows; ``$e\cA_\Gamma e \otimes_{H_\Gamma} H_l$'' is not defined, as $H_l \subset H_\Gamma$.
\item Remark~4.14: the dimensions of the algebras $\Delta_{\mathrm{con}}$ are $1, 4, 12, 24, 40, 60$, as in \cite[Remark~3.17]{DW} ($9$ is printed in length $3$).
\end{itemize}

\paragraph{The explicit presentations.} Table~\ref{tab:Kmisprints} lists the misprints in the presentations of \cite[Thm~1.3, Examples~4.19--4.27 and~5.3--5.4]{K}, grouped by type, with the printed form in brackets where it helps.

\begin{table}[ht]
\caption{Misprints in the explicit presentations of \cite{K}.}\label{tab:Kmisprints}
\small
\renewcommand{\arraystretch}{1.15}
\begin{tabular}{@{}>{\raggedright\arraybackslash}p{2.4cm} >{\raggedright\arraybackslash}p{\dimexpr\textwidth-2.4cm-2\tabcolsep\relax}@{}}
\toprule
Place in \cite{K} & Should read \\
\midrule
\multicolumn{2}{@{}l}{\emph{Invariants and constraints}} \\
Ex.~4.19, 4.25 & $t_0 + 2t_1 + 3t_2 + 2t_3 + t_4 + 2t_5 + 3t_6 + 4t_7 = 0$ (both examples use the labelling of Ex.~4.25, with centre $7$; printed $t_0 + 2t_1 + 4t_2 + 3t_3 + 2t_4 + 3t_5 + 2t_6 + t_7$) \\
Ex.~4.24 & $t = (2t_0 + t_1)/2$ (printed first as $2t_0 + t_1$) \\
Ex.~4.26 & $t = (4t_0 + 3t_1 + 2t_2 + t_3)/4$ (printed $/5$) in the intermediate step \\ \addlinespace[4pt]
\multicolumn{2}{@{}l}{\emph{Changes of variables and parameters}} \\
Ex.~4.23 & $\gamma = c - t_2/2$, $\delta = d - t_3/2$ (printed $\gamma = c - t_1/2$, $\delta = d - t_1/2$) \\
Ex.~4.26 & $B = b + \tau_5$, $C = c + \tau_3$ (printed $\tau_3$, $\tau_5$); $T_3, T_2, T_1, T_0 = +\sigma_2, +\sigma_3, +\sigma_4, +\sigma_5$ of $\tau_1, \dots, \tau_5$ (printed $T_3 := -\sigma_2, \dots, T_0 := -\sigma_5$), as \cite{K}'s identities $T_3 = R_2 + S_2 - R_1^2$, \dots, $T_0 = R_3S_2$ require \\
Ex.~4.27 & $T_0^\beta = -\sigma_2(\tau^\beta)$, $T_1^\gamma = -\sigma_2(\tau^\gamma)$, $T_0^\gamma = -\sigma_3(\tau^\gamma)$ (printed $-\sigma_2(\tau^\beta)/4$, $-\sigma_3(\tau^\gamma)$, $-\sigma_2(\tau^\gamma)/27$), following lengths $3$ and $4$ \\ \addlinespace[4pt]
\multicolumn{2}{@{}l}{\emph{Relations}} \\
Thm~1.3(3), Ex.~4.24 & $e_6$ for $e_4$ throughout (vertex labels of \cite[Ex.~4.24]{K}): the relations $a^*a = \delta^2 - T_0^\delta$, $\beta^3 = T_1^\beta\beta + T_0^\beta$, $\gamma^3 = T_1^\gamma\gamma + T_0^\gamma$, $\beta + \gamma + \delta = t/3$ hold at vertex $6$ \\
Ex.~4.24 & $da^*$ for $\delta a^*$ in the $H_\Gamma$ presentation \\
Ex.~4.25 & $ad = (t_0 + (t_1 + t_2))a$ in the $H_\Gamma$ presentation (printed without the final $a$) \\
Ex.~4.27 & $t_8e_8$ for $t_8e_5$ \\
Ex.~5.3 & second relation $(\beta + \gamma)^2 - \gamma^3$ (printed $(\beta + \gamma)^3 - \gamma^2$), as $a^*a = (\beta + \gamma)^2 - \beta^3 = (\beta + \gamma)^2 - \gamma^3$ just above it requires \\
Ex.~5.4 & $aa^* = -Te_0$ (printed $Te_0$), as in Ex.~5.3 and~5.5: by \cite[Thm~1.3(4), (6)]{K}, $aa^* = (t^3 - T_1^\delta t - T_0^\delta)e_0$ in length $4$, which the arc sends to $-Te_0$, and similarly in length $6$; also $a(a^*a) = a\delta^m - Ta$, with $m$ the size of the arm of $\delta$ and $a\delta = 0$, forces the sign \\
\bottomrule
\end{tabular}
\end{table}

\paragraph{Example~5.5.} The displayed relations $\gamma\beta\gamma + \gamma^2\beta + \gamma\beta^3 = Te_4$ and $(\gamma + \beta^2)^2 + \beta\gamma\beta = 0$, and the presentation derived from them, are those of the classifying arc with $T_1 = 0$ and $T_0 = -T$. For the stated arc, with $T_1 = T_0 = T$, \cite[Thm~1.3(5)]{K} gives $\gamma\beta\gamma + \gamma^2\beta + \gamma\beta^3 = (\gamma + \beta^2)^2 + \beta\gamma\beta = -Te_4$. Both arcs classify flops with GV invariants $(8, 6, 4, 2, 1)$ (Lemmas~\ref{lem:flopcrit}(1) and~\ref{lem:arcflop}(2) and the proof of Proposition~\ref{prop:typeE}).

\section*{Acknowledgements}

We thank Michael Wemyss for helpful correspondence. This work was supported by the Beijing Natural Science Foundation (Grant Nos.~1262017, 1252011).

\section*{Use of AI}

The ideas in this paper, and the directions the work took, are the authors'. They include the overall approach, pulling the centre back from the root combinatorics of deformed preprojective algebras along the universal flops and building on \cite{K}, and the key steps of the proofs. We used AI language models (principally Anthropic's Claude) extensively, under our close review and guidance, to write the paper and to help fill in its technical details. The AI drafted the \LaTeX{} from our outline and wrote most of the prose through many rounds of revision. It filled in details of the proofs from our sketches and proposed simplifications of them, and it handled much of the Lie-theoretic bookkeeping. It also carried out repeated referee-style reviews of the drafts, which found gaps in our arguments and in our use of the literature and proposed repairs, and it checked our citations and bibliographic data against the sources.

We also made wide use of it to write and run computer-algebra calculations. These helped us to experiment, to draw the results together and to validate them. No theorem, proposition or lemma relies on them. Only Remarks~\ref{rem:noER} and~\ref{rem:reduced}, the comparison of derivations in Example~\ref{ex:lauferarc} and an aside in the proof of Lemma~\ref{lem:hilbTangent} rest on such calculations, as stated there. The authors chose what to prove and how to frame it, directed every revision, decided on every substantive change and checked the arguments. We take full responsibility for the content.

\end{document}